\magnification=1100
\overfullrule0pt

\input amssym.def
\input prepictex
\input pictex
\input postpictex



\def\CC{{\Bbb C}}

\def\RR{{\Bbb R}}
\def\ZZ{{\Bbb Z}}

\def\fb{\frak{b}}
\def\fg{\frak{g}}
\def\fh{\frak{h}}
\def\fn{\frak{n}}

\def\End{\hbox{End}}
\def\Hom{\hbox{Hom}}
\def\Ind{\hbox{Ind}}
\def\Res{\hbox{Res}}

\def\id{\hbox{id}}
\def\mt{\hbox{mt}}
\def\rad{\hbox{rad}}
\def\tr{\hbox{tr}}
\def\Tr{\hbox{Tr}}
\def\wt{\hbox{wt}}


\font\smallcaps=cmcsc10
\font\titlefont=cmr10 scaled \magstep1

\font\sectionfont=cmbx10
\font\tinyrm=cmr10 at 8pt


\newcount\sectno
\newcount\subsectno
\newcount\resultno

\def\section #1. #2\par{
\sectno=#1
\resultno=0
\bigskip\noindent{\sectionfont #1.  #2}~\medbreak}

\def\subsection #1\par{\bigskip\noindent{\it  #1} \medbreak}


\def\prop{ \global\advance\resultno by 1
\bigskip\noindent{\bf Proposition \the\sectno.\the\resultno. }\sl}
\def\lemma{ \global\advance\resultno by 1
\bigskip\noindent{\bf Lemma \the\sectno.\the\resultno. }
\sl}
\def\remark{ \global\advance\resultno by 1
\bigskip\noindent{\bf Remark \the\sectno.\the\resultno. }}
\def\example{ \global\advance\resultno by 1
\bigskip\noindent{\bf Example \the\sectno.\the\resultno. }\sl}
\def\cor{ \global\advance\resultno by 1
\bigskip\noindent{\bf Corollary \the\sectno.\the\resultno. }\sl}
\def\thm{ \global\advance\resultno by 1
\bigskip\noindent{\bf Theorem \the\sectno.\the\resultno. }\sl}
\def\defn{ \global\advance\resultno by 1
\bigskip\noindent{\it Definition \the\sectno.\the\resultno. }\slrm}
\def\endthm{\rm\bigskip}

\def\pf{\rm\medskip\noindent{\it Proof. }}
\def\endpf{\qed\hfil\bigskip}


\def\qed{\hbox{\hskip 1pt\vrule width4pt height 6pt depth1.5pt \hskip 1pt}}

\def\sqr#1#2{{\vcenter{\vbox{\hrule height.#2pt
\hbox{\vrule width.#2pt height#1pt \kern#1pt
\vrule width.2pt}
\hrule height.2pt}}}}
\def\square{\mathchoice\sqr54\sqr54\sqr{3.5}3\sqr{2.5}3}


\def\formula{\global\advance\resultno by 1
\eqno{(\the\sectno.\the\resultno)}}
\def\formulano{\global\advance\resultno by 1 (\the\sectno.\the\resultno)}
\def\tableno{\global\advance\resultno by 1
\the\sectno.\the\resultno. }
\def\lformula{\global\advance\resultno by 1
\leqno(\the\sectno.\the\resultno)}


\def\monthname {\ifcase\month\or January\or February\or March\or April\or
May\or June\or
July\or August\or September\or October\or November\or December\fi}

\newcount\mins  \newcount\hours  \hours=\time \mins=\time
\def\now{\divide\hours by60 \multiply\hours by60 \advance\mins by-\hours
     \divide\hours by60         
     \ifnum\hours>12 \advance\hours by-12
       \number\hours:\ifnum\mins<10 0\fi\number\mins\ P.M.\else
       \number\hours:\ifnum\mins<10 0\fi\number\mins\ A.M.\fi}


\def\mapright#1{\smash{\mathop
        {\longrightarrow}\limits^{#1}}}


\nopagenumbers
\def\runningtitle{\smallcaps braids and jantzen filtrations}
\headline={\ifnum\pageno>1\eoheadline\else\firstheadline\fi}
\def\names{\smallcaps rosa orellana and arun ram}
\def\firstheadline{}
\def\eoheadline{\ifodd\pageno\oddheadline\else\evenheadline\fi}
\def\oddheadline{\tenrm\hfil\runningtitle\hfil\folio}
\def\evenheadline{\tenrm \folio\hfil{\names}\hfil}

\vphantom{$ $}  
\vskip.75truein

\centerline{\titlefont Affine braids, Markov traces and the
category ${\cal O}$}
\bigskip
\centerline{\rm Rosa Orellana}
\centerline{Department of Mathematics}
\centerline{Dartmouth College}
\centerline{Hanover, NH 03755-3551}
\centerline{{\tt Rosa.C.Orellana@dartmouth.edu}}
\medskip
\centerline{\rm and}
\medskip
\centerline{\rm Arun Ram${}^\ast$}
\centerline{Department of Mathematics}
\centerline{University of Wisconsin, Madison}
\centerline{Madison, WI 53706}
\centerline{{\tt ram@math.wisc.edu}}
\smallskip
\footnote{}{\tinyrm \noindent
${}^\ast$ Research supported in part by National
Science Foundation grant DMS-9971099, the National Security Agency
and EPSRC grant GR K99015.}

\bigskip



\section 1. Introduction

This paper provides a unified approach to results on representations
of affine Hecke algebras, cyclotomic Hecke algebras, 
affine BMW algebras, cyclotomic BMW algebras,
Markov traces, Jacobi-Trudi type identities, 
dual pairs [Ze], and link invariants [Tu2].
The key observation in the genesis of this paper was that
the technical tools used to obtain the
results in Orellana [Or] and Suzuki [Su], two 
a priori unrelated papers,
are really the same.  Here we develop this method and explain
how to apply it to obtain results similar to those in [Or]
and [Su] in more general settings.
Some specific new results which are obtained are the following:
\smallskip\noindent
\itemitem{(a)} A generalization of the results on Markov traces
obtained by Orellana [Or] to centralizer algebras coming from 
quantum groups of all Lie types.
\smallskip\noindent
\itemitem{(b)} A generalization of the results of Suzuki [Su]
to show that Kazhdan-Lusztig polynomials of all finite Weyl groups
occur as decomposition numbers in the representation theory of affine
braid groups of type A,
\smallskip\noindent
\itemitem{(c)} A generalization of the functors used by
Zelevinsky [Ze] to representations of affine braid groups
of type A,
\smallskip\noindent
\itemitem{(d)} We define the affine BMW algebra
(Birman-Murakami-Wenzl) and show that it has a representation
theory analogous to that of affine Hecke algebras.  In particular
there are ``standard modules'' for these algebras which
have composition series where multiplicites of
the factors are given by Kazhdan-Lusztig polynomials for 
Weyl groups of types A,B,and C.
\smallskip\noindent
\itemitem{(e)} We generalize the results of Leduc and Ram [LR]
to affine centralizer algebras. 
\smallskip\noindent

Let $U_h\fg$ be the Drinfel'd-Jimbo quantum group associated to
a finite dimensional complex semisimple Lie algebra $\fg$.
If $M$ is a (possibly infinite dimensional) $U_h\fg$-module
in the category ${\cal O}$ and $V$ is a finite dimensional
$U_h\fg$-module then we show that the affine
braid group $\tilde {\cal B}_k$ acts on the $U_h\fg$-module
$M\otimes V^{\otimes k}$.  Fix $V$ and define 
$$F_\lambda(M)
=\left(\matrix{\hbox{the vector space
of highest weight vectors} \cr
\hbox{of weight $\lambda$ in 
$M\otimes V^{\otimes k}$.} \cr}
\right).
$$
Then $F_\lambda$ is a functor from
$U_h\fg$ modules in category ${\cal O}$ to
finite dimensional modules for the affine braid group
$\tilde {\cal B}_k$ which takes
\smallskip\noindent
\itemitem{(1)} finite dimensional $U_h\fg$ modules to 
``calibrated'' $\tilde {\cal B}_k$ modules,
\smallskip\noindent
\itemitem{(2)} Verma modules to ``standard'' modules, and
\smallskip\noindent
\itemitem{(3)} under appropriate conditions, irreducible
$U_h\fg$ modules to irreducible $\tilde {\cal B}_k$ modules.
\smallskip\noindent
Applying the functor $F_\lambda$ to a Jantzen filtration of Verma modules
of $U_h\fg$ provides a ``Jantzen filtration'' of the standard
modules of $\tilde {\cal B}_k$ and shows that the irreducible
$\tilde {\cal B}_k$ modules appear in a composition series
of the standard module with multiplicities given by the 
Kazhdan-Lusztig polynomials of the Weyl group of $\fg$.  
Though $\tilde {\cal B}_k$ is always the affine braid group of
type A, the Weyl group of $\fg$ is {\it not necessarily} of type A.

Applying the functor $F_\lambda$ to the BGG resolution of
an irreducible highest weight module provides a BGG resolution
for the corresponding $\tilde {\cal B}_k$-modules and a corresponding
``Jacobi-Trudi'' identities for the characters of 
$\tilde {\cal B}_k$ modules.  Once again, it is interesting to
note that, though $\tilde {\cal B}_k$ is the affine braid group of
type A, it is the Weyl group of a different type which appears
in this Jacobi-Trudi identity.

Using the general formulation for constructing Markov traces on 
braid groups, given for example in [Tu1], we obtain a Markov
trace on the affine braid group $\tilde {\cal B}_k$ for every
choice of $\fg$ and $U_h\fg$ modules $M$ and $V$.
\itemitem{(a)}  If $\fg={\frak{sl}}_{n+1}$, $M=L(0)$ and
$V=L(\omega_1)$ this gives the Markov trace
on the Hecke algebra studied in [Jo1] and [Wz].
\itemitem{(b)}  If $\fg={\frak{sl}}_{2}$, $M=L(0)$ and
$V=L(\omega_1)$ this gives the Markov trace
on the Temperley-Lieb algebra used by Jones [Jo2].
\itemitem{(c)}  If $\fg={\frak{sl}}_{n+1}$, $M=L(k\omega_\ell)$
with $k$ and $\ell$ large and $n$ very large,
and $V=L(\omega_1)$ this gives the 
Markov traces on the Hecke algebra of type B studied by
[GL], [Lb], [Ic] and [Or].
\itemitem{(d)}  If $\fg={\frak{sl}}_{n+1}$, $M=L(\lambda)$,
where $\lambda$ is ``large'', and $V=L(\omega_1)$ this gives the 
Markov traces on the cyclotomic Hecke algebras
introduced by Lambropoulou [Lb] and studied in [GIM].
\itemitem{(e)}  If $\fg={\frak{so}}_{n}$ or 
$\fg={\frak{sp}}_{2n}$, $M=L(0)$ and
$V=L(\omega_1)$ this gives the Markov traces
used to construct Kauffman polynomials.
\smallskip\noindent
For general $\fg$, general $V$, and $M=L(0)$, this mechanism
gives the traces necessary to compute the Reshetikhin-Turaev link
invariants [RT].  In some sense, this paper is a study of
the representation theory behind the generalization
of the Reshetikhin-Turaev method given in [Tu2].

In the final section of this paper we describe precisely
the combinatorics of the representations $F_\lambda(M)$
in the cases when $\fg$ is type $A_n, B_n, C_n$ or $D_n$
and $V$ is the fundamental representation.  In these
cases the representations can be constructed with partitions,
standard tableaux, up-down tableaux, multisegments and the combinatorics
of Young diagrams.  In particular, in type A,
the functor $F_\lambda$ naturally constructs the standard modules
and irreducible modules of affine Hecke algebras of type A
in terms of multisegments (a classification originally
obtained by Zelevinsky [Ze2] by different methods).  We then specify
explicitly the correspondence between the decomposition numbers
of the affine Hecke algebra and Kazhdan-Lusztig polynomials for
the symmetric group.
Using the recent results of Polo [Po] we can show that 
every polynomial in $1+{\tt v}\ZZ_{\ge 0}[{\tt v}]$
is a decomposition number for the affine Hecke algebra.

\medskip\noindent
{\bf Acknowledgements.}  A. Ram thanks P. Littelmann,
the Department of Mathematics at the University of Strasbourg
and the Isaac Newton Institute for the Mathematical Sciences
at Cambridge University for hospitality and support
during residencies when this paper was written.

\section 2. Preliminaries on quantum groups

Let $U_h\fg$ be the Drinfel'd-Jimbo quantum group corresponding
to a finite dimensional complex semisimple Lie algebra $\fg$.  
Let us fix some notations.  In particular, fix a triangular
decomposition
$$\fg = \fn^-\oplus\fh\oplus\fn^+,
\qquad \fn^+ = \bigoplus_{\alpha>0} \fg_\alpha,
\qquad \fb^+ = \fh \oplus \fn^+,
$$
and let $W$ be the Weyl group of $\fg$.
Let $\langle,\rangle$ be the usual inner product on $\fh^*$
so that, if $\alpha$ is a root, the corresponding reflection 
$s_\alpha$ in $W$ is given by 
$$s_\alpha\lambda = \lambda-\langle\lambda,\alpha^\vee\rangle \alpha,
\qquad\hbox{where}\quad \alpha^\vee={2\alpha\over \langle\alpha,\alpha\rangle}.
\qquad\hbox{The element}\quad
\rho = {1\over 2}\sum_{\alpha>0} \alpha$$
is often viewed as an element of $\fh$ by using the form $\langle,\rangle$ to identify
$\fh$ and $\fh^*$.  We shall use the conventions for quantum groups
as in [Dr] and [LR] so that 
$$q=e^{h/2},
\qquad \fh\subseteq U_h\fg,
\qquad\hbox{and} \quad
U_h\fg\cong U\fg[[h]],
\hbox{\ \ as algebras.}
$$
The quantum group has a triangular decomposition corresponding to that
of $\fg$,
$$U_h\fg = U_h\fn^- \otimes U_h\fh \otimes U_h\fn^+
\qquad\hbox{and}\qquad
U_h\fb^+ = U_h\fh\otimes U_h\fn^+.
$$

\subsection The category ${\cal O}$

If $M$ is a $U_h\fg$ module and $\lambda\in \fh^*$ the {$\lambda$
weight space} of $M$ is 
$$M_\lambda = \{m\in M \ |\ am=\lambda(a)m,
\ \hbox{for all $a\in \fh$} \}.
$$
The category ${\cal O}$ is the category of $U_h\fg$ modules $M$
such that
\smallskip\noindent
\itemitem{(a)} $M = \bigoplus_{\lambda\in\fh^*} M_\lambda$,
\smallskip\noindent
\itemitem{(b)} For all $m\in M$, $\dim(U_h\fn^+ m)$ is finite,
\smallskip\noindent
\itemitem{(c)} $M$ is finitely generated as a $U_h\fg$ module.
\smallskip\noindent
For $\mu\in \fh^*$ let
\smallskip\noindent
\itemitem{} $M(\mu)$ be the Verma module of highest weight $\mu$, and let
\smallskip\noindent
\itemitem{} $L(\mu)$ be the irreducible module of highest weight $\mu$.
\smallskip\noindent
The irreducible module $L(\mu)$ is the quotient of $M(\mu)$ by a maximal
proper submodule and $M(\mu)=U_h\fg\otimes_{U_h\fb^+} \CC v_\mu^+$ where
$\CC v_\mu^+$ is the one dimensional $U_h\fb^+$ module spanned by a vector
$v_\mu^+$ such that $av^+_\mu=\mu(a)v^+_\mu$ for $a\in\fh$ and $U_h\fn^+ v^+_\mu=0$.
Every module $M\in {\cal O}$ has a finite composition series with
factors $L(\mu)$, $\mu\in\fh^*$.  Each of the sets
$$
\{ [L(\lambda)] \ |\ \lambda\in \fh^*\}
\qquad\hbox{and}\qquad
\{ [M(\lambda)] \ |\ \lambda\in \fh^*\}
$$
(where $[M]$ denotes the isomorphism class of the module $M$)
are bases of the Grothendieck group of the category ${\cal O}$.

If $M$ is a $U_h\fg$ module generated by a highest weight vector of weight 
$\lambda$ (i.e., a vector $v^+$ such that $av^+=\lambda(a)v^+$
for $a\in \fh$ and $U_h\fn^+=0$)
then any element of the center $Z(U_h\fg)$ acts on $M$
by a constant,
$$
zm=\chi^\lambda(z)m,
\quad\hbox{for $z\in Z(U_h\fg)$, $m\in M$, $\chi^\lambda(z)\in \CC$.}
$$
For each $U_h\fg$ module $M\in {\cal O}$ let
$$M^{[\lambda]} = \bigoplus_{\nu\in Q} M^{[\lambda]}_{\lambda+\nu},
\qquad\hbox{where}\quad
Q=\sum_{i=1}^n \ZZ\alpha_i,$$
$\alpha_1,\ldots, \alpha_n$ are the simple roots and
$$M^{[\lambda]}_{\lambda+\nu} = \{
m\in M_{\lambda+\nu} \ |\ \hbox{there is $k\in \ZZ_{>0}$ such that
$(z-\chi^\lambda(z))^km=0$ for all $z\in Z(U_h\fg)$} \}.$$
Then
$$M=\bigoplus_\lambda M^{[\lambda]},$$
where the sum is over all {\it integrally dominant} weights
$\lambda\in \fh^*$ i.e., $\lambda\in \fh^*$ such that 
$\langle \lambda+\rho,\alpha^\vee\rangle\not\in \ZZ_{<0}$
for all $\alpha\in R^+$.

The dot action of the Weyl group $W$ on $\fh^*$ is given by
$$w\circ \lambda = w(\lambda+\rho)-\rho,
\qquad w\in W, \lambda\in \fh^*.
$$
For a fixed $\lambda\in \fh^*$ the
stabilizer of the dot action of the integral Weyl group
$$W^\lambda = < s_\alpha\ |\ \langle\lambda+\rho,\alpha^\vee\rangle
\in \ZZ>
\formula$$
is the subgroup
$$W_{\lambda+\rho} = \{ w\in W\ |\ w(\lambda+\rho)=\lambda+\rho\}$$
and the elements of $W^\lambda\circ\lambda$ are exactly the 
$w\circ\lambda$ such that $w\in W^\lambda$ is the longest
element of the coset $wW_{\lambda+\rho}$ in $W^\lambda$.
To summarize, there is a decomposition of the category ${\cal O}$,
$${\cal O} = \bigoplus_\lambda {\cal O}^{[\lambda]},
\formula$$
where the sum is over all integrally dominant weights $\lambda\in \fh^*$
and ${\cal O}^{[\lambda]}$ is the full subcategory of modules
$M\in {\cal O}$ such that $M=M^{[\lambda]}$.  The Grothendieck group
of the category ${\cal O}^{[\lambda]}$ has bases
$$
\{ [L(\mu)]\ |\ \mu\in W^\lambda\circ\lambda\}
\qquad\hbox{and}\qquad
\{ [M(\mu)]\ |\ \mu\in W^\lambda\circ\lambda\}.
\formula$$

\vfill\eject

\subsection Jantzen filtrations

Following the notations for the quantum group used in
[LR, \S 2], let $\fh$, $X_1,\ldots, X_r$ and $Y_1,\ldots, Y_r$
be the standard generators of the quantum group $U_h\fg$ which 
satisfy the quantum Serre relations.
The {\it Cartan involution} $\theta\colon U_h\fg\to U_h\fg$
is the algebra anti-involution defined by
$$
\theta(X_i)=Y_i,
\qquad
\theta(Y_i)=X_i,
\qquad\hbox{and}\qquad
\theta(a)=a,
\quad\hbox{for $a\in \fh$}.
\formula$$
A {\it contravariant form} on a $U_h\fg$ module $M$ is a
symmetric bilinear form $\langle,\rangle\colon M\times M\to \CC$
such that
$$\langle um_1,m_2\rangle = \langle m_1, \theta(u)m_2\rangle,
u\in U_h\fg, \ \ m_1,m_2\in M.
$$

Fix $\lambda\in \fh^*$ and $\delta\in \fh^*$ such that $\lambda+t\delta$
is integrally dominant for all small positive real numbers $t$.
Consider $t$ as an indeterminate and consider the Verma module
$$M(\lambda+t\delta)
=U_h\fg[t]\otimes_{U_h\fb^+[t]} \CC_{\lambda+t\delta}
$$
as the module for $U_h\fg[t]=\CC[t]\otimes_\CC U_h\fg$ generated by a vector
$v^+$ such that $av^+=(\lambda+t\delta)(a)v^+$ for $a\in \fh$ and
$U_h\fn^+[t]v^+=0$.  There is a unique contravariant form
$\langle,\rangle_t\colon M(\lambda+t\delta)\times M(\lambda+t\delta)
\to \CC[t]$ such that $\langle v^+,v^+\rangle_t = 1$.  Define
$$
M(\lambda+t\delta)(j) = \{ m\in M(\lambda+t\delta)\ |\ 
\langle m,n\rangle_t\in t^jM(\lambda+t\delta) \hbox{\ for all\ }
n\in M(\lambda+t\delta) \}.
$$ 
The ``specialization of $M(\lambda+t\delta)(j)$ at $t=0$'' is
$$M(\lambda)^{(j)} = 
\hbox{image of $M(\lambda+t\delta)(j)$ in $M(\lambda+t\delta)\otimes_{\CC[t]}
\CC[t]/t\CC[t]$}
$$
and the
{\it Jantzen filtration} of $M(\lambda)$ is
$$M(\lambda) = M(\lambda)^{(0)}\supseteq M(\lambda)^{(1)}
\supseteq\cdots\, . \formula
$$
By [Jz, Theorem 5.3],
the Jantzen filtration is a filtration of $M(\lambda)$ by
$U_h\fg$ modules, the module $M(\lambda)^{(1)}$ is a maximal proper
submodule of $M(\lambda)$ and each quotient 
$M(\lambda)^{(i)}/M(\lambda)^{(i+1)}$ has a nondegenerate contravariant form.
It is known [Bb] that the Jantzen filtration does not depend on the
choice of $\delta$.
It is a deep theorem [BB] that the quotients 
$M(\lambda)^{(i)}/M(\lambda)^{(i+1)}$ 
are semisimple and that if $w\in W^\mu$ and $y\in W^\mu$ are maximal length 
in their cosets $wW_{\mu+\rho}$ and $yW_{\mu+\rho}$, respectively,
then the Kazhdan-Lusztig polynomial for $W^\mu$ is
$$\sum_{j\ge 0} [M(w\circ\mu)^{(j)}/M(w\circ\mu)^{(j+1)}:L(y\circ \mu)]
{\tt v}^{{1\over2}(\ell(y)-\ell(w)-j)}
=P_{wy}({\tt v}), \formula
$$
where $\ell$ is the length function on $W^\mu$ and
$[M(w\circ\mu)^{(j)}/M(w\circ\mu)^{(j+1)}:L(y\circ \mu)]$
is the multiplicity of the simple module $L(y\circ\mu)$ 
in the $j$th factor of the Jantzen filtration of $M(\lambda)$.

\vfill\eject

\subsection The BGG resolution

Not all simple modules $L(\lambda)$ in the category ${\cal O}$
have a BGG resolution.  The general form of the BGG resolution
given by Gabber and Joseph [GJ] is as follows.

Let $\mu\in \fh^*$ be such that $-(\mu+\rho)$ is dominant and regular
and let $W^\nu_J$ be a parabolic subgroup of the integral Weyl group
$W^\mu$.  Let $w_0$ be the longest element of $W_J^\mu$ and
fix $\nu = w_0\circ\mu$.   Define a resolution 
$$
0\longrightarrow C_{\ell(w_0)} \longrightarrow \cdots \longrightarrow 
C_2\ \mapright{d_2}\  C_1\ \mapright{d_2}\  C_0\longrightarrow
L(\nu)\longrightarrow 0
\formula
$$
of the simple module $L(\nu)$ by Verma modules by setting
$$C_j = \bigoplus_{\ell(w)=j} M(w\circ\nu),$$
where the sum is over all $w\in W^\mu_J$ of length $j$, and
defining the map
$$d_j\colon C_j\to C_{j-1},
\quad\hbox{by the matrix}\quad
(d_j)_{v,w} = \cases{
\varepsilon_{v,w}\iota_{v,w}, &if $v\to w$, \cr
0, &otherwise, \cr},$$
where $v\to w$ means that there is a (not necessarily
simple) root $\alpha$ such that $w=s_\alpha v$ and
$\ell(w)=\ell(v)-1$,  the maps $\iota_{v,w}$ are fixed choices
of inclusions
$$\iota_{v,w}\colon M(v\circ\nu)\hookrightarrow
M(w\circ\nu),
\qquad\hbox{and}\qquad
\varepsilon_{v,w}=\pm1,$$
are fixed choices of signs such that
$$\varepsilon_{u,v}\varepsilon_{v,w}
=-\varepsilon_{u,v'}\varepsilon_{v',w}
\qquad\hbox{if}\quad
u\to v\to w,\quad u\to v'\to w\quad\hbox{and}\quad
v\ne v'.$$
Gabber and Joseph [GJ] prove that the sequence (2.7)
is exact in this general setting.  See [BGG] and [Dx,7.8.14]
for the original form of the BGG resolution.
From the exactness of (2.7) it follows that if $-(\mu+\rho)$ is dominant and regular
then, in the Grothendieck group of the category ${\cal O}$,
$$
[L(\nu)] = \sum_{w\in W_J^\mu} (-1)^{\ell(w)}[M(w\circ\nu)],
\formula$$
where $\nu=w_0\circ \mu$ and $w_0$ is the longest element of
$W^\mu_J$.

\subsection $\check R_{MN}$ matrices and the quantum Casimir $C_M$

Let $U_h\fg$ be the Drinfeld-Jimbo quantum group corresponding
to a finite dimensional complex semisimple Lie algebra $\fg$.  There
is an invertible element ${\cal R}=\sum a_i\otimes b_i$ in (a suitable
completion of) $U_h\fg\otimes U_h\fg$ such that, for any two $U_h\fg$
modules $M$ and $N$, the map
$$
\beginpicture
\setcoordinatesystem units <1cm,.5cm>         
\setplotarea x from -8 to 2, y from 0 to 1.5    
\put{$
\matrix{
\check R_{MN}\colon &M\otimes N &\longrightarrow &N\otimes M\cr
&m\otimes n &\longmapsto &\displaystyle{
\sum b_in\otimes a_i m }\cr
}$} at -5 1
\put{$M\otimes N$} at 0.5 2.4
\put{$N\otimes M$} at 0.5 0.1
\put{$\bullet$} at  0.1 1.9      
\put{$\bullet$} at  0.9 1.9      %
\put{$\bullet$} at  0.1 .6          
\put{$\bullet$} at  0.9 .6          %
\setquadratic
\plot  0.1 .6  .15 .9  .4 1.15 /
\plot  .6 1.35  .85 1.6  0.9 1.9 /
\plot 0.1 1.9  .15 1.6  .5 1.25  .85 .9  0.9 .6 /
\endpicture
$$
is a $U_h\fg$ module isomorphism.  There is also a
{\it quantum Casimir element} $e^{-h\rho}u$ in the center of $U_h\fg$
and, for a $U_h\fg$ module $M$ we define
$$
\beginpicture
\setcoordinatesystem units <1cm,.6cm>         
\setplotarea x from -6 to 1, y from 0 to 1.5    
\put{$
\matrix{
C_M\colon &M &\longrightarrow &M \cr
&m &\longmapsto &(e^{-h\rho}u)m \cr
}$} at -5 1.2
\put{$M$} at 0.7 2.2
\put{$M$} at 0.7 0.3
\put{$C_M$} at 1.2 1.3
\put{$\bullet$} at  0.7 1.8      
\put{$\bullet$} at  0.7 .7          
\plot  0.7 1.8   0.7 .7 /
\endpicture
$$
The elements ${\cal R}$ and $e^{-h\rho}u$ satisfy relations
(see [LR, (2.1-2.12)]) which imply that,
for $U_h\fg$ modules $M,N,P$ and a $U_h\fg$ module isomorphism
$\tau_M\colon M\to M$,
$$\eqalign{
\beginpicture
\setcoordinatesystem units <1cm,.5cm>         
\setplotarea x from 0 to 2, y from -1 to 2    
\put{$M\otimes N$} at 0.5 2.4
\put{$N\otimes M$} at 0.5 -1.2
\put{$\tau_M$} at 1.2 -0.1
\put{$\bullet$} at  0.1 1.9      
\put{$\bullet$} at  0.9 1.9      %
\put{$\bullet$} at  0.1 -0.7          
\put{$\bullet$} at  0.9 -0.7          %
\plot  0.9 .6   0.9 -0.7 /
\plot  0.1 .6   0.1 -0.7 /
\setquadratic
\plot  0.1 .6  .15 .9  .4 1.15 /
\plot  .6 1.35  .85 1.6  0.9 1.9 /
\plot 0.1 1.9  .15 1.6  .5 1.25  .85 .9  0.9 .6 /
\endpicture
~~&=~~
\beginpicture
\setcoordinatesystem units <1cm,.5cm>         
\setplotarea x from 0 to 2, y from 0 to 2    
\put{$M\otimes N$} at 0.5 2.4
\put{$N\otimes M$} at 0.5 -1.2
\put{$\tau_M$} at -0.2 1.3
\put{$\bullet$} at  0.1 1.9      
\put{$\bullet$} at  0.9 1.9      %
\put{$\bullet$} at  0.1 -0.7          
\put{$\bullet$} at  0.9 -0.7          %
\plot  0.9 .6   0.9 1.9 /
\plot  0.1 .6   0.1 1.9 /
\setquadratic
\plot  0.1 -0.7  .15 -0.4  .4 -0.15 /
\plot  .6 0.05  .85 .3  0.9 .6 /
\plot 0.1 .6  .15 0.3  .5 -0.05  .85 -0.4  0.9 -0.7 /
\endpicture
\cr
\check R_{MN}(\tau_M\otimes \id_N) &= (\tau_M\otimes \id_N)\check R_{MN}, \cr
}
\formula$$
\medskip
$$
\beginpicture
\setcoordinatesystem units <1cm,.5cm>         
\setplotarea x from -2 to 2, y from 0 to 1.5    
\put{$M\otimes\, (N\otimes P)$} at -2.1 1.8
\put{$(N\otimes P)\ \otimes M$} at -2.15 -0.6
\put{$\check R_{M,N\otimes P}
= (\id_N\otimes \check R_{MP})
(\check R_{MN}\otimes \id_P)$} at -.7 -2.5
\put{$=$} at -.7 0.6
\put{$\bullet$} at  -1.65 1.2      
\put{$\bullet$} at  -2.6 1.2      %
\put{$\bullet$} at  -1.65 0          
\put{$\bullet$} at  -2.6  0          %
\put{$M\otimes N\otimes P$} at 0.9 2.5
\put{$N\otimes P\otimes M$} at 0.9 -1.4
\put{$\bullet$} at  0.15 1.9      
\put{$\bullet$} at  0.9 1.9      %
\put{$\bullet$} at  1.65 1.9      %
\put{$\bullet$} at  0.15 -0.7          
\put{$\bullet$} at  0.9 -0.7          %
\put{$\bullet$} at  1.65 -0.7          %
\plot  1.65 1.9   1.65 .6 /
\plot  0.15 .6   0.15 -0.7 /
\setquadratic
\plot  0.15 .6  .2 .9  .45 1.15 /
\plot  .65 1.35  .85 1.6  0.9 1.9 /
\plot 0.15 1.9  .2 1.6  .55 1.25  .85 .9  0.9 .6 /
\plot  0.9 -0.7  0.95 -0.4  1.15 -0.15 /
\plot  1.35 .05  1.6 .3  1.65 .6 /
\plot 0.9 .6  0.95 .3  1.25 -0.05  1.6 -0.4  1.65 -0.7 /
\plot  -2.6 0  -2.55 .3  -2.2 0.55 /
\plot  -2 .75  -1.7 1  -1.65 1.3 /
\plot -2.6 1.3  -2.55 1  -2.1 .65  -1.7 .3  -1.65 0 /
\endpicture
\qquad\qquad
\beginpicture
\setcoordinatesystem units <1cm,.5cm>         
\setplotarea x from -2 to 2, y from 0 to 1.5    
\put{$(M\otimes N)\,\otimes P$} at -2.2 1.8
\put{$P\otimes (M\otimes N)$} at -2.05 -0.6
\put{$\check R_{M\otimes N,P}
= (\check R_{MP}\otimes \id_N)
(\id_M\otimes \check R_{NP}),$} at -.7 -2.5
\put{$=$} at -.7 0.6
\put{$\bullet$} at  -1.65 1.2      
\put{$\bullet$} at  -2.6 1.2      %
\put{$\bullet$} at  -1.65 0          
\put{$\bullet$} at  -2.6  0          %
\put{$M\otimes N\otimes P$} at 0.9 2.5
\put{$P\otimes M\otimes N$} at 0.9 -1.4
\put{$\bullet$} at  0.15 1.9      
\put{$\bullet$} at  0.9 1.9      %
\put{$\bullet$} at  1.65 1.9      %
\put{$\bullet$} at  0.15 -0.7          
\put{$\bullet$} at  0.9 -0.7          %
\put{$\bullet$} at  1.65 -0.7          %
\plot  1.65 0.6   1.65 -0.7 /
\plot  0.15 1.9   0.15 0.6 /
\setquadratic
\plot  0.15 -0.7  .2 -0.4  .45 -0.15 /
\plot  .65 .05  .85 0.3  0.9 0.6 /
\plot 0.15 0.6  .2 0.3  .55 -0.05  .85 -0.4  0.9 -0.7 /
\plot  0.9 0.6  0.95 0.9  1.15 1.15 /
\plot  1.35 1.35  1.6 1.6  1.65 1.9 /
\plot 0.9 1.9  0.95 1.6  1.25 1.25  1.6 0.9  1.65 0.6 /
\plot  -2.6 0  -2.55 .3  -2.2 0.55 /
\plot  -2 .75  -1.7 1  -1.65 1.3 /
\plot -2.6 1.3  -2.55 1  -2.1 .65  -1.7 .3  -1.65 0 /
\endpicture
\formula
$$
$$
C_{M\otimes N} =
(\check R_{NM}\check R_{MN})^{-1}
(C_M\otimes C_N).
\formula
$$
The relations (2.9) and (2.10) together imply the braid relation
$$\eqalign{
\beginpicture
\setcoordinatesystem units <1cm,.5cm>         
\setplotarea x from 0 to 2, y from -3 to 2.6    
\put{$M\otimes N\otimes P$} at 0.9 2.3
\put{$P\otimes N\otimes M$} at 0.9 -2.5
\put{$\bullet$} at  0.15 1.9      
\put{$\bullet$} at  0.9 1.9      %
\put{$\bullet$} at  1.65 1.9      %
\put{$\bullet$} at  0.15 -2          
\put{$\bullet$} at  0.9 -2          %
\put{$\bullet$} at  1.65 -2          %
\plot  1.65 1.9   1.65 .6 /
\plot  0.15 .6   0.15 -0.7 /
\plot  1.65 -0.7   1.65 -2 /
\setquadratic
\plot  0.15 .6  .2 .9  .45 1.15 /
\plot  .65 1.35  .85 1.6  0.9 1.9 /
\plot 0.15 1.9  .2 1.6  .55 1.25  .85 .9  0.9 .6 /
\plot  0.9 -0.7  0.95 -0.4  1.15 -0.15 /
\plot  1.35 .05  1.6 .3  1.65 .6 /
\plot 0.9 .6  0.95 .3  1.25 -0.05  1.6 -0.4  1.65 -0.7 /
\plot  0.15 -2  .2 -1.7  .45 -1.45 /
\plot  .65 -1.25  .85 -1  0.9 -0.7 /
\plot 0.15 -0.7  .2 -1  .55 -1.35  .85 -1.7  0.9 -2 /
\endpicture
~~&=~~
\beginpicture
\setcoordinatesystem units <1cm,.5cm>         
\setplotarea x from 0 to 2, y from -3 to 2.6    
\put{$M\otimes N\otimes P$} at 0.9 2.3
\put{$P\otimes N\otimes M$} at 0.9 -2.5
\put{$\bullet$} at  0.15 1.9      
\put{$\bullet$} at  0.9 1.9      %
\put{$\bullet$} at  1.65 1.9      %
\put{$\bullet$} at  0.15 -2          
\put{$\bullet$} at  0.9 -2          %
\put{$\bullet$} at  1.65 -2          %
\plot  0.15 1.9   0.15 .6 /
\plot  1.65 .6   1.65 -0.7 /
\plot  0.15 -0.7   0.15 -2 /
\setquadratic
\plot  0.9 0.6  0.95 0.9  1.15 1.15 /
\plot  1.35 1.35  1.6 1.6  1.65 1.9 /
\plot 0.9 1.9  0.95 1.6  1.25 1.25  1.6 0.9  1.65 0.6 /
\plot  0.15 -0.7  .2 -0.4  .45 -0.15 /
\plot  .65 .05  .85 .3  0.9 .6 /
\plot 0.15 .6  .2 .3  .55 -0.05  .85 -0.4  0.9 -0.7 /
\plot  0.9 -2  0.95 -1.7  1.15 -1.45 /
\plot  1.35 -1.25  1.6 -1  1.65 -0.7 /
\plot 0.9 -0.7  0.95 -1  1.25 -1.35  1.6 -1.7  1.65 -2 /
\endpicture
\cr
(\check R_{NP}\otimes \id_M)
(\id_N\otimes \check R_{MP})
(\check R_{MN}\otimes \id_P)
&=
(\id_P\otimes \check R_{MN})
(\check R_{MP}\otimes \id_N)
(\id_M\otimes \check R_{NP}), \cr
}
\formula
$$
If $M$ is a highest weight module of weight $\lambda$
($M$ is generated by a highest weight vector $v^+$ of weight $\lambda$)
then, by [Dr, Prop. 3.2],
$$C_M = q^{-\langle \lambda,\lambda+2\rho\rangle} \id_M.
\formula
$$
Note that $\langle \lambda,\lambda+2\rho\rangle
=\langle \lambda+\rho,\lambda+\rho\rangle - \langle\rho,\rho\rangle$
are the eigenvalues of the classical Casimir operator [Dx, 7.8.5].
If $M$ is a finite dimensional $U_h\fg$ module then $M$ is a direct
sum of the irreducible modules $L(\lambda)$, $\lambda\in P^+$, and
$$C_M = \bigoplus_{\lambda\in P^+}
q^{-\langle \lambda,\lambda+2\rho\rangle}P_\lambda,
$$
where $P_\lambda\colon M\to M$ is the projection onto 
$M^{[\lambda]}$ in $M$.
From the relation (2.11) it follows that if $M=L(\mu)$,
$N=L(\nu)$ are finite dimensional irreducible $U_h\fg$ modules
then $\check R_{NM}\check R_{MN}$ acts on the
$\lambda$ isotypic component 
$L(\lambda)^{\oplus c_{\mu\nu}^\lambda}$
of the decomposition
$$L(\mu)\otimes L(\nu) = \bigoplus_\lambda L(\lambda)^{\oplus c_{\mu\nu}^\lambda}
\qquad\hbox{by the constant}\qquad
q^{\langle\lambda,\lambda+2\rho\rangle 
-\langle\mu,\mu+2\rho\rangle 
-\langle\nu,\nu+2\rho\rangle}.
\formula$$ 

Suppose that $M$ and $N$ are $U_h\fg$ modules with contravariant forms
$\langle ,\rangle_M$ and $\langle ,\rangle_N$, respectively.
Define a contravariant form on $M\otimes N$ by
$$\langle m_1\otimes n_1,m_2\otimes n_2\rangle
=\langle m_1,m_2\rangle_M \langle n_1,n_2\rangle_N,
\formula$$
for $m_1, m_2\in M$, $n_1,n_2\in N$.  If $\theta$ is the Cartan
involution defined in (2.4) then a formula of Drinfeld [Dr, Prop. 4.2]
states
$$(\theta\otimes\theta)({\cal R}) = \sum_i b_i\otimes a_i,
$$
from which it follows that
$$\eqalign{
\langle \check R_{MN}(m_1\otimes n_1), n_2\otimes m_2\rangle
&= \sum_i \langle (b_i\otimes a_i)(n_1\otimes m_1), n_2\otimes m_2\rangle  \cr
&=\sum_i \langle n_1\otimes m_1, (\theta(b_i)\otimes\theta(a_i))(n_2\otimes m_2)\rangle
\cr
&=\sum_i \langle n_1\otimes m_1, (a_i\otimes b_i)(n_2\otimes m_2)\rangle \cr
&=\sum_i \langle m_1\otimes n_1, b_im_2\otimes a_in_2\rangle. \cr
}$$
Thus
$$
\langle \check R_{MN}(m_1\otimes n_1), n_2\otimes m_2\rangle
=\langle m_1\otimes n_1, \check R_{NM}(n_2\otimes m_2)\rangle. 
\formula$$

\section 3. Affine braid group representations and the functors $F_\lambda$

There are three common ways of depicting affine braids
[Cr], [GL], [Jo3]:
\smallskip\noindent
\itemitem{(a)} As braids in a (slightly thickened) cylinder,
\smallskip\noindent
\itemitem{(b)} As braids in a (slightly thickened) annulus,
\smallskip\noindent
\itemitem{(c)} As braids with a flagpole.
\smallskip\noindent
See Figure 1.  The multiplication is by placing one cylinder
on top of another, placing one annulus inside another, or placing
one flagpole braid on top of another.  These are equivalent formulations:
an annulus can be made into a cylinder by turning up the edges,
and a cylindrical braid can be made into a flagpole braid by putting
a flagpole down the middle of the cylinder and pushing the pole
over to the left so that the strings begin and end to its right.

The {\it affine braid group} is the group
$\tilde {\cal B}_k$ formed by the affine braids with $k$ strands.
The affine braid group $\tilde {\cal B}_k$ can be presented
by generators $T_1,T_2,\ldots,T_{k-1}$ and $X^{\varepsilon_1}$
$$
T_i = 
\beginpicture
\setcoordinatesystem units <.5cm,.5cm>         
\setplotarea x from -5 to 3.5, y from -2 to 2    
\put{${}^i$} at 0 1.2      %
\put{${}^{i+1}$} at 1 1.2      %
\put{$\bullet$} at -3 .75      %
\put{$\bullet$} at -2 .75      %
\put{$\bullet$} at -1 .75      %
\put{$\bullet$} at  0 .75      
\put{$\bullet$} at  1 .75      %
\put{$\bullet$} at  2 .75      %
\put{$\bullet$} at  3 .75      %
\put{$\bullet$} at -3 -.75          %
\put{$\bullet$} at -2 -.75          %
\put{$\bullet$} at -1 -.75          %
\put{$\bullet$} at  0 -.75          
\put{$\bullet$} at  1 -.75          %
\put{$\bullet$} at  2 -.75          %
\put{$\bullet$} at  3 -.75          %
\plot -4.5 1.25 -4.5 -1.25 /
\plot -4.25 1.25 -4.25 -1.25 /
\ellipticalarc axes ratio 1:1 360 degrees from -4.5 1.25 center 
at -4.375 1.25
\put{$*$} at -4.375 1.25  
\ellipticalarc axes ratio 1:1 180 degrees from -4.5 -1.25 center 
at -4.375 -1.25 
\plot -3 .75  -3 -.75 /
\plot -2 .75  -2 -.75 /
\plot -1 .75  -1 -.75 /
\plot  2 .75   2 -.75 /
\plot  3 .75   3 -.75 /
\setquadratic
\plot  0 -.75  .05 -.45  .4 -0.1 /
\plot  .6 0.1  .95 0.45  1 .75 /
\plot 0 .75  .05 .45  .5 0  .95 -0.45  1 -.75 /
\endpicture
\qquad\hbox{and}\qquad
X^{\varepsilon_1} = 
~~\beginpicture
\setcoordinatesystem units <.5cm,.5cm>         
\setplotarea x from -5 to 3.5, y from -2 to 2    
\put{$\bullet$} at -3 0.75      %
\put{$\bullet$} at -2 0.75      %
\put{$\bullet$} at -1 0.75      %
\put{$\bullet$} at  0 0.75      
\put{$\bullet$} at  1 0.75      %
\put{$\bullet$} at  2 0.75      %
\put{$\bullet$} at  3 0.75      %
\put{$\bullet$} at -3 -0.75          %
\put{$\bullet$} at -2 -0.75          %
\put{$\bullet$} at -1 -0.75          %
\put{$\bullet$} at  0 -0.75          
\put{$\bullet$} at  1 -0.75          %
\put{$\bullet$} at  2 -0.75          %
\put{$\bullet$} at  3 -0.75          %
\plot -4.5 1.25 -4.5 -0.13 /
\plot -4.5 -0.37   -4.5 -1.25 /
\plot -4.25 1.25 -4.25  -0.13 /
\plot -4.25 -0.37 -4.25 -1.25 /
\ellipticalarc axes ratio 1:1 360 degrees from -4.5 1.25 center 
at -4.375 1.25
\put{$*$} at -4.375 1.25  
\ellipticalarc axes ratio 1:1 180 degrees from -4.5 -1.25 center 
at -4.375 -1.25 
\plot -2 0.75  -2 -0.75 /
\plot -1 0.75  -1 -0.75 /
\plot  0 0.75   0 -0.75 /
\plot  1 0.75   1 -0.75 /
\plot  2 0.75   2 -0.75 /
\plot  3 0.75   3 -0.75 /
\setlinear
\plot -3.3 0.25  -4.1 0.25 /
\ellipticalarc axes ratio 2:1 180 degrees from -4.65 0.25  center 
at -4.65 0 
\plot -4.65 -0.25  -3.3 -0.25 /
\setquadratic
\plot  -3.3 0.25  -3.05 .45  -3 0.75 /
\plot  -3.3 -0.25  -3.05 -0.45  -3 -0.75 /
\endpicture
\formula$$
with relations
\smallskip
\itemitem{(\global\advance\resultno by 1
\the\sectno.\the\resultno a)} $T_iT_j=T_jT_i$, \qquad if $|i-j|>1$,
\smallskip
\itemitem{(\the\sectno.\the\resultno b)} 
$T_iT_{i+1}T_i=T_{i+1}T_iT_{i+1}$, \qquad for $1\le i\le k-2$, 
\smallskip
\itemitem{(\the\sectno.\the\resultno c)} $X^{\varepsilon_1}T_1 X^{\varepsilon_1}T_1
=T_1 X^{\varepsilon_1}T_1 X^{\varepsilon_1}$,
\smallskip
\itemitem{(\the\sectno.\the\resultno d)} 
$X^{\varepsilon_1}T_i=T_i X^{\varepsilon_1}$, \qquad for $2\le i\le k-1$.
\smallskip
Define
$$
X^{\varepsilon_i}=T_{i-1}T_{i-2}\cdots T_2T_1 
X^{\varepsilon_1}T_1T_2\cdots T_{i-1},
\qquad 1\le i\le k.
\formula$$
By drawing pictures of the corresponding affine braids
it is easy to check that the $X^{\varepsilon_i}$ all 
commute with each other and so 
$X=\langle X^{\varepsilon_i}\ |\ 1\le i\le k\rangle$
is an abelian subgroup of $\tilde {\cal B}_k$.  Let
$L\cong \ZZ^k$ be the free abelian group generated by
$\varepsilon_1,\ldots, \varepsilon_k$.  Then
$$
L=\{ \lambda_1\varepsilon_1+\cdots+\lambda_k\varepsilon_k
\ |\ \lambda_i\in \ZZ\}
\qquad\hbox{and}\qquad
X=\{ X^\lambda \ |\ \lambda\in L\},
\formula$$
where $X^\lambda=(X^{\varepsilon_1})^{\lambda_1}
(X^{\varepsilon_2})^{\lambda_2} \cdots 
(X^{\varepsilon_k})^{\lambda_k}$, for $\lambda\in L$.

\subsection The $\tilde {\cal B}_k$ module $M\otimes V^{\otimes k}$

Let $U_h\fg$ be the Drinfeld-Jimbo quantum group
associated to a finite dimensional complex semisimple
Lie algebra $\fg$.   Let $M$ be a $U_h\fg$-module in the category
${\cal O}$ and let $V$ be a finite dimensional $U_h\fg$ module.
Define $\check R_i$, $1\le i\le k-1$, and $\check R_0^2$
in $\End_{U_h\fg}(M\otimes V^{\otimes k})$ by
$$\check R_i = \id_M\otimes \id_V^{\otimes (i-1)}
\otimes \check R_{VV}\otimes \id_V^{\otimes (k-i-1)}
\qquad\hbox{and}\qquad
\check R_0^2 = (\check R_{VM}\check R_{MV})\otimes \id_V^{\otimes (k-1)}.
$$

\prop The map defined by
$$\matrix{
\Phi\colon &\tilde {\cal B}_k &\longrightarrow
&\End_{U_h\fg}(M\otimes V^{\otimes k}) \cr
&T_i &\longmapsto &\check R_i, &\qquad &1\le i\le k-1,\cr
&X^{\varepsilon_1} &\longmapsto &\check R_0^2,\cr
}$$
makes $M\otimes V^{\otimes k}$ into a $\tilde {\cal B}_k$
module.
\pf
It is necessary to show that
\smallskip\noindent
\itemitem{(a)} $\check R_i\check R_j = \check R_j\check R_i$,
\qquad if $|i-j|>1$,
\smallskip\noindent
\itemitem{(b)}
$\check R_0^2\check R_i = \check R_i\check R_0^2$,
\quad $i>2$,
\smallskip\noindent
\itemitem{(c)}
$\check R_i\check R_{i+1}\check R_i = 
\check R_{i+1}\check R_i\check R_{i+1}$,
$1\ge i\ge k-2$,
\smallskip\noindent
\itemitem{(d)}
$\check R_0^2\check R_1\check R_0^2\check R_1 
= \check R_1\check R_0^2\check R_1\check R_0^2$.
\smallskip\noindent
The relations (a) and (b) follow immediately from the definitions
of $\check R_i$ and $\check R_0^2$ and (c)
is a particular case of the braid relation (2.12).
The relation (d) is also a consequence of the braid relation:
$$\eqalign{
\check R_0^2\check R_1\check R_0^2\check R_1 
&=(\check R_{VM}\check R_{MV}\otimes \id)
(\id\otimes \check R_{VV})
(\check R_{VM}\check R_{MV}\otimes \id)
(\id\otimes \check R_{VV}) \cr
&=(\check R_{VM}\otimes \id)
\underbrace{(\id\otimes \check R_{VV})
(\check R_{VV}\otimes \id)
(\id\otimes \check R_{MV})}
(\check R_{MV}\otimes \id)
(\id\otimes \check R_{VV}) \cr
&=
\underbrace{
(\id\otimes \check R_{VV}) 
(\check R_{MV}\otimes \id)
(\id\otimes \check R_{VM}) 
}
\underbrace{
(\check R_{VV}\otimes \id)
(\id\otimes \check R_{MV}) 
(\check R_{MV}\otimes \id)
} \cr
&=
(\id\otimes \check R_{MV}) 
(\check R_{VM} \overbrace{\check R_{MV}\otimes \id)
(\id\otimes \check R_{VV}) 
(\check R_{VM}} \check R_{MV}\otimes \id) \cr
&= \check R_1\check R_0^2\check R_1\check R_0^2, \cr
}
$$
or equivalently,
$$
\check R_0^2\check R_1\check R_0^2\check R_1 
=
\beginpicture
\setcoordinatesystem units <.25cm,.25cm>         
\setplotarea x from 0 to 2, y from -2 to 2    
\put{$\cdot$} at  0.15 3.9      
\put{$\cdot$} at  0.9 3.9      %
\put{$\cdot$} at  1.65 3.9      %
\put{$\cdot$} at  0.15  -4          
\put{$\cdot$} at  0.9  -4          %
\put{$\cdot$} at  1.65  -4          %
\plot  0.15  -2.7   0.15  -4 /
\plot  1.65  3.9   1.65  2.6 /
\plot  1.65  2.6   1.65  1.3 /
\plot  0.15  1.3     0.15  -0.1 /
\plot  1.65  -0.1   1.65  -1.4 /
\plot  1.65  -1.4   1.65  -2.7 /
\setsolid
\setquadratic
\plot  0.15 -1.4  .2 -1.1  .45 -0.85 /
\plot  .65 -0.7  .85 -0.3  0.9 0 /
\plot 0.15 0  .2 -0.3  .55 -.8  .85 -1.1  0.9 -1.4 /
\plot  0.9 0  0.95 0.3  1.15 0.55 /
\plot  1.35 0.75  1.6 1  1.65 1.3 /
\plot 0.9 1.3  0.95 1  1.25 .65  1.6 .3  1.65 0 /
\plot  0.15 -2.7  .2 -2.4  .45  -2.15 /
\plot  .65 -2.05  .85 -1.7  0.9 -1.3 /
\plot 0.15 -1.4  .2 -1.7  .55 -2.05  .85 -2.4  0.9 -2.7 /
\plot  0.9 -4  0.95 -3.7  1.15 -3.45 /
\plot  1.35 -3.25  1.6 -3  1.65 -2.7 /
\plot 0.9 -2.7  0.95 -3  1.25 -3.35  1.6 -3.7  1.65 -4 /
\plot  0.15 2.6  .2 2.9  .45 3.15 /
\plot  .65 3.35  .85 3.6  0.9 3.9 /
\plot 0.15 3.9  .2 3.6  .55 3.25  .85 2.9  0.9 2.6 /
\plot  0.15 1.3  .2 1.6  .45 1.85 /
\plot  .65 2.05  .85 2.3  0.9 2.6 /
\plot 0.15 2.6  .2 2.3  .55 1.95  .85 1.6  0.9 1.3 /
\endpicture
=
\beginpicture
\setcoordinatesystem units <.25cm,.25cm>         
\setplotarea x from 0 to 2, y from -2 to 2    
\put{$\cdot$} at  0.15 3.9      
\put{$\cdot$} at  0.9 3.9      %
\put{$\cdot$} at  1.65 3.9      %
\put{$\cdot$} at  0.15  -4          
\put{$\cdot$} at  0.9  -4          %
\put{$\cdot$} at  1.65  -4          %
\plot  1.65  3.9   1.65  2.6 /
\plot  0.15  2.6   0.15  1.3 /
\plot  1.65  1.3   1.65  0 /
\plot  0.15  0   0.15  -1.4 /
\plot  1.65  -1.4   1.65  -2.7 /
\plot  0.15  -2.7   0.15  -4 /
\setdashes
\plot  0  2.6   2  2.6 /
\plot  0  -1.4   2  -1.4 /
\setsolid
\setquadratic
\plot  0.15 2.6  .2 2.9  .45 3.15 /
\plot  .65 3.35  .85 3.6  0.9 3.9 /
\plot 0.15 3.9  .2 3.6  .55 3.25  .85 2.9  0.9 2.6 /
\plot  0.9 1.3  0.95 1.6  1.15  1.85 /
\plot  1.35 2.05  1.6 2.3  1.65 2.6 /
\plot 0.9 2.6  0.95 2.3  1.25 1.95  1.6 1.6  1.65 1.3 /
\plot  0.15 0  .2 .3  .45  0.55 /
\plot  .65 .75  .85 1  0.9 1.3 /
\plot 0.15 1.3  .2 1  .55 .65  .85 .3  0.9 0 /
\plot  0.9 -1.4  0.95 -1.1  1.15 -0.85 /
\plot  1.35 -0.7  1.6 -0.3  1.65 0 /
\plot 0.9 -0  0.95 -0.3  1.25 -.8  1.6 -1.1  1.65 -1.4 /
\plot  0.15 -2.7  .2 -2.4  .45 -2.15 /
\plot  .65 -1.95  .85 -1.7  0.9 -1.4 /
\plot 0.15 -1.4  .2 -1.7  .55 -2.05  .85 -2.4  0.9 -2.7 /
\plot  0.9 -4  0.95 -3.7  1.15 -3.45 /
\plot  1.35 -3.25  1.6 -3  1.65 -2.7 /
\plot 0.9 -2.7  0.95 -3  1.25 -3.35  1.6 -3.7  1.65 -4 /
\endpicture
=
\beginpicture
\setcoordinatesystem units <.25cm,.25cm>         
\setplotarea x from 0 to 2, y from -2 to 2    
\put{$\cdot$} at  0.15 3.9      
\put{$\cdot$} at  0.9 3.9      %
\put{$\cdot$} at  1.65 3.9      %
\put{$\cdot$} at  0.15  -4          
\put{$\cdot$} at  0.9  -4          %
\put{$\cdot$} at  1.65  -4          %
\plot  0.15  3.9   0.15  2.6 /
\plot  1.65  2.6   1.65  1.3 /
\plot  0.15  1.3   0.15  0 /
\plot  1.65  0   1.65  -1.4 /
\plot  0.15  -1.4   0.15  -2.7 /
\plot  1.65  -2.7   1.65  -4 /
\setdashes
\plot  0  0   2  0 /
\setsolid
\setquadratic
\plot  0.15 -1.4  .2 -1.1  .45 -0.85 /
\plot  .65 -0.7  .85 -0.3  0.9 0 /
\plot 0.15 0  .2 -0.3  .55 -.8  .85 -1.1  0.9 -1.4 /
\plot  0.9 -2.7  0.95 -2.4  1.15  -2.15 /
\plot  1.35 -1.95  1.6 -1.7  1.65 -1.4 /
\plot 0.9 -1.4  0.95 -1.7  1.25 -2.05  1.6 -2.4  1.65 -2.7 /
\plot  0.15 -4  .2 -3.7  .45  -3.45 /
\plot  .65 -3.25  .85 -3  0.9 -2.7 /
\plot 0.15 -2.7  .2 -3  .55 -3.35  .85 -3.7  0.9 -4 /
\plot  0.9 2.6  0.95 2.9  1.15 3.15 /
\plot  1.35 3.3  1.6 3.7  1.65 4 /
\plot 0.9 4  0.95 3.7  1.25 3.2  1.6 2.9  1.65 2.6 /
\plot  0.15 1.3  .2 1.6  .45 1.85 /
\plot  .65 2.05  .85 2.3  0.9 2.6 /
\plot 0.15 2.6  .2 2.3  .55 1.95  .85 1.6  0.9 1.3 /
\plot  0.9 0  0.95 0.3  1.15 0.55 /
\plot  1.35 0.75  1.6 1  1.65 1.3 /
\plot 0.9 1.3  0.95 1  1.25 .65  1.6 .3  1.65 0 /
\endpicture
=
\beginpicture
\setcoordinatesystem units <.25cm,.25cm>         
\setplotarea x from 0 to 2, y from -2 to 2    
\put{$\cdot$} at  0.15 3.9      
\put{$\cdot$} at  0.9 3.9      %
\put{$\cdot$} at  1.65 3.9      %
\put{$\cdot$} at  0.15  -4          
\put{$\cdot$} at  0.9  -4          %
\put{$\cdot$} at  1.65  -4          %
\plot  0.15  3.9   0.15  2.6 /
\plot  1.65  2.6   1.65  1.3 /
\plot  1.65  1.3   1.65  0 /
\plot  0.15  0     0.15  -1.4 /
\plot  1.65  -1.4   1.65  -2.7 /
\plot  1.65  -2.7   1.65  -4 /
\setdashes
\plot  0  1.3   2  1.3 /
\plot  0  -2.7   2  -2.7 /
\setsolid
\setquadratic
\plot  0.15 -2.7  .2 -2.4  .45 -2.15 /
\plot  .65 -1.95  .85 -1.7  0.9 -1.4 /
\plot 0.15 -1.4  .2 -1.7  .55 -2.05  .85 -2.4  0.9 -2.7 /
\plot  0.9 -1.4  0.95 -1.1  1.15  -0.85 /
\plot  1.35 -0.7  1.6 -0.3  1.65 0 /
\plot 0.9 0  0.95 -0.3  1.25 -0.8  1.6 -1.1  1.65 -1.4 /
\plot  0.15 -4  .2 -3.7  .45  -3.45 /
\plot  .65 -3.25  .85 -3  0.9 -2.7 /
\plot 0.15 -2.7  .2 -3  .55 -3.35  .85 -3.7  0.9 -4 /
\plot  0.9 2.6  0.95 2.9  1.15 3.15 /
\plot  1.35 3.3  1.6 3.7  1.65 4 /
\plot 0.9 4  0.95 3.7  1.25 3.2  1.6 2.9  1.65 2.6 /
\plot  0.15 1.3  .2 1.6  .45 1.85 /
\plot  .65 2.05  .85 2.3  0.9 2.6 /
\plot 0.15 2.6  .2 2.3  .55 1.95  .85 1.6  0.9 1.3 /
\plot  0.15 0  .2 0.3  .45 0.55 /
\plot  .65 .75  .85 1  0.9 1.3 /
\plot 0.15 1.3  .2 1  .55 .65  .85 .3  0.9 0 /
\endpicture
= \check R_1\check R_0^2\check R_1\check R_0^2. 
\qquad\hbox{\qed}
$$
\medskip

A $\tilde {\cal B}_k$ module $N$ is {\it calibrated} if the abelian 
group $X$ defined in (3.4) acts semisimply on $N$, i.e. if $N$
has a basis of simultaneous eigenvectors for the
action of $X^{\varepsilon_1},\ldots,X^{\varepsilon_k}$.

\prop If $M$ and $V$ are finite dimensional $U_h\fg$ modules
then the $\tilde {\cal B}_k$ module $M\otimes V^{\otimes k}$
defined in Proposition 3.5 is calibrated.
\pf
Let $P^+$ be the set of dominant integral weights.
Since $M$ and $V$ are finite dimensional the $U_h\fg$-module
$M\otimes V^{\otimes i}$
is semisimple for every $1\le i\le k$ and
$$M\otimes V^{\otimes i} =
\bigoplus_{\lambda\in P^+} (M\otimes V^{\otimes i})^{[\lambda]}
\cong \bigoplus_{\lambda\in P^+} L(\lambda)^{\oplus m_\lambda},$$
where $m_\lambda\in \ZZ_{\ge 0}$ and
$(M\otimes V^{\otimes i})^{[\lambda]}
\cong \bigoplus_{\lambda\in P^+} L(\lambda)^{\oplus m_\lambda}.$
Given a basis of $M\otimes V^{\otimes (i-1)}$ which respects the
decomposition $M\otimes V^{\otimes (i-1)} = \bigoplus_\mu
(M\otimes V^{\otimes (i-1)})^{[\mu]}$  one can construct a basis
of $M\otimes V^{\otimes i}$ which respects the decomposition
$$M\otimes V^{\otimes i}=(M\otimes V^{\otimes (i-1)})\otimes V
=\bigoplus_{\lambda,\mu,\nu}
((M\otimes V^{\otimes (i-1)})^{[\mu]}\otimes V^{[\nu]})^{[\lambda]}.
$$
Since 
$((M\otimes V^{\otimes (i-1)})^{[\mu]}\otimes V^{[\nu]})^{[\lambda]}
\subseteq
(M\otimes V^{\otimes i})^{[\lambda]}$ this new basis respects the
decomposition
$M\otimes V^{\otimes i} = \bigoplus_\lambda
(M\otimes V^{\otimes i})^{[\lambda]}$.
This procedure produces, inductively, a basis $B$ of
$M\otimes V^{\otimes k}$ which respects the decompositions
$$M\otimes V^{\otimes k}
=(M\otimes V^{\otimes i})\otimes V^{\otimes (k-i)}
=\bigoplus_\lambda (M\otimes V^{\otimes i})^{[\lambda]}
\otimes V^{\otimes (k-i)},$$
for all $0\le i\le k$.  The central element $e^{-h\rho}u$
in $U_h\fg$ acts on $(M\otimes V^{\otimes i})^{[\lambda]}$
by the constant $q^{-\langle \lambda,\lambda+2\rho\rangle}$.
From (2.10), (2.11) and (2.14) it follows that $X^{\varepsilon_i}$ acts on
$M\otimes V^{\otimes k}$ by
$$\eqalign{
\check R_{i-1}\cdots \check R_1\check R_0^2
\check R_1\cdots \check R_{i-1}
&=\check R_{V,M\otimes V^{\otimes (i-1)}}
\check R_{V,M\otimes V^{\otimes (i-1)}}
\otimes \id_V^{\otimes (k-i)} \cr
&=(C_{M\otimes V^{\otimes (i-1)}}\otimes C_V)
C^{-1}_{M\otimes V^{\otimes i}}
\otimes \id_V^{\otimes (k-i)} \cr 
&= \sum_{\lambda,\mu,\nu}
q^{\langle \lambda,\lambda+2\rho\rangle
-\langle \mu,\mu+2\rho\rangle-\langle \nu,\nu+2\rho\rangle}
P_{\mu\nu}^\lambda\otimes \id_V^{\otimes (k-i)} \cr
}$$
where $P_{\mu\nu}^\lambda\colon M\otimes \id_V^{\otimes i}\to
M\otimes \id_V^{\otimes i}$ is the projection onto
$((M\otimes V^{\otimes (i-1)})^{[\mu]}\otimes V^{[\nu]})^{[\lambda]}$.
Thus $X^{\varepsilon_i}$ acts diagonally on the basis $B$.
\endpf

Define an anti-involution on $\tilde {\cal B}_k$ by
$$\tilde\theta(T_i)=T_i
\qquad\hbox{and}\qquad
\tilde\theta(X^\lambda)=X^\lambda,$$
for $1\le i\le k-1$ and $\lambda\in L$.
A {\it contravariant form} on
a $\tilde {\cal B}_k$ module $N$ is a symmetric bilinear
form $\langle,\rangle\colon N\times N\to \CC$ such that
$$\langle bn_1,n_2\rangle =
\langle n_1, \tilde\theta(b)n_2\rangle
\qquad\hbox{for $n_1,n_2\in N$, $b\in \tilde {\cal B}_k$.}
$$
Suppose $M$ is a $U_h\fg$-module in the category ${\cal O}$ and
$V$ is a finite dimensional $U_h\fg$ module.  Let
$\langle ,\rangle_{{}_M}$ and $\langle , \rangle_{{}_V}$ be 
$U_h\fg$-contravariant forms on $M$ and $V$ respectively.  By (2.16),
$$\langle \check R_{VV}(v_1\otimes v_2), v'_1\otimes v'_2\rangle
=\langle v_1\otimes v_2, \check R_{VV}(v'_1\otimes v'_2)\rangle
$$
for $v_1,v_2,v_1',v_2'\in V$, and
$$\langle \check R_{VM}\check R_{MV}(m\otimes v), m'\otimes v'\rangle
=\langle \check R_{MV}(m\otimes v), \check R_{MV}(m'\otimes v')\rangle
=\langle m\otimes v, \check R_{VM}\check R_{MV}(m'\otimes v')\rangle
$$
for $m,m'\in M$, $v,v'\in V$.  Thus it follows that the form
$\langle,\rangle$ on $M\otimes V^{\otimes k}$ given by 
$$
\langle m\otimes v_1\otimes \cdots v_k,
m'\otimes v'_1\otimes \cdots v'_k\rangle
=\langle m,m'\rangle_{{}_M}
\langle v_1,v'_1\rangle_{{}_V}\langle v_2,v'_2\rangle_{{}_V}\cdots\langle v_k,v'_k\rangle_{{}_V},
\formula$$
for $m,m'\in M$, $v_i,v'_i\in V$ is a $\tilde {\cal B}_k$
contravariant form on the $\tilde {\cal B}_k$ module
$M\otimes V^{\otimes k}$.

\vfill\eject

\subsection The functor $F_\lambda$

Fix a finite dimensional $U_h\fg$ module $V$ and an
integrally dominant weight $\lambda$ in $\fh^*$.  Let
$\tilde {\cal O}_k$ be the category of finite dimensional
$\tilde {\cal B}_k$ modules and define a functor
$$\matrix{
F_\lambda\colon & {\cal O} &\longrightarrow &\tilde {\cal O}_k \hfill\cr
&M &\longmapsto &\Hom_{U_h\fg}(M(\lambda), M\otimes V^{\otimes k}).\cr
}\formula$$

\prop Let $\lambda$ be a integrally dominant weight in $\fh^*$.
The functor $F_\lambda$ is exact.
\pf  The functor
$F_\lambda$ is the composition of two functors:
the functor $\cdot\otimes V^{\otimes k}$ and the functor
$\Hom_U(M(\lambda),\cdot)$.  The first is exact since
$V^{\otimes k}$ is finite dimensional and the second is exact
because when $\lambda$ is integrally dominant 
$M(\lambda)$ is projective, see [Jz, p. 72].
\endpf

The following proposition gives equivalent ways of
expressing the $\tilde{\cal B}_k$-module $F_\lambda(M)$.
We use the notation
$$\fn^-(M\otimes V^{\otimes k}) = \sum_i Y_i(M\otimes V^{\otimes k}),
\formula$$
where the $Y_i$ are the Chevalley generators of $\fn^-$.
In the case of $U\fg$-modules, the notation 
$\fn^-(M\otimes V^{\otimes k})$ is self explanatory---the notation
in (3.10) is simply a way to define the same object for the 
quantum group $U_h\fg$.

\prop  Let $M$ be a $U_h\fg$ module in the category
${\cal O}$ and let $V$ be a finite dimensional $U_h\fg$-module.  
Let $\lambda$ be an integrally dominant weight.  Then
$$\Hom_{U_h\fg}(M(\lambda),M\otimes V^{\otimes k})
\cong ((M\otimes V^{\otimes k})^{[\lambda]})_\lambda
\cong \left({M\otimes V^{\otimes k}\over
\fn^-(M\otimes V^{\otimes k})}\right)_\lambda$$
as $\tilde {\cal B}_k$ modules.
\pf
Since the action of $\tilde {\cal B}_k$ on
$M\otimes V^{\otimes k}$ commutes with the action of 
$U_h\fg$ on $M\otimes V^{\otimes k}$, all three vector
spaces in the statement are $\tilde {\cal B}_k$ modules,
and in all three cases, the $\tilde {\cal B}_k$ action comes from
the $\tilde {\cal B}_k$ action on $M\otimes V^{\otimes k}$.
The isomorphisms come from the fact that these vector spaces
are naturally identified with the vector space of highest weight
vectors of weight $\lambda$ in $M\otimes V^{\otimes k}$.
This identification is done as follows.

\smallskip\noindent
(a)  If $m\otimes n$ is a highest weight vector of weight
$\lambda$ in $M\otimes V^{\otimes k}$ and $v_\lambda^+$ is the
highest weight vector of weight $\lambda$ in the Verma module
$M(\lambda)$ then
$$\matrix{
\phi\colon &M(\lambda) &\to &M\otimes V^{\otimes k} \cr
&v^+ &\longmapsto &m\otimes n \cr
}$$
uniquely determines a homomorphism in 
$\Hom_U(M(\lambda),M\otimes V^{\otimes k})$.
So $\Hom_U(M(\lambda),M\otimes V^{\otimes k})$ can be identified with the
space of highest weight vectors of weight $\lambda$ in
$M\otimes V^{\otimes k}$.

\smallskip\noindent
(b)  If $m\otimes n$ is a highest weight vector of weight $\mu$
in $M\otimes V^{\otimes k}$ then there is a unique integrally
dominant weight $\nu$ such that $\mu\in W\circ \nu$.
Since $\lambda$ is integrally dominant any highest weight
vector of weight $\lambda$ in $M\otimes V^{\otimes k}$ is an
element of $(M\otimes V^{\otimes k})^{[\lambda]}$.
Furthermore, 
$$(M\otimes V^{\otimes k})^{[\lambda]}
=\bigoplus_{\mu\le\lambda} 
((M\otimes V^{\otimes k})^{[\lambda]})_\mu
\formula$$
where the sum is over all $\mu\le \lambda$ in dominance
i.e., over all $\mu$ such that $\mu=\lambda-\nu$ with
$\nu$ a nonnegative linear combination of positive
roots.  Thus $((M\otimes V^{\otimes k})^{[\lambda]})_\lambda$
consists exactly of the highest weight vectors of weight
$\lambda$.

\smallskip\noindent
(c) It follows from (3.12) that $(\fn^-M^{[\lambda]})_\lambda=0$.
So the canonical surjection $M^{[\lambda]}\to (M^{[\lambda]}/
\fn^-M^{[\lambda]})$ produces a vector space isomorphism
$$
((M\otimes V^{\otimes k})^{[\lambda]})_\lambda
\mapright{\sim} \left( {(M\otimes V^{\otimes k})^{[\lambda]}
\over \fn^-(M\otimes V^{\otimes k})^{[\lambda]}
}\right)_\lambda.$$
The last isomorphism in the statement of the proposition
now follows from
$$
\left({(M\otimes V^{\otimes k})
\over \fn^-(M\otimes V^{\otimes k})
}\right)_\lambda
=\bigoplus_\mu \left( {(M\otimes V^{\otimes k})^{[\mu]}
\over \fn^-(M\otimes V^{\otimes k})^{[\mu]}
}\right)_\lambda
=\left( {(M\otimes V^{\otimes k})^{[\lambda]}
\over \fn^-(M\otimes V^{\otimes k})^{[\lambda]}
}\right)_\lambda,
$$
where the direct sum is over all integrally dominant weights
$\mu$. 
\endpf

\section 4. The $\tilde {\cal B}_k$ modules 
${\cal M}^{\lambda/\mu}$ and ${\cal L}^{\lambda/\mu}$

Let $\lambda$ be integrally dominant and let $\mu\in \fh^*$.
Define $\tilde {\cal B}_k$ modules
$$
{\cal M}^{\lambda/\mu}
=F_\lambda(M(\mu))
\qquad\hbox{and}\qquad
{\cal L}^{\lambda/\mu}
=F_\lambda(L(\mu)).
\formula$$ 
The following lemma is the main tool for studying the structure of
these $\tilde{\cal B}_k$ modules.

\lemma  
([Jz, Theorem 2.2], [Dx, Lemma 7.6.14])
Let $E$ be a finite dimensional $U_h\fg$ module
and let $\{e_i\}$ be a basis of $E$ consisting of weight
vectors ordered so that $i<j$ if $\wt(e_i)<\wt(e_j)$.
Suppose $M$ is a $U_h\fg$ module generated by a highest weight
vector $v^+_\mu$ of weight $\mu$.  Set
$$
M_i = \sum_{j\ge i} U_h\fn^-(v^+_\mu\otimes e_j).$$
Then
\item{(a)} $M\otimes E=M_1\supseteq M_2\supseteq
\cdots$ is a filtration of $U_h\fg$ modules such that
$M_i/M_{i+1}$ is $0$ or is a highest weight module of
highest weight $\mu+\wt(e_i)$.
\smallskip\noindent
\item{(b)} If $M=M(\mu)$ then
$M_i/M_{i+1}\cong M(\mu+\wt(e_i))$.
\endthm

The braid group ${\cal B}_k$ is the subgroup of 
$\tilde {\cal B}_k$ generated by $T_1,\ldots, T_{k-1}$.
By restriction, both ${\cal M}^{\lambda/\mu}$ and
$V^{\otimes k} = L(0)\otimes V^{\otimes k}$ are ${\cal B}_k$
modules.

There is a unique $U_h\fg$ contravariant form 
$\langle , \rangle_M$ on the Verma module 
$M(w\circ\mu)$ determined by $\langle v_{w\circ\mu}^+,v_{w\circ\mu}^+\rangle_M=1$
where $v_{w\circ\mu}^+$ is the generating highest weight vector of 
$M(w\circ\mu)$.  As in (2.15), this form together with a
nondegenerate $U_h\fg$ contravariant form
$\langle , \rangle_V$ on $V$ gives a $U_h\fg$ contravariant forms
$\langle , \rangle_{V^{\otimes k}}$ and 
$\langle,\rangle$ on 
$V^{\otimes k}$ and $M(w\circ\mu)\otimes V^{\otimes k}$, respectively.

With these notations at hand we use Lemma 4.2 to
prove the fundamental facts about the 
$\tilde{\cal B}_k$ modules ${\cal M}^{\lambda/\mu}$ and 
${\cal L}^{\lambda/\mu}$ defined in (4.1).

\prop  Let $\lambda,\mu$ be integrally dominant weights and
$w\in W$.
\smallskip\noindent
\item{(a)} As ${\cal B}_k$ modules,
${\cal M}^{\lambda/w\circ\mu} \cong (V^{\otimes k})_{\lambda-w\circ\mu}$
\smallskip\noindent
\item{(b)}  ${\cal M}^{\lambda/w\circ\mu}\cong {\cal M}^{\lambda/y\circ\mu}$
if $W_{\lambda+\rho}wW_{\mu+\rho}=W_{\lambda+\rho}yW_{\mu+\rho}$.
\smallskip\noindent
\item{(c)} Use the same notation $\langle,\rangle$ 
for the $U_h\fg$ contravariant form $\langle,\rangle$ on 
$M(w\circ\mu)\otimes V^{\otimes k}$ and the $\tilde{\cal B}_k$ contravariant
form on ${\cal M}^{\lambda/(w\circ\mu)}$ obtained by restriction
of $\langle,\rangle$ to the subspace 
$(M(w\circ\mu)\otimes V^{\otimes k})^{[\lambda]}_\lambda$.  Then
$${\cal L}^{\lambda/w\circ\mu} \cong 
{{\cal M}^{\lambda/(w\circ\mu)}\over {\rm rad} \langle,\rangle}.$$
\smallskip\noindent
\item{(d)}  Assume $w$ is maximal length in $wW_{\mu+\rho}$.  If 
${\cal L}^{\lambda/(w\circ\mu)}\ne 0$ then
\itemitem{(1)} $\lambda-w\circ\mu$ is a weight of $V^{\otimes k}$,
\itemitem{(2)} $w$ is maximal length in $W_{\lambda+\rho}wW_{\mu+\rho}$.
\smallskip\noindent
\item{(e)}  If $\mu$ is a dominant integral weight then
$${\cal L}^{\lambda/\mu} \cong
\left\{ v\in (V^{\otimes k})_{\lambda-\mu}\ |\ 
X_i^{\langle \mu+\rho,\alpha_i^\vee\rangle}v=0,
\hbox{for all $1\le i\le n$.}\right\}.$$
\pf
(a) Let $v_{w\circ\mu}^+$ be the generating highest weight vector of
$M(w\circ\mu)$ and, for $n\in V^{\otimes k}$ let 
${\sl pr}(v_{w\circ\mu}^+\otimes n)$ be the image of 
$v_{w\circ\mu}^+\otimes n$ in 
$(M\otimes V^{\otimes k})/\fn^-(M\otimes V^{\otimes k})$.  Then,
since $\lambda$ is integrally dominant, Lemma 4.2 shows that 
$$\matrix{
(V^{\otimes k})_{\lambda-w\circ\mu} 
&\longrightarrow 
&{\cal M}^{\lambda/(w\circ\mu)}  \cr
n &\longmapsto &{\sl pr}(v^+_{w\circ\mu}\otimes n) \cr
}\formula$$
is a vector space isomorphism.  This is a ${\cal B}_k$-module
isomorphism since the ${\cal B}_k$ action on $M(w\circ\mu)\otimes V^{\otimes k}$
commutes with $\fn^-$ and fixes $v_{w\circ\mu}^+$.

\smallskip\noindent
(b) It is sufficient to show that ${\cal M}^{\lambda/w\circ\mu}\cong
{\cal M}^{\lambda/(s_iw\circ\mu)}$ for all simple reflections
$s_i\in W_{\lambda+\rho}$ such that $s_iw>w$.  
Applying the exact functor $F_\lambda$ to the Verma module inclusion
$$M(s_iw\circ\mu)\hookrightarrow M(w\circ\mu)
\qquad\hbox{gives}\qquad
{\cal M}^{\lambda/s_iw\circ\mu}\hookrightarrow 
{\cal M}^{\lambda/w\circ\mu},$$
an inclusion of $\tilde {\cal B}_k$-modules.
Since $s_i(\lambda-w\circ\mu)
=s_i(\lambda+\rho)-s_iw(\mu+\rho)
=\lambda+\rho-s_iw(\mu+\rho)
=\lambda-(zw)\circ\mu$ there is a (vector space) isomorphism of weight spaces 
$$(V^{\otimes})_{\lambda-w\circ\mu}\cong
V^{\otimes k}_{\lambda-s_iw\circ\mu}.$$
(This isomorphism can be realized by Lusztig's braid group action 
[CP, \S 8.1-8.2] $T_i\colon (V^{\otimes k})_{\lambda-w\circ\mu}
\to (V^{\otimes k})_{s_i(\lambda-w\circ\mu)}$).
Thus, by part (a), the $\tilde {\cal B}_k$-module inclusion 
${\cal M}^{\lambda/s_iw\circ\mu}\hookrightarrow {\cal M}^{\lambda/w\circ\mu}$
is an isomorphism.

\smallskip\noindent
(c) 
Use the notations for the bilinear forms on $M(w\circ\mu)$
and $V^{\otimes k}$ as given in the paragraph before the statement 
of the proposition.  
Let $\{b_i\}$ be an orthonormal basis of
$V^{\otimes k}$ with respect to $\langle,\rangle_{V^{\otimes k}}$.
If $r\in \rad \langle,\rangle_M$ then
$$\langle r\otimes b, s\otimes b'\rangle
=\langle r,s\rangle_M \langle b,b'\rangle_{V^{\otimes k}} =0,
\qquad\hbox{for all $s\in M(w\circ\mu)$, $b,b'\in V^{\otimes k}$,}$$
and so $(\rad\langle,\rangle_M)\otimes V^{\otimes k}
\subseteq \rad\langle,\rangle$.
Conversely, if $r_i\in M(w\circ\mu)$ such that $\sum r_i\otimes b_i
\in \rad\langle,\rangle$ then
$$0=\big\langle \sum_i r_i\otimes b_i, s\otimes b_j\big\rangle
=\sum_i \langle r_i,s\rangle_M\delta_{ij} = \langle r_i,s\rangle,
\qquad\hbox{for all $s\in M(w\circ\mu)$.}$$
So $r_i\in \rad\langle,\rangle_M$ and thus 
$\rad\langle,\rangle\subseteq
\rad\langle,\rangle_M\otimes V^{\otimes k}$.
By the $U_h\fg$ contravariance of $\langle,\rangle$
$$
(M(w\circ\mu)\otimes V^{\otimes k})^{[\lambda]}_\lambda
\perp (M(w\circ\mu)\otimes V^{\otimes k})^{[\mu]}_\mu 
$$
for integrally dominant weights $\lambda$, $\mu$ with $\lambda\ne \mu$.
Thus
$$\rad\langle,\rangle 
= \big(\rad\langle,\rangle_M\otimes V^{\otimes k}\big)^{[\lambda]}_\lambda,
\formula$$
where $\langle,\rangle$ is the restriction of the form
on $M(w\circ\mu)\otimes V^{\otimes k}$ to 
$(M(w\circ\mu)\otimes V^{\otimes k})^{[\lambda]}_\lambda$.
Thus
$$\eqalign{
{(M(w\circ\mu)\otimes V^{\otimes k})^{[\lambda]}_\lambda
\over 
\rad\langle,\rangle }
&= {(M(w\circ\mu)\otimes V^{\otimes k})^{[\lambda]}_\lambda
\over (\rad\langle,\rangle_M\otimes V^{\otimes k})^{[\lambda]}_\lambda } \cr
&\cong \left( {M(w\circ\mu)\over \rad\langle,\rangle_M}
\otimes V^{\otimes k}\right)^{[\lambda]}_\lambda
=(L(w\circ\mu)\otimes V^{\otimes k})^{[\lambda]}_\lambda
= {\cal L}^{\lambda/(w\circ\mu)}, \cr
}$$
where the isomorphism is a consequence of the fact that,
because $\lambda$ is an integrally dominant weight,
the functor $(\cdot\otimes V^{\otimes k})^{[\lambda]}_\lambda$
is exact (Prop. 3.9).

\smallskip\noindent
(d)  If $\lambda-w\circ\mu$ is not a weight of $V^{\otimes k}$ then,
by part (a), ${\cal M}^{\lambda/w\circ\mu}=0$.  Since the
functor $F_\lambda$ is exact and 
$L(w\circ\mu)$ is a quotient of $M(w\circ \mu)$, ${\cal L}^{\lambda/w\circ\mu}$
is a quotient of ${\cal M}^{\lambda/w\circ\mu}$.  Thus
${\cal M}^{\lambda/w\circ\mu}=0$ implies ${\cal L}^{\lambda/w\circ\mu}=0$.

Assume that $w$ is not the longest element of 
$W_{\lambda+\rho}w W_{\mu+\rho}$.  Then there is a positive
root $\alpha>0$ such that $s_\alpha\in W_{\lambda+\rho}$ and $s_\alpha w>w$.
Since $s_\alpha w W_{\mu+\rho}\ne w W_{\mu+\rho}$ there is an inclusion of 
Verma modules $M(s_\alpha w\circ \mu) \subseteq M(w\circ\mu)$ and
$F_\lambda(L(\mu))$ is a quotient of 
$F_\lambda(M(w\circ\mu))/F_\lambda(M(s_\alpha w\circ\mu))$.
On the other hand, by part (b), 
$${\cal M}^{\lambda/s_\alpha w\circ\mu} \cong
{\cal M}^{\lambda/w\circ\mu},
\qquad\hbox{and so}\qquad
{{\cal M}^{\lambda/w\circ\mu}\over {\cal M}^{\lambda/s_\alpha w\circ\mu} }
={F_\lambda(M(w\circ\mu))\over F_\lambda(M(s_\alpha w\circ\mu))=0}.$$
Thus $F_\lambda(L(w\circ \mu))=0$.

\smallskip\noindent
(e)  When $\mu$ is a dominant integral weight 
$$\rad\langle,\rangle_M = 
U_h\fn^-\{ Y_i^{\langle\mu+\rho,\alpha_i^\vee\rangle}v_\mu^+\ |\ 1\le i\le n\}
=\sum_i U_h\fn^-Y_i^{\langle\mu+\rho,\alpha_i^\vee\rangle} v_\mu^+,$$
see [Dx, 7.2.7].  Thus, by (c) and the vector space isomorphism (4.4)
it follows that, as vector spaces,
$${\cal L}^{\lambda/\mu} \cong
\left(
\big(
\hbox{span-}\big\{ {\sl pr}(v_\mu^+\otimes n)
\ |\ n\in V^{\otimes k}\big\}
\big) \Big/ \big(
\hbox{span-}\big\{ {\sl pr}(Y_i^{\langle \mu+\rho,\alpha_i^\vee\rangle}
v_\mu^+\otimes n)
\ |\ n\in V^{\otimes k}\big\}
\big)\right)_{\lambda-\mu}.$$
For any $k\ge 0$, 
${\sl pr}(Y_i^{k+1}\otimes n)
=
{\sl pr}(Y_i(Y_i^k v_\mu^+\otimes n) - Y_i^kv_\mu^+\otimes Y_i n)
=
-\;{\sl pr}(Y_i^kv_\mu^+\otimes Y_in),$
and so, by induction, ${\sl pr}(Y_i^{k+1}v_\mu^+\otimes n)
~=~\xi\cdot\;{\sl pr}(v_\mu^+\otimes Y_i^{k+1}n)$ for some $\xi\in \CC$,
$\xi\ne 0$.  Thus ${\cal L}^{\lambda/\mu}$ is isomorphic to the vector space
$$\left(
V^{\otimes k}\Big/
\Big(\sum_i Y_i^{\langle \mu+\rho,\alpha_i^\vee\rangle}
V^{\otimes k} \Big)\right)_{\lambda-\mu}.$$
If $b\in 
(Y_i^{\langle\mu+\rho,\alpha_i^\vee\rangle}V^{\otimes k})^\perp$
then the $U_h\fg$ contravariance of $\langle,\rangle_{V^{\otimes k}}$
gives that 
$$0=
\big\langle Y_i^{\langle\mu+\rho,\alpha_i^\vee\rangle} n,
b\big\rangle_{V^{\otimes k}}
=\big\langle n,
X_i^{\langle\mu+\rho,\alpha_i^\vee\rangle}b\big\rangle_{V^{\otimes k}},
\qquad\hbox{for all $n\in V^{\otimes k}$.}$$
Thus, by the nondegeneracy of $\langle,\rangle_{V^{\otimes k}}$,
$${\cal L}^{\lambda/\mu}\cong
\left(\sum_i Y_i^{\langle\mu+\rho,\alpha_i^\vee\rangle}
V^{\otimes k}\right)^{\perp}_{\lambda-\mu}
=\{ b\in (V^{\otimes k})_{\lambda-\mu} \ |\ 
X_i^{\langle \mu+\rho,\alpha_i^\vee\rangle} b=0\}.
\qquad\hbox{\qed}
$$

\remark
In the case when $\fg$ is type $A_{n-1}$ and $V=L(\omega_1)$ is the 
$n$-dimensional fundamental representation the converse to Proposition
4.3b also holds (see [Su] Prop. 2.3.4 and [Ze2] Th. 6.1b).  
The following example shows that this is not
true in general.  In the notation of Section 6, let 
$$\hbox{$\fg$ be of type $D_n$},\qquad
\mu = \varepsilon_1+\cdots+\varepsilon_{n-1},
\qquad\hbox{and}\qquad
V=L(\omega_1),$$
the $2n$-dimensional fundamental representation.
If 
$\lambda^\pm= \varepsilon_1+\cdots+\varepsilon_{n-1}\pm\varepsilon_n$
then ${\cal M}^{\lambda^+/\mu}$ and ${\cal M}^{\lambda^-/\mu}$ 
are isomorphic (one dimensional and simple) $\tilde {\cal B}_1$ modules.

\medskip

Proposition 4.3d gives a necessary condition on 
$\lambda/w\circ\mu$ for the $\tilde {\cal B}_k$-module 
${\cal L}^{\lambda/w\circ \mu}$ to be nonzero.  The following lemma gives an 
alternative characterization of this condition.  This will be useful for 
analyzing the combinatorics of the examples in Section 6.

\lemma  
Let $P$ be the weight lattice and let $\lambda$ be an integrally dominant
weight.  Then
$W_{\lambda+\rho}$ acts on $\lambda-P$ by the dot action.
This action has fundamental domain
$$C_{\lambda+\rho}^- = \{ \mu\in \lambda-P\ |\ 
\langle \mu+\rho,\alpha^\vee\rangle\in \ZZ_{\le 0}
\hbox{\ for all $\alpha>0$ such that 
$\langle \lambda+\rho,\alpha^\vee\rangle=0$}
\}.$$
The following are equivalent.
\itemitem{(a)} $\mu\in C_{\lambda+\rho}^-$,
\itemitem{(b)}
$\mu=w^\lambda\circ \nu$ with $\nu$ integrally dominant and $w^\lambda$ longest
in $W_{\lambda+\rho}w^\lambda$.
\itemitem{(c)}
$\mu=w^\lambda_{\tilde\mu}\circ \tilde\mu$ with $w^\lambda_{\tilde\mu}$ longest
in $W_{\lambda+\rho}w^\lambda_{\tilde\mu} W_{\mu+\rho}$.
\pf
(b) and (c) are equivalent since $W_{\mu+\rho}$ is the stabilizer of 
$\mu$ under the $\circ$ action.  
\smallskip\noindent
(b) $\Longrightarrow$ (a):  If $s_\alpha\in W_{\lambda+\rho}$ then
$\langle \lambda+\rho,\alpha^\vee\rangle = 0$ and $\ell(s_\alpha w^\lambda)<\ell(w^\lambda)$.
So $(w^\lambda)^{-1}\alpha<0$ and $\langle \mu+\rho,\alpha^\vee\rangle\in \ZZ$ since
$\lambda-\nu\in P$.  Thus,
$$\langle\mu+\rho,\alpha^\vee\rangle=\langle w^\lambda\circ\mu + \rho,\alpha^\vee\rangle
=\langle w^\lambda(\mu+\rho),\alpha^\vee\rangle
=\langle \nu+\rho, (w^\lambda)^{-1}\alpha^\vee\rangle \in \ZZ_{\le 0},$$
since $\nu$ is integrally dominant.  So $\mu\in C_{\lambda+\rho}^-$.

\smallskip\noindent
(a) $\Longrightarrow$ (b):  Let $\mu\in C_{\lambda+\rho}^-$ and fix
$\nu$ integrally dominant and $w\in W$ such that $\mu= w\circ \nu$.
Let $\alpha>0$ such that $\langle \lambda+\rho,\alpha^\vee\rangle = 0$.
Then
$$\langle \nu+\rho,w^{-1}\alpha^\vee\rangle = \langle \mu+\rho,\alpha^\vee\rangle
\in \ZZ_{\le 0}$$
and so $w^{-1}\alpha<0$.  So $\ell(s_\alpha w)<\ell(w)$.
Since this is true for all $\alpha>0$ such that 
$\langle \lambda+\rho,\alpha^\vee\rangle = 0$ it follows that 
$w$ is maximal length in its coset $W_{\lambda+\rho}w$.
\endpf

In the classical case, when $\fg$ is type $A_n$ and $V=L(\omega_1)$
is the $n+1$ dimensional fundamental representation the
$\tilde {\cal B}_k$-module ${\cal L}^{\lambda/w\circ\mu}$
is a simple $\tilde {\cal B}_k$-module whenever it is nonzero
(see [Su]).  As the following Proposition shows,
this is a very special phenomenon.

\prop  Assume that $V=L(\nu)$ for a dominant integral weight $\nu$.
If the $\tilde{\cal B}_k$-module $F_\lambda(\mu)$ is irreducible
(or $0$) for all $k$, all dominant integral weights $\mu$, and all 
integrally dominant weights $\lambda$ then
\smallskip
\item{(a)} $\fg$ is type $A_n$, $B_n$, $C_n$ or $G_2$ and $V = L(\omega_1)$, and
\smallskip
\item{(b)} the action of the subgroup ${\cal B}_k$ of $\tilde {\cal B}_k$ generates
$\End_{U_h\fg}(V^{\otimes k})$.
\pf
(a) If $\mu$ is large dominant integral weight (for example,
we may take $\mu=n\rho$, $n>>0$) then, as a $U_h\fg$-module,
$$L(\mu)\otimes V \cong \bigoplus_b L(\mu+{\rm wt}(b)),$$
where the sum is over a basis of $V$ consisting of
weight vectors and ${\rm wt}(b)$ is the weight of the vector $b$.
The group $\tilde {\cal B}_1$ is generated by the element
$X^{\varepsilon_1}$ which acts on a summand $L(\lambda)$ in 
$L(\mu)\otimes V$ by the constant 
$q^{\langle \lambda,\lambda+2\rho\rangle
-\langle \mu,\mu+2\rho\rangle - \langle \nu,\nu+2\rho\rangle}$.
Then $F_\lambda(L(\mu))$ is the $L(\lambda)$ isotypic component
of $L(\mu)\otimes V$ and these are simple $\tilde {\cal B}_1$ modules
only if all the values
$$
\eqalign{
\langle \mu+{\rm wt}(b),\mu+{\rm wt}(b)+2\rho\rangle
&-\langle \mu,\mu+2\rho\rangle - \langle \nu,\nu+2\rho\rangle \cr
&=2\langle \mu+\rho,{\rm wt}(b)\rangle +
\langle {\rm wt}(b),{\rm wt}(b)\rangle-\langle\nu,\nu+2\rho\rangle, \cr
}
\formula$$
as $b$ ranges over a weight basis of $V$, are distinct.  
It follows that all weight spaces of $V$ must be one
dimensional.  This means that 
\itemitem{(a)} $\fg$ is type $A_n$, $B_n$, $C_n$, $D_n$, $E_6$, $E_7$ or $G_2$
and $V=L(\omega_1)$, or
\itemitem{(b)} $\fg$ is type $A_n$ and $V=L(k\omega_1)$ or $V=(k\omega_n)$
for some $k$, or
\itemitem{(c)} $\fg$ is type $B_n$ and $V=L(\omega_n)$, or
\itemitem{(d)} $\fg$ is type $D_n$ and $V=L(\omega_{n-1})$ or $V=L(\omega_n)$.
\smallskip\noindent
Most of the weights of these representations lie in a single $W$-orbit.
If $\gamma$ and $\gamma'$ are two distinct weights of $V$ which
are in the same $W$-orbit then
$\langle \gamma,\gamma\rangle = \langle \gamma',\gamma'\rangle$.
If $\mu=n\rho$ with $n>>0$ then the condition that all the values in (4.9)
be distinct forces that 
$$(2n+1)\langle \rho,\gamma\rangle =2\langle \mu+\rho,\gamma\rangle 
~~\ne~~
2\langle \mu+\rho,\gamma'\rangle=(2n+1)\langle \rho,\gamma'\rangle.$$
Writing $\gamma = \nu-\sum_i c_i\alpha_i$ and $\gamma'=\nu-\sum_i c_i'\alpha_i$
with $c_i,c_i'\in \ZZ_{>0}$ the last equation becomes
$$(2n+1)\cdot\sum_i c_i \ne (2n+1)\cdot\sum_i c_i'.$$
Finally, an easy case by case check verifies that the only choices of $V$ 
in (a-d) above which satisfy this last condition for all weights in the $W$-orbit 
of the highest weight are those listed in the statement of the proposition.

\smallskip\noindent
(b) Let ${\cal Z}_k=\End_{U_h\fg}(V^{\otimes k})$.  
As a $(U_h\fg,{\cal Z}_k)$-bimodule
$$V^{\otimes k} \cong \bigoplus_{\lambda} L(\lambda)\otimes {\cal Z}_k^\lambda,$$
where ${\cal Z}^\lambda_k$ is an irreducible ${\cal Z}_k$-module and
the sum is over all dominant integral weights for which the
irreducible $U_h\fg$-module $L(\lambda)$ appears in $V^{\otimes k}$.
By restriction ${\cal Z}_k^\lambda$ is an ${\cal B}_k$-module and this is the 
$\tilde {\cal B}_k$-module $F_\lambda(L(0))$ which, by assumption, is
simple.  Since $L(0)$ is the trivial module $X^{\varepsilon_1}$ acts on 
$F_\lambda(L(0))$ by the identity and so $F_\lambda(L(0))$ is simple as
a ${\cal B}_k$-module.  Thus the simple ${\cal Z}_k$-modules in $V^{\otimes k}$
coincide exactly with the simple ${\cal B}_k$-modules in $V^{\otimes k}$
and it follows that ${\cal B}_k$ generates 
${\cal Z}_k = \End_{U_h\fg}(V^{\otimes k})$.  
\endpf

\vfill\eject

\subsection Jantzen filtrations for affine braid group representations

Applying the functor $F_\lambda$ to the Jantzen filtration
of $M(\mu)$ produces a filtration of ${\cal M}^{\lambda/\mu}$,
$$
{\cal M}^{\lambda/\mu} = F_\lambda(M(\mu))
=F_\lambda(M(\mu)^{(0)}) \supseteq
F_\lambda(M(\mu)^{(1)}) \supseteq \cdots .\formula$$
An argument of Suzuki [Su, Thm. 4.3.5] shows that this filtration can be 
obtained directly from the
$\tilde {\cal B}_k$-contravariant form $\langle,\rangle_t$
on
$${\cal M}^{\lambda+t\delta/\mu+t\delta}
=F_{\lambda+t\delta}(M(\mu+t\delta)) =
(M(\mu+t\delta)\otimes V^{\otimes k}
)^{[\lambda+t\delta]}_{\lambda+t\delta}
$$
which is the restriction of the $U_h\fg$ contravariant
form $\langle,\rangle_t$ on 
$(M(\mu+t\delta)\otimes V^{\otimes k})$, see (2.5) and (3.7).  To do this
define
$$
{\cal M}^{\lambda+t\delta/\mu+t\delta}(j) 
=\{ m\in {\cal M}^{(\lambda+t\delta)/(\mu+t\delta)} \ |\ 
\langle m,n\rangle_t = t^j\CC[t]
\hbox{\ for all\ } n\in {\cal M}^{\lambda+t\delta/\mu+t\delta}
\}$$
and
$$
\left({\cal M}^{\lambda/\mu}\right)^{(j)}
= \hbox{image of ${\cal M}^{\lambda+t\delta/\mu+t\delta}(j)$
in ${\cal M}^{\lambda+t\delta/\mu+t\delta}\otimes_{\CC[t]}
\CC[t]/t\CC[t]$}$$
to obtain a filtration
$$
{\cal M}^{\lambda/\mu}
=\left({\cal M}^{\lambda/\mu}\right)^{(0)}
\supseteq
\left({\cal M}^{\lambda/\mu}\right)^{(1)}
\supseteq \cdots
\formula$$
such that the quotients
$\big({\cal M}^{\lambda/\mu}\big)^{(j)}\big)/
\big({\cal M}^{\lambda/\mu}\big)^{(j+1)}\big)$
carry nondegenerate $\tilde {\cal B}_k$ contravariant forms.
Since, for different $\lambda$, the subspaces
$(M(\mu+t\delta)\otimes V^{\otimes k}
)^{[\lambda+t\delta]}_{\lambda+t\delta}$
are mutually orthogonal with respect to the $U_h\fg$
contravariant form $\langle , \rangle_t$ on
$(M(\mu+t\delta)\otimes V^{\otimes k}
)$,
$$
(M(\mu+t\delta)(j)\otimes V^{\otimes k}
)^{[\lambda+t\delta]}_{\lambda+t\delta}
\subseteq (M(\mu+t\delta)\otimes V^{\otimes k}
)^{[\lambda+t\delta]}_{\lambda+t\delta}(j)
={\cal M}^{(\lambda+t\delta)/(\mu+t\delta)}(j).
$$
On the other hand, if $u\in (M(\mu+t\delta)\otimes V^{\otimes k}
)^{[\lambda+t\delta]}_{\lambda+t\delta}(j)$ then write
$u=\sum_i a_i\otimes b_i$ where $a_i\in M(\mu+t\delta)$ and
${b_i}$ is an orthonormal basis of $V^{\otimes k}$.
Then, for all $v\in M(\mu+t\delta)$, and all $k$,
$\langle a_k,v\rangle_t = \langle u, v\otimes b_k\rangle_t
\in t^j\CC[t]$
and so $u\in (M(\mu+t\delta)(j)\otimes V^{\otimes k}
)^{[\lambda+t\delta]}_{\lambda+t\delta}$.  So
$$F_{\lambda+t\delta}(M(\mu+t\delta))^{(j)}
={\cal M}^{(\lambda+t\delta)/(\mu+t\delta)}(j)
$$
and the filtrations in (4.10) and (4.11) are identical.  

\prop
Let $\lambda$ and $\mu$ be integrally dominant weights and
let $w,y\in W^\mu$ be elements of maximal length in 
$W_{\lambda+\rho}wW_{\mu+\rho}$ and $W_{\lambda+\rho}yW_{\mu+\rho}$
respectively.  Then multiplicities of ${\cal L}^{\lambda/y\circ \mu}$
in the filtration (4.11) are given by
$$
\sum_{j\ge 0} \left[ 
{({\cal M}^{\lambda/w\circ\mu})^{(j)}\over
({\cal M}^{\lambda/w\circ\mu})^{(j+1)} }:
{\cal L}^{\lambda/(y\circ \mu)}\right]
{\tt v}^{{1\over2}(\ell(y)-\ell(w)+j)}
=P_{wy}({\tt v}).$$
where $P_{wy}({\tt v})$ is the Kazhdan-Lusztig polynomial for the
Weyl group $W^\mu$.
\pf
Since the functor $F_\lambda$ is exact this result follows from the
Beilinson-Bernstein theorem (2.6). The condition on $y$ is necessary for the
module ${\cal L}^{\lambda/y\circ \mu}$ to be nonzero.
\endpf

\vfill\eject

\subsection  The BGG resolution for affine braid groups

Let $\mu\in \fh^*$ be such that $-(\mu+\rho)$ is dominant and regular
and let $W^\mu_J$ be a parabolic subgroup of the integral Weyl group
$W^\mu$.  Let $w_0$ be the longest element of $W_J^\mu$ and
fix $\nu = w_0\circ\mu$.  
Applying the exact functor $F_\lambda$ to the 
BGG resolution in (2.7) produces
an exact sequence of $\tilde {\cal B}_k$-modules
$$
0\to {\cal C}_N\longrightarrow \cdots \longrightarrow
{\cal C}_1\longrightarrow {\cal C}_0\longrightarrow {\cal L}^{\lambda/\nu}
\longrightarrow 0
\qquad\hbox{where}\qquad
{\cal C}_k = \bigoplus_{\ell(w)=j}
{\cal M}^{\lambda/w\circ\nu},
\formula$$
and the sum is over all $w\in W^\mu_J$ of length $j$ (in $W^\mu_J$).
Thus, in the Grothendieck
group of the category $\tilde {\cal O}_k$ of finite dimensional
$\tilde {\cal B}_k$-modules
$$
[{\cal L}^{\lambda/\nu}]
= \sum_{w\in W_J^\mu} (-1)^{\ell(w)}
[{\cal M}^{\lambda/(w\circ\nu)}]
\formula$$
where $\nu=w_0\circ \mu$ and $w_0$ is the longest element of
$W^\mu_J$.
This identity is a generalization of the classical Jacobi-Trudi
identity [Mac\ I (5.4)] for expanding Schur functions in terms of
homogeneous symmetric functions
$$
s_{\lambda/\nu}
=\sum_{w\in S_n} (-1)^{\ell(w)} h_{\lambda+
\delta-w(\nu+\delta)}.
\formula$$

\subsection Restriction of ${\cal L}^{\lambda/\mu}$
to the braid group

The braid group is the subgroup ${\cal B}_k$ of 
$\tilde {\cal B}_k$ generated by $T_1,\ldots, T_{k-1}$.
The following proposition determines the structure
of $F_\lambda(L(\mu))$ as a ${\cal B}_k$ module
when $L(\mu)$ is finite dimensional.

\prop  Let $P^+$ be the set of dominant integral weights.  Define 
the tensor product multiplicities 
$c_{\mu\nu}^\lambda$, $\lambda,\mu,\nu\in P^+$, by the 
$U_h\fg$-module decompositions
$$L(\mu)\otimes L(\nu) 
\cong \bigoplus_{\lambda\in P^+} L(\lambda)^{\oplus c_{\mu\nu}^\lambda}.$$
Then 
$$\Res^{\tilde{\cal B}_k}_{{\cal B}_k} ({\cal L}^{\lambda/\mu}) 
= \bigoplus_{\nu\in P^+}
({\cal L}^{\nu})^{\oplus c_{\mu\nu}^\lambda},
\qquad\hbox{where}\qquad {\cal L}^{\nu}={\cal L}^{\nu/0}.$$
\pf
Let us abuse notation slightly and write sums instead of direct sums.
Then, as a $(U_h\fg, {\cal B}_k)$ bimodule
$$L(\mu)\otimes V^{\otimes k} =
\sum_\lambda L(\lambda)\otimes {\cal L}^{\lambda/\mu},$$
where ${\cal L}^{\lambda/\mu}=F_\lambda(L(\mu))$.
As a $(U_h\fg, {\cal B}_k)$ bimodule
$$L(\mu)\otimes V^{\otimes k} =
L(\mu)\otimes\left(\sum_\nu L(\nu)\otimes {\cal L}^{\nu/0}\right)
=\sum_{\lambda,\nu} c_{\mu\nu}^\lambda
L(\lambda)\otimes {\cal L}^{\nu/0}.$$
Comparing coefficients of $L(\lambda)$ in these two identities
yields the formula in the statement.
\endpf

\vfill\eject

\section 5. Markov traces

A Markov trace on the affine braid group is a trace functional
which respects the inclusions $\tilde {\cal B}_1
\subseteq \tilde {\cal B}_2\subseteq \cdots$ where
$$\matrix{
\tilde {\cal B}_k &\hookrightarrow &\tilde {\cal B}_{k+1} \cr
\beginpicture
\setcoordinatesystem units <.5cm,.5cm>         
\setplotarea x from -2 to 2, y from -2 to 2.4    
\put{${}_1$} at -.8 1.8
\put{$\ldots$} at 0.2 1.8
\put{${}_k$} at 1.2 1.8
\put{$b$} at 0 0
\plot  -1.5 1   1.5 1 /
\plot  -1.5 -1   1.5 -1 /
\plot  -1.5 -1   -1.5 1 /
\plot  1.5 -1   1.5 1 /
\plot  -.8 1.5    -.8  1 /
\plot  -.4 1.5    -.4  1 /
\plot    0 1.5     0 1 /
\plot  .4 1.5    .4  1 /
\plot  .8 1.5    .8 1 /
\plot  1.2 1.5    1.2 1 /
\plot  -.8 -1.5    -.8  -1 /
\plot  -.4 -1.5    -.4  -1 /
\plot    0 -1.5     0 -1 /
\plot  .4 -1.5    .4  -1 /
\plot  .8 -1.5    .8 -1 /
\plot  1.2 -1.5  1.2 -1 /
\linethickness=1pt
\putrule  from -1.2 1.5 to  -1.2  1 
\putrule  from  -1.2 -1.5  to  -1.2  -1 
\endpicture
&\longmapsto 
&
\beginpicture
\setcoordinatesystem units <.5cm,.5cm>         
\setplotarea x from -2 to 2, y from -2 to 2.4    
\put{$b$} at 0 0
\put{${}_1$} at -.8 1.8
\put{$\ldots$} at 0.2 1.8
\put{${}_k$} at 1.2 1.8
\put{${}_{k+1}$} at 2.0 1.8
\plot  -1.5 1   1.5 1 /
\plot  -1.5 -1   1.5 -1 /
\plot  -1.5 -1   -1.5 1 /
\plot  1.5 -1   1.5 1 /
\plot  -.8 1.5    -.8  1 /
\plot  -.4 1.5    -.4  1 /
\plot    0 1.5     0 1 /
\plot  .4 1.5    .4  1 /
\plot  .8 1.5    .8 1 /
\plot  1.2 1.5    1.2 1 /
\plot  -.8 -1.5    -.8  -1 /
\plot  -.4 -1.5    -.4  -1 /
\plot    0 -1.5     0 -1 /
\plot  .4 -1.5    .4  -1 /
\plot  .8 -1.5    .8 -1 /
\plot  1.2 -1.5    1.2 -1 /
\plot  2.0 1.5    2.0 -1.5 /
\linethickness=1pt
\putrule  from -1.2 1.5 to  -1.2  1 
\putrule  from  -1.2 -1.5  to  -1.2  -1 
\endpicture
\cr
}\formula$$
More precisely,  a {\it Markov trace} on the affine braid
group with parameters $z, Q_1,Q_2,\ldots \in \CC$ is a sequence of 
functions
$$\mt_k\colon \tilde {\cal B}_k \longrightarrow \CC
\qquad\hbox{such that}$$
\itemitem{(1)} $\mt_1(1)=1$,
\smallskip\noindent
\itemitem{(2)} $\mt_{k+1}(b)=\mt_k(b)$, for $b\in \tilde {\cal B}_k$,
\smallskip\noindent
\itemitem{(3)} $\mt_k(b_1b_2)=\mt_k(b_2b_1)$, for 
$b_1,b_2\in \tilde {\cal B}_k$,
\smallskip\noindent
\itemitem{(4)} $\mt_{k+1}(bT_k)=z\mt_k(b)$, for 
$b\in \tilde {\cal B}_k$,
\smallskip\noindent
\itemitem{(5)} $\mt_{k+1}(b(\tilde X^{\varepsilon_{k+1}})^r)
=Q_r\mt_k(b)$, for 
$b\in \tilde {\cal B}_k$,
\smallskip\noindent
where
$$ \tilde X^{\varepsilon_{k+1}}
=T_kT_{k-1}\cdots T_2X^{\varepsilon_1}T_2^{-1}\cdots
T_{k-1}^{-1}T_k^{-1} = 
\enspace
\beginpicture
\setcoordinatesystem units <.5cm,.5cm>         
\setplotarea x from -5 to 3.5, y from -2 to 2    
\put{${}^1$} at -3 1.2
\put{${}^2$} at -2 1.2
\put{${}^{\cdot\ \cdot\ \cdot}$} at 0 1.2
\put{${}^{k+1}$} at 3 1.2
\put{$\bullet$} at -3 0.75      %
\put{$\bullet$} at -2 0.75      %
\put{$\bullet$} at -1 0.75      %
\put{$\bullet$} at  0 0.75      
\put{$\bullet$} at  1 0.75      %
\put{$\bullet$} at  2 0.75      %
\put{$\bullet$} at  3 0.75      %
\put{$\bullet$} at -3 -0.75          %
\put{$\bullet$} at -2 -0.75          %
\put{$\bullet$} at -1 -0.75          %
\put{$\bullet$} at  0 -0.75          
\put{$\bullet$} at  1 -0.75          %
\put{$\bullet$} at  2 -0.75          %
\put{$\bullet$} at  3 -0.75          %
\plot -4.5 1.25 -4.5 -0.13 /
\plot -4.5 -0.37 -4.5 -1.25 /
\plot -4.25 1.25 -4.25 -0.13 /
\plot -4.25 -0.37 -4.25 -1.25 /
\ellipticalarc axes ratio 1:1 360 degrees from -4.5 1.25 center 
at -4.375 1.25
\put{$*$} at -4.375 1.25  
\ellipticalarc axes ratio 1:1 180 degrees from -4.5 -1.25 center 
at -4.375 -1.25 
\plot -3 .75  -3 -0.75 /
\plot -2 .75  -2 -0.75 /
\plot -1 .75  -1 -0.75 /
\plot  0 .75   0 -0.75 /
\plot  1 .75   1 -0.75 /
\plot  2 .75   2 -0.75 /
\setlinear
\plot -3.2 0.25  -4.2 0.25 /
\plot -2.8 0.25  -2.2 0.25 /
\plot -1.8 0.25  -1.2 0.25 /
\plot -0.8 0.25  -0.2 0.25 /
\plot 0.2 0.25  0.8 0.25 /
\plot 1.2 0.25  1.8 0.25 /
\plot 2.2 0.25  2.7 0.25 /
\ellipticalarc axes ratio 2:1 180 degrees from -4.65 0.25  center 
at -4.65 0 
\plot -4.65 -0.25  -3.2 -0.25 /
\plot -2.8 -0.25  -2.2 -0.25 /
\plot -1.8 -0.25  -1.2 -0.25 /
\plot -0.8 -0.25  -0.2 -0.25 /
\plot 0.2 -0.25  0.8 -0.25 /
\plot 1.2 -0.25  1.8 -0.25 /
\plot 2.2 -0.25  2.7 -0.25 /
\setquadratic
\plot  2.7 0.25  2.95 0.45  3 0.75 /
\plot  2.7 -0.25  2.95 -0.45  3 -0.75 /
\endpicture
$$

If $M$ is a finite dimensional $U=U_h\fg$ module
and $a\in \End_U(M)$ the {\it quantum trace} of $a$ on $M$ (see
[LR \S3] and [CP Def. 4.2.9]) is 
the trace of the action of $e^{h\rho}a$ on $M$,
$$\tr_q(a)=\Tr(e^{h\rho}a,M),
\qquad\hbox{and}\qquad
\dim_q(M) = \tr_q(\id_M)=\Tr(e^{h\rho},M)
\formula$$
is the {\it quantum dimension} of $M$.
The first step of the standard argument for proving Weyl's dimension
formula [B-tD, VI Lemma 1.19] 
shows that the quantum dimension of the finite dimensional
irreducible $U_h\fg$-module $L(\mu)$ is
$$\dim_q(L(\mu)) = \Tr(e^{h\rho}, L(\mu))
=\prod_{\alpha>0} 
{e^{{h\over2}\langle \mu+\rho, \alpha^\vee\rangle}
-e^{-{h\over2}\langle \mu+\rho, \alpha^\vee\rangle}
\over
e^{{h\over2}\langle \rho, \alpha^\vee\rangle}
-e^{-{h\over2}\langle \rho, \alpha^\vee\rangle} }
=\prod_{\alpha>0} {[\langle\mu+\rho,\alpha^\vee]\over
[\langle\rho,\alpha^\vee\rangle]},
\formula$$
where $q=e^{h/2}$ and $[d] = (q^d-q^{-d})/(q-q^{-1})$ for a positive integer $d$.

\thm  Let $\mu,\nu\in P^+$ be dominant integral weights.
Let $M=L(\mu)$ and $V=L(\nu)$ and let $\Phi_k$ be
the representation of $\tilde {\cal B}_k$ defined in Proposition 3.5.  Then
the functions
$$\matrix{
\mt_k\colon &\tilde {\cal B}_k &\longrightarrow &\CC \cr
&b &\longmapsto &\displaystyle{
\tr_q(\Phi_k(b))\over
\dim_q(M)\dim_q(V)^k} \cr
}$$
form a Markov trace on the affine braid group with
parameters
$$z={q^{\langle\nu,\nu+2\rho\rangle}\over \dim_q(V)}
\quad\hbox{and}\quad
Q_r= \sum_\lambda
q^{r(\langle \lambda,\lambda+2\rho\rangle
-\langle \mu,\mu+2\rho\rangle
-\langle \nu,\nu+2\rho\rangle)}
{\dim_q(L(\lambda))c_{\mu\nu}^\lambda\over
\dim_q(L(\mu))\dim_q(L(\nu))},$$
where the positive integers $c_{\mu\nu}^\lambda$
and the sum in the expression for $Q_r$ are as in the
tensor product decomposition
$$L(\mu)\otimes L(\nu)
= \bigoplus_\lambda L(\lambda)^{\oplus c_{\mu\nu}^\lambda}.$$
\pf
The fact that $\hbox{\sl mt}_k$ as defined in the statement of the Theorem satisfies
(1)-(4) in the definition of a Markov trace follows exactly as in 
[LR] Theorem 3.10c.  The formula for the parameter $z$
is derived in [LR, (3.9) and Thm. 3.10(2)].

It remains to check (5).  The proof is a combination
of the argument used in [Or] Theorem 5.3 and the argument in the proof of
[LR] Theorem 3.10c.
Let $\varepsilon_k\colon \End_U(M\otimes V^{\otimes k})
\to \End_U(M\otimes V^{\otimes(k-1)})$ be given by
$$\varepsilon_k(z) = (\id_{M\otimes V^{\otimes k}}\otimes \check e)(z\otimes \id)
\qquad \hbox{where}\qquad
\matrix{
\check e\colon &V\otimes V^* &\longrightarrow &\CC\cr
&x\otimes \phi &\longmapsto &\dim_q(V)^{-1}\phi(e^{h\rho} x).
\cr}
\formula
$$
If $V$ is simple then $\check e$ is the unique $U_h\fg$-invariant projection
onto the invariants in $V\otimes V^*$.  Pictorially,
$$\varepsilon_k\left(
\beginpicture
\setcoordinatesystem units <.5cm,.5cm>         
\setplotarea x from -2 to 2, y from -2 to 2    
\put{${}_1$} at -.8 1.8
\put{$\ldots$} at 0.2 1.8
\put{${}_k$} at 1.2 1.8
\put{$z$} at 0 0
\plot  -1.5 1   1.5 1 /
\plot  -1.5 -1   1.5 -1 /
\plot  -1.5 -1   -1.5 1 /
\plot  1.5 -1   1.5 1 /
\plot  -.8 1.5    -.8  1 /
\plot  -.4 1.5    -.4  1 /
\plot    0 1.5     0 1 /
\plot  .4 1.5    .4  1 /
\plot  .8 1.5    .8 1 /
\plot  1.2 1.5    1.2 1 /
\plot  -.8 -1.5    -.8  -1 /
\plot  -.4 -1.5    -.4  -1 /
\plot    0 -1.5     0 -1 /
\plot  .4 -1.5    .4  -1 /
\plot  .8 -1.5    .8 -1 /
\plot  1.2 -1.5    1.2 -1 /
\linethickness=1pt
\putrule  from -1.2 1.5 to  -1.2  1 
\putrule  from  -1.2 -1.5  to  -1.2  -1 
\endpicture
\right) = 
\beginpicture
\setcoordinatesystem units <.5cm,.5cm>         
\setplotarea x from -2 to 2, y from -2 to 2    
\put{$z$} at 0 0
\plot  -1.5 1   1.5 1 /
\plot  -1.5 -1   1.5 -1 /
\plot  -1.5 -1   -1.5 1 /
\plot  1.5 -1   1.5 1 /
\plot  -.8 1.5    -.8  1 /
\plot  -.4 1.5    -.4  1 /
\plot    0 1.5     0 1 /
\plot  .4 1.5    .4  1 /
\plot  .8 1.5    .8 1 /
\plot  1.2 1.5    1.2 1 /
\plot  -.8 -1.5    -.8  -1 /
\plot  -.4 -1.5    -.4  -1 /
\plot    0 -1.5     0 -1 /
\plot  .4 -1.5    .4  -1 /
\plot  .8 -1.5    .8 -1 /
\plot  1.2 -1.5    1.2 -1 /
\plot  1.9 1.5    1.9 -1.5 /
\setquadratic
\plot  1.2 1.5  1.55 1.85  1.9 1.5 /
\plot  1.2 -1.5 1.55 -1.85  1.9 -1.5 /
\linethickness=1pt
\putrule  from -1.2 1.5 to  -1.2  1 
\putrule  from  -1.2 -1.5  to  -1.2  -1 
\endpicture
= ~
\beginpicture
\setcoordinatesystem units <.5cm,.5cm>         
\setplotarea x from -2 to 2, y from -2 to 2    
\put{${}_1$} at -.8 1.9
\put{$\ldots$} at -0.2 1.9
\put{${}_{k-1}$} at .9 1.9
\put{$\varepsilon_k(z)$} at 0 0
\plot  -1.5 1   1.2 1 /
\plot  -1.5 -1   1.2 -1 /
\plot  -1.5 -1   -1.5 1 /
\plot  1.2 -1   1.2 1 /
\plot  -.8 1.5    -.8  1 /
\plot  -.4 1.5    -.4  1 /
\plot    0 1.5     0 1 /
\plot  .4 1.5    .4  1 /
\plot  .8 1.5    .8 1 /
\plot  -.8 -1.5    -.8  -1 /
\plot  -.4 -1.5    -.4  -1 /
\plot    0 -1.5     0 -1 /
\plot  .4 -1.5    .4  -1 /
\plot  .8 -1.5    .8 -1 /
\linethickness=1pt
\putrule  from -1.2 1.5 to  -1.2  1 
\putrule  from  -1.2 -1.5  to  -1.2  -1 
\endpicture
.$$
The argument of [LR] Theorem 3.10b shows that
$$\mt_k(b)=\mt_{k-1}(\varepsilon_{k-1}(b)),
\qquad\hbox{if $b\in \tilde {\cal B}_k$.}
\formula$$
Since $\varepsilon_1((X^{\varepsilon_1})^r)$ is a $U_h\fg$-module homomorphism
from $M$ to $M$ and, since $M$ is simple, Schur's lemma implies that
$$
\beginpicture
\setcoordinatesystem units <.5cm,.5cm>         
\setplotarea x from -4 to 1.5, y from -2 to 2    
\put{$\hbox{$r$ loops\ }\Bigg\{$} at -3.8 .4   %
\put{$\bullet$} at  0 2      
\put{$\bullet$} at  0 -1.5          
\plot -1.5 2.2 -1.5 1.12 /
\plot -1.5 .88 -1.5 0.12 /
\plot -1.5 -0.12 -1.5 -0.88 /
\plot -1.5 -1.12 -1.5 -1.75 /
\plot -1.25 2.2 -1.25 1.12 /
\plot -1.25 .88 -1.25 0.12 /
\plot -1.25 -0.12 -1.25 -0.88 /
\plot -1.25 -1.12 -1.25 -1.75 /
\ellipticalarc axes ratio 1:1 360 degrees from -1.5 2.2 center 
at -1.375 2.2
\put{$*$} at -1.375 2.2  
\ellipticalarc axes ratio 1:1 180 degrees from -1.5 -1.75 center 
at -1.375 -1.75 
\plot  1 2   1 -1.5 /
\setlinear
\plot -0.3 1.5  -1.1 1.5 /
\ellipticalarc axes ratio 2:1 180 degrees from -1.65 1.5  center 
at -1.65 1.25 
\plot -1.65 1  -1.1 1 /
\ellipticalarc axes ratio 2:1 -180 degrees from -1.1 1  center 
at -1.1 .75 
\ellipticalarc axes ratio 2:1 180 degrees from -1.65 .5  center 
at -1.65 .25 
\ellipticalarc axes ratio 2:1 -180 degrees from -1.1 0  center 
at -1.1 -0.25 
\plot -1.65 0  -1.1 0 /
\ellipticalarc axes ratio 2:1 180 degrees from -1.65 -.5  center 
at -1.65 -.75 
\plot -1.65 -1  -0.3 -1 /
\setquadratic
\plot  -0.3 1.5  -0.05 1.7  -0 2 /
\plot  -0.3 -1  -0.05 -1.2  -0 -1.5 /
\ellipticalarc axes ratio 1:1 180 degrees from 1 2 center 
at 0.5 2
\ellipticalarc axes ratio 1:1 180 degrees from 0 -1.5 center 
at 0.5 -1.5
\endpicture
~ = \varepsilon_1((X^{\varepsilon_1})^r)=\xi\cdot\id_M,
\qquad\hbox{for some $\xi\in \CC$.}$$
Let
$\tilde R_i = \id_V^{\otimes i}\otimes \check R_{VM}
\otimes \id_V^{(k+1)-i}$.
Then
$(\tilde X^{\varepsilon_{k+1}})^r
=
(\tilde R_1\cdots \tilde R_k)^{-1}
(X^{\varepsilon_1})^r(\tilde R_1\cdots \tilde R_k)$
and
$$\eqalign{
\mt_{k+1}(b(\tilde X^{\varepsilon_{k+1}})^r)
&=\mt_k(\varepsilon_k(b(\tilde X^{\varepsilon_{k+1}})^r)) \cr
&=\mt_k(\varepsilon_k(b(\tilde R_1\cdots\tilde R_k)^{-1}
(X^{\varepsilon_1})^r(\tilde R_1\cdots \tilde R_k)) \cr
&=\mt_k(b(\tilde R_1\cdots\tilde R_k)^{-1}
\varepsilon_1((X^{\varepsilon_1})^r)
(\tilde R_1\cdots \tilde R_k)) \cr
&=\mt_k(b(\tilde R_k\cdots \tilde R_1)^{-1}\xi\cdot\id_M\tilde R_1\cdots \tilde R_k)
=
\xi\cdot\mt_k(b). \cr
}$$
This last calculation is more palatable in a pictorial format,
$$\eqalign{
\mt_{k+1}\left(
\beginpicture
\setcoordinatesystem units <.5cm,.5cm>         
\setplotarea x from -1.8 to 2.2, y from -2.4 to 1.5   
\put{${}_1$} at -.8 1.8
\put{$\ldots$} at 0.2 1.8
\put{${}_k$} at 1.2 1.8
\put{$b$} at 0 0.5
\plot  -1.5 1   1.5 1 /
\plot  -1.5 0   1.5 0 /
\plot  -1.5 0   -1.5 1 /
\plot  1.5 0   1.5 1 /
\plot  -.8 1.5    -.8  1 /
\plot  -.4 1.5    -.4  1 /
\plot    0 1.5     0 1 /
\plot  .4 1.5    .4  1 /
\plot  .8 1.5    .8 1 /
\plot  1.2 1.5    1.2 1 /
\plot  -.8 -0.5    -.8  -0 /
\plot  -.4 -0.5    -.4  -0 /
\plot    0 -0.5     0 -0 /
\plot  .4 -0.5    .4  -0 /
\plot  .8 -0.5    .8 -0 /
\plot  1.2 -0.5  1.2 -0 /
\plot  1.8 1.5    1.8 -.5 /
\plot  1.8 -1.5    1.8 -2 /
\put{$(\tilde X^{\varepsilon_{k+1}})^r$} at 0.2 -1
\plot -1.5 -.5   2 -.5 /
\plot -1.5 -1.5   2 -1.5 /
\plot -1.5 -.5   -1.5 -1.5 /
\plot 2 -.5     2 -1.5 /
\plot  -.8 -1.5    -.8  -2 /
\plot  -.4 -1.5    -.4  -2 /
\plot    0 -1.5     0 -2 /
\plot  .4 -1.5    .4  -2 /
\plot  .8 -1.5    .8 -2 /
\plot  1.2 -1.5    1.2 -2 /
\linethickness=1pt
\putrule  from -1.2 1.5 to  -1.2  1 
\putrule  from  -1.2 -0.5  to  -1.2  -0 
\putrule  from  -1.2 -1.5  to  -1.2  -2 
\endpicture
\right)
&=\mt_{k+1}\left(
\beginpicture
\setcoordinatesystem units <.5cm,.5cm>         
\setplotarea x from -1.8 to 4.2, y from -2 to 3   
\put{${}_1$} at -.8 2.3
\put{$\ldots$} at 0.2 2.3
\put{${}_k$} at 1.2 2.3
\put{$b$} at 0 1
\plot  -1.5 1.5   1.5 1.5 /
\plot  -1.5 0.5   1.5 0.5 /
\plot  -1.5 0.5   -1.5 1.5 /
\plot  1.5 0.5   1.5 1.5 /
\plot  -.8 2    -.8  1.5 /
\plot  -.4 2    -.4  1.5 /
\plot    0 2     0 1.5 /
\plot  .4 2    .4  1.5 /
\plot  .8 2    .8 1.5 /
\plot  1.2 2    1.2 1.5 /
\plot  -.8 0.5    -.8  -2.5 /
\plot  -.4 0.5    -.4  -2.5 /
\plot    0 0.5     0 -2.5 /
\plot  .4 0.5    .4  -2.5 /
\plot  .8 0.5    .8 -2.5 /
\plot  1.2 0.5  1.2 -2.5 /
\plot  3.5 2    3.5 -0.5 /
\plot  3.5 -1.5    3.5 -2.5 /
\put{$(\tilde X^{\varepsilon_1})^r$} at 2.8 -1
\plot 1.5 -0.5   4 -0.5 /
\plot 1.5 -1.5   4 -1.5 /
\plot 1.5 -0.5   1.5 -1.5 /
\plot 4 -0.5    4 -1.5 /
\linethickness=1pt
\putrule  from -1.2 2 to  -1.2  1.5 
\putrule from -0.7 0   to  -0.5  0
\putrule from -0.3 0   to  -0.1  0
\putrule from 0.1 0   to  0.3  0
\putrule from 0.5 0   to  0.7  0
\putrule from 0.9 0   to  1.1  0
\putrule from 1.3 0   to  1.7  0
\putrule from -0.7 -2   to  -0.5  -2
\putrule from -0.3 -2   to  -0.1  -2
\putrule from 0.1 -2   to  0.3  -2
\putrule from 0.5 -2   to  0.7  -2
\putrule from 0.9 -2   to  1.1  -2
\putrule from 1.3 -2   to  1.7  -2
\setquadratic
\plot  1.7 -2  1.95 -1.8  2 -1.5 /
\plot  1.7 0  1.95 -0.2  2 -0.5 /
\plot  -0.9 0  -1.15 0.2  -1.2 0.5 /
\plot  -0.9 -2  -1.15 -2.2  -1.2 -2.5 /
\endpicture
\right) 
=\mt_k\left(
\beginpicture
\setcoordinatesystem units <.5cm,.5cm>         
\setplotarea x from -2 to 4.5, y from -2 to 3   
\put{${}_1$} at -.8 2.3
\put{$\ldots$} at 0.2 2.3
\put{${}_k$} at 1.2 2.3
\put{$b$} at 0 1
\plot  -1.5 1.5   1.5 1.5 /
\plot  -1.5 0.5   1.5 0.5 /
\plot  -1.5 0.5   -1.5 1.5 /
\plot  1.5 0.5   1.5 1.5 /
\plot  -.8 2    -.8  1.5 /
\plot  -.4 2    -.4  1.5 /
\plot    0 2     0 1.5 /
\plot  .4 2    .4  1.5 /
\plot  .8 2    .8 1.5 /
\plot  1.2 2    1.2 1.5 /
\plot  -.8 0.5    -.8  -2.5 /
\plot  -.4 0.5    -.4  -2.5 /
\plot    0 0.5     0 -2.5 /
\plot  .4 0.5    .4  -2.5 /
\plot  .8 0.5    .8 -2.5 /
\plot  1.2 0.5  1.2 -2.5 /
\put{$(\tilde X^{\varepsilon_1})^r$} at 2.8 -1
\plot 1.5 -0.5   4 -0.5 /
\plot 1.5 -1.5   4 -1.5 /
\plot 1.5 -0.5   1.5 -1.5 /
\plot 4 -0.5    4 -1.5 /
\plot 4.4 -0.5   4.4 -1.5 /
\linethickness=1pt
\putrule  from -1.2 2 to  -1.2  1.5 
\putrule from -0.7 0   to  -0.5  0
\putrule from -0.3 0   to  -0.1  0
\putrule from 0.1 0   to  0.3  0
\putrule from 0.5 0   to  0.7  0
\putrule from 0.9 0   to  1.1  0
\putrule from 1.3 0   to  1.7  0
\putrule from -0.7 -2   to  -0.5  -2
\putrule from -0.3 -2   to  -0.1  -2
\putrule from 0.1 -2   to  0.3  -2
\putrule from 0.5 -2   to  0.7  -2
\putrule from 0.9 -2   to  1.1  -2
\putrule from 1.3 -2   to  1.7  -2
\setquadratic
\plot  1.7 -2  1.95 -1.8  2 -1.5 /
\plot  1.7 0  1.95 -0.2  2 -0.5 /
\plot  -0.9 0  -1.15 0.2  -1.2 0.5 /
\plot  -0.9 -2  -1.15 -2.2  -1.2 -2.5 /
\plot  3.7 -0.5  4.05 -0.15  4.4 -0.5 /
\plot  3.7 -1.5  4.05 -1.85  4.4 -1.5 /
\endpicture
\right) \cr
&=\xi\cdot\mt_k\left(
\beginpicture
\setcoordinatesystem units <.5cm,.5cm>         
\setplotarea x from -2 to 3, y from -2 to 3   
\put{${}_1$} at -.8 2.3
\put{$\ldots$} at 0.2 2.3
\put{${}_k$} at 1.2 2.3
\put{$b$} at 0 1
\plot  -1.5 1.5   1.5 1.5 /
\plot  -1.5 0.5   1.5 0.5 /
\plot  -1.5 0.5   -1.5 1.5 /
\plot  1.5 0.5   1.5 1.5 /
\plot  -.8 2    -.8  1.5 /
\plot  -.4 2    -.4  1.5 /
\plot    0 2     0 1.5 /
\plot  .4 2    .4  1.5 /
\plot  .8 2    .8 1.5 /
\plot  1.2 2    1.2 1.5 /
\plot  -.8 0.5    -.8  -2.5 /
\plot  -.4 0.5    -.4  -2.5 /
\plot    0 0.5     0 -2.5 /
\plot  .4 0.5    .4  -2.5 /
\plot  .8 0.5    .8 -2.5 /
\plot  1.2 0.5  1.2 -2.5 /
\linethickness=1pt
\putrule  from -1.2 2 to  -1.2  1.5 
\putrule  from 2 -0.5 to  2 -1.5 
\putrule from -0.7 0   to  -0.5  0
\putrule from -0.3 0   to  -0.1  0
\putrule from 0.1 0   to  0.3  0
\putrule from 0.5 0   to  0.7  0
\putrule from 0.9 0   to  1.1  0
\putrule from 1.3 0   to  1.7  0
\putrule from -0.7 -2   to  -0.5  -2
\putrule from -0.3 -2   to  -0.1  -2
\putrule from 0.1 -2   to  0.3  -2
\putrule from 0.5 -2   to  0.7  -2
\putrule from 0.9 -2   to  1.1  -2
\putrule from 1.3 -2   to  1.7  -2
\setquadratic
\plot  1.7 -2  1.95 -1.8  2 -1.5 /
\plot  1.7 0  1.95 -0.2  2 -0.5 /
\plot  -0.9 0  -1.15 0.2  -1.2 0.5 /
\plot  -0.9 -2  -1.15 -2.2  -1.2 -2.5 /
\endpicture
\right) 
=\xi\cdot \mt_k\left(
\beginpicture
\setcoordinatesystem units <.5cm,.5cm>         
\setplotarea x from -1.5 to 1.5, y from 0.5 to 2   
\put{${}_1$} at -.8 1.3
\put{$\ldots$} at 0.2 1.3
\put{${}_k$} at 1.2 1.3
\put{$b$} at 0 0
\plot  -1.5 0.5   1.5 0.5 /
\plot  -1.5 -0.5   1.5 -0.5 /
\plot  -1.5 -0.5   -1.5 0.5 /
\plot  1.5 -0.5   1.5 0.5 /
\plot  -.8 1    -.8  0.5 /
\plot  -.4 1    -.4  0.5 /
\plot    0 1     0 0.5 /
\plot  .4 1    .4  0.5 /
\plot  .8 1    .8 0.5 /
\plot  1.2 1    1.2 0.5 /
\plot  -.8 -1    -.8  -0.5 /
\plot  -.4 -1    -.4  -0.5 /
\plot    0 -1     0 -0.5 /
\plot  .4 -1    .4  -0.5 /
\plot  .8 -1    .8 -0.5 /
\plot  1.2 -1  1.2 -0.5 /
\linethickness=1pt
\putrule  from -1.2 1 to  -1.2  0.5 
\putrule  from  -1.2 -0.5  to  -1.2  -1 
\endpicture
\right). \cr
}$$
It remains to calculate the constant $\xi$.  By (2.14),
$$(X^{\varepsilon_1})^r =(\check R_0^2)^r
=\left(
\sum_\lambda q^{c(\lambda)}P_{\mu\nu}^\lambda\right)^r
= \sum_\lambda q^{rc(\lambda)}P_{\mu\nu}^\lambda,$$
where $c(\lambda)=\langle \lambda,\lambda+2\rho\rangle
-\langle \mu,\mu+2\rho\rangle-\langle \nu,\nu+2\rho\rangle$
and $P_{\mu\nu}^\lambda$ is the projection onto the
$L(\lambda)^{\oplus c_{\mu\nu}^\lambda}$ component in the decomposition
of $M\otimes V = L(\mu)\otimes L(\nu)$.  Thus
$$\eqalign{
\xi&=\mt_0(\xi\cdot\id_M)
=\mt_1((X^{\varepsilon_1})^r)
={1\over \dim_q(M)\dim_q(V)}\tr_q\left(
\sum_\lambda q^{rc(\lambda)}P_{\mu\nu}^\lambda\right)  \cr
&= {1\over \dim_q(M)\dim_q(V)}\tr_q\left(
\sum_\lambda q^{rc(\lambda)}c_{\mu\nu}^\lambda \id_{L(\lambda)} \right)  \cr
&= 
\sum_\lambda q^{rc(\lambda)}c_{\mu\nu}^\lambda 
{\dim_q(L(\lambda))\over \dim_q(L(\mu))\dim_q(L(\nu))}. 
\qquad\hbox{\qed}\cr
}$$

\remark
There is another formula [TW, Lemma (3.51)] for the constant
$Q_1$ in Theorem 5.4, namely,
$$
Q_1 = {\sum_{w\in W} 
(-1)^{\ell(w)}q^{\langle \nu+\rho, w(\mu+\rho)\rangle} 
\over
\sum_{w\in W} (-1)^{\ell(w)}q^{\langle \nu+\rho,w\rho\rangle}
},\formula$$
where $W$ is the Weyl group of $\fg$.

\medskip
Let $\mt_k$ be as in Theorem 5.4 and let $\tilde {\cal Z}_k
=\End_U(M\otimes V^{\otimes k})$.  Then $\mt_k$ is the restriction
of the linear functional
$$\matrix{
\mt_k\colon &\tilde {\cal Z}_k &\longrightarrow &\CC\cr
&a &\longmapsto &\displaystyle{
\tr_q(a)\over \dim_q(M)\dim_q(V^{\otimes k}) }
\cr}
\formula$$
to $\Phi_k(\tilde {\cal B}_k)$.  Since $M\otimes V^{\otimes k}$
is a finite dimensional semisimple module $\tilde {\cal Z}_k$
is a finite dimensional semisimple algebra.  The {\it weights}
of the Markov trace $\mt$ are the constants $t_{\lambda/\mu}$
defined by
$$\mt_k = \sum_\lambda t_{\lambda/\mu} 
\chi^{\lambda/\mu}_{\tilde {\cal Z}_k},
\formula$$
where $\chi^\lambda_{\tilde {\cal Z}_k}$ are the irreducible characters
of $\tilde {\cal Z}_k$.

\thm  Let $M=L(\mu)$ and $V=L(\nu)$ be finite dimensional irreducible 
$U_h\fg$-modules.
The weights of the Markov trace on the affine braid group
defined in Theorem 5.4 are
$$t_{\lambda/\mu} = {\dim_q(L(\lambda))\over
\dim_q(L(\mu))\dim_q(V)^k}.$$
\pf
Since $M\otimes V^{\otimes k}$ is finite dimensional and semisimple
the algebra $\tilde {\cal Z}_k=\End_U(M\otimes V^{\otimes k})$
is a finite dimensional semisimple algebra.  Schur's lemma
can be used to show that, as a $(U_h\fg, \tilde {\cal Z}_k)$
bimodule,
$$M\otimes V^{\otimes k} \cong
\bigoplus_\lambda L(\lambda)\otimes {\cal L}^{\lambda/\mu},
\formula$$
where the ${\cal L}^{\lambda/\mu}$ are the irreducible
$\tilde {\cal Z}_k$ modules.  
In the notation of (4.1),  
${\cal L}^{\lambda/\mu}=F_\lambda(L(\mu))$
and $\chi^{\lambda/\mu}_{\tilde {\cal Z}_k}$ is the character of
${\cal L}^{\lambda/\mu}$.  Taking the quantum trace on both sides
of (5.12) gives
$$\tr_q(a) = \Tr(e^{-h\rho}a) = \sum_\lambda
\Tr(e^{-h\rho}, L(\lambda)) \chi^{\lambda/\mu}_{\tilde {\cal Z}_k}(a)
=\sum_\lambda
\tr_q(L(\lambda))\chi^{\lambda/\mu}_{\tilde {\cal Z}_k}(a).
$$
The result follows by dividing both sides by $\dim_q(L(\mu))\dim_q(V)^k$.
\endpf

\section 6.  Examples

\subsection Affine and cyclotomic Hecke algebras

Let $q\in \CC^*$.  The {\it affine Hecke algebra} $\tilde H_k$
is the quotient of the group algebra $\CC \tilde {\cal B}_k$ of
the affine braid group by the relations
$$T_i^2=(q-q^{-1})T_i+1, \qquad\qquad 1\le i\le k-1.
\formula$$
The affine Hecke algebra $\tilde H_k$ is an infinite dimensional
algebra with a very interesting representation theory (see
[KL] and [CG]).  With $X$ as in (3.4) the subalgebra 
$$\CC[X] = \CC[X^{\pm \varepsilon_1},\ldots,X^{\pm \varepsilon_k}]
=\hbox{span\ }\{ X^\lambda\ |\ \lambda\in L\}$$ 
is a commutative subalgebra of $\tilde H_k$.  It is a theorem
of Bernstein and Zelevinsky (see [RR, Theorem 4.12]) that the center 
of $\tilde H_k$ is the ring of symmetric (Laurent) polynomials in 
$X^{\pm \varepsilon_1},\ldots,X^{\pm \varepsilon_k}$,
$$Z(\tilde H_k) = \CC[X]^{S_k}
=\CC[X^{\pm \varepsilon_1},\ldots,X^{\pm \varepsilon_k}]^{S_k}.$$
If $w\in S_k$ define $T_w=T_{i_1}\cdots T_{i_p}$ if
$w=s_{i_1}\cdots s_{i_p}$ is a reduced word for $w$ in terms
of the generating reflections $s_i=(i,i+1)$, $1\le i\le k-1$,
of $S_k$.  Then, with $X^\lambda$ as in (3.4)
$$\{X^\lambda T_w\ |\ \lambda\in L, w\in S_k\}
\quad\hbox{is a basis of $\tilde H_k$.}$$

Let $u_1,\ldots, u_r\in \CC$.  The {\it cyclotomic Hecke algebra}
$H_{r,1,n}$ with parameters $u_1,\ldots, u_r,q$ is the quotient of 
the affine Hecke algebra by the relation
$$(X^{\varepsilon_1}-u_1)(X^{\varepsilon_1}-u_2)\cdots
(X^{\varepsilon_1}-u_r)=0.
\formula$$
The algebra $H_{r,1,n}$ is a deformation of the group algebra
of the complex reflection group $G(r,1,n)=(\ZZ/r\ZZ)\wr S_n$ and 
is of dimension $r^nn!$.  It was introduced by Ariki and Koike [AK]
and its representations and its connection to the affine Hecke algebra
have been well studied ([Ar],[AK],[Gk]).

\vfill\eject

\subsection The affine and cyclotomic BMW algebras

Fix $q,z\in \CC^*$ and an infinite number of values
$Q_1,Q_2,\ldots$ in $\CC$.  The {\it affine BMW 
(Birman-Murakami-Wenzl) algebra} $\tilde {\cal Z}_k$ is the quotient 
of the group algebra $\CC \tilde {\cal B}_k$ of the affine braid 
group by the relations
\smallskip
\itemitem{(\global\advance\resultno by 1
\the\sectno.\the\resultno a)}\ \ $(T_i-z^{-1})(T_i+q^{-1})(T_i-q)=0$, 
\smallskip
\itemitem{(\the\sectno.\the\resultno b)}\ \ 
$E_iT_i^{\pm 1} = T_i^{\pm 1}E_i = z^{\mp1}E_i$, 
\smallskip
\itemitem{(\the\sectno.\the\resultno c)}\ \ 
$E_iT_{i-1}^{\pm 1}E_i=z^{\pm1}E_i$ \quad and \quad
$E_iT_{i+1}^{\pm 1}E_i=z^{\pm1}E_i$,
\smallskip
\itemitem{(\the\sectno.\the\resultno d)}\ \ 
$E_1(X^{\varepsilon_1})^rE_1 =Q_rE_1$, 
\smallskip
\itemitem{(\the\sectno.\the\resultno e)}\ \ 
$E_1X^{\varepsilon_1}T_1X^{\varepsilon_1}=z^{-1}E_1$, 
\smallskip\noindent
where the $E_i$, $1\le i\le k-1$, are defined by the equations
$$
{T_i-T_i^{-1}\over q-q^{-1}} = 1-E_i, \qquad 1\le i\le k-1.
\formula$$
It follows that 
$$E_i^2 = xE_i
\qquad\hbox{where}\qquad x = {z-z^{-1}\over q-q^{-1}}+1.
\formula$$
The classical {\it BMW algebra} is the subalgebra ${\cal Z}_k$
of the affine BMW algebra which is generated by $T_1,\ldots, T_{k-1}$, and 
$E_1,\ldots, E_{k-1}$.  
Fix $u_1,\ldots, u_r\in \CC$. 
The {\it cyclotomic BMW} algebra ${\cal Z}_{r,1,k}$ is the quotient
of the affine BMW algebra by the relation
$$(X^{\varepsilon_1}-u_1)(X^{\varepsilon_1}-u_2)\cdots (X^{\varepsilon_1}-u_r)=0.
\formula$$

Although the affine BMW algebras have been ``in the air''
for some time we are not aware of any existing literature.  
The ``degenerate'' version of these algebras
were defined by Nazarov [Nz] who called them
``degenerate affine Wenzl algebras''.  The relation between his algebras
and the affine BMW algebras $\tilde {\cal Z}_k$ is analogous to
the relation between the graded Hecke algebras (sometimes called
the degenerate affine Hecke algebras) and the affine Hecke
algebras (see [Lu]).
The cyclotomic BMW algebras have been defined and studied by [H\"a1-2].  
They are quotients of the 
affine BMW algebras in the same way that cyclotomic Hecke algebras
are quotients of affine Hecke algebras.
The classical BMW algebras ${\cal Z}_k={\cal Z}_{1,1,k}$ have been
studied in [Wz2], [HR], [Mu], [LR] and many other works.

Elements of the affine BMW algebra can be viewed as linear combinations
of affine tangles.  An affine tangle has $k$ strands and a flagpole
just as in the case of an affine braid, but there is no restriction
that a strand must connect an upper vertex to a lower vertex.
Let $X^{\varepsilon_1}$ and $T_i$ be the affine braids given in (3.1)
and let
$$E_i = 
\beginpicture
\setcoordinatesystem units <.5cm,.5cm>         
\setplotarea x from -5 to 3.5, y from -1 to 1    
\put{$\bullet$} at -3 0.75      %
\put{$\bullet$} at -2 0.75      %
\put{$\bullet$} at -1 0.75      %
\put{$\bullet$} at  0 0.75      
\put{$\bullet$} at  1 0.75      %
\put{$\bullet$} at  2 0.75      %
\put{$\bullet$} at  3 0.75      %
\put{$\bullet$} at -3 -0.75          %
\put{$\bullet$} at -2 -0.75          %
\put{$\bullet$} at -1 -0.75          %
\put{$\bullet$} at  0 -0.75          
\put{$\bullet$} at  1 -0.75          %
\put{$\bullet$} at  2 -0.75          %
\put{$\bullet$} at  3 -0.75          %
\plot -4.5 1.25 -4.5 -1.25 /
\plot -4.25 1.25 -4.25 -1.25 /
\ellipticalarc axes ratio 1:1 360 degrees from -4.5 1.25 center 
at -4.375 1.25
\put{$*$} at -4.375 1.25  
\ellipticalarc axes ratio 1:1 180 degrees from -4.5 -1.25 center 
at -4.375 -1.25 
\plot -3 0.75  -3 -0.75 /
\plot -2 0.75  -2 -0.75 /
\plot -1 0.75  -1 -0.75 /
\plot  2 0.75   2 -0.75 /
\plot  3 0.75   3 -0.75 /
\plot 0.3 0.25  0.7 0.25 /
\plot 0.3 -0.25    0.7 -0.25   /
\setquadratic
\plot  0.7 0.25  0.95 0.45  1 0.75 /
\plot  0 0.75  0.05 0.45  0.3 0.25 /
\plot  0.7 -0.25  0.95 -0.45  1 -0.75 /
\plot  0 -0.75  0.05 -0.45  0.3 -0.25 /
\endpicture
\formula$$
Then $\tilde {\cal Z}_k$ is the algebra of linear combinations of
tangles generated by $X^{\varepsilon_1}, T_1,\ldots, T_{k-1},
E_1,\ldots, E_{k-1}$ and the relations in (6.3) expressed in the 
form
$$
\beginpicture
\setcoordinatesystem units <.5cm,.5cm>         
\setplotarea x from -1 to 1, y from -1 to 1    
\setquadratic
\plot  -.5 -.75  -.45 -0.45  -0.1 -0.1 /
\plot  .1 .15  .45 .45  .5 .75 /
\plot -.5 .75  -.45 .45  0 0  .45 -.45  .5 -.75 /
\endpicture
-
\beginpicture
\setcoordinatesystem units <.5cm,.5cm>         
\setplotarea x from -1 to 1, y from -1 to 1    
\setquadratic
\plot  -.5 .75  -.45 0.45  -0.1 0.1 /
\plot  .1 -.15  .45 -.45  .5 -.75 /
\plot -.5 -.75  -.45 -.45  0 0  .45 .45  .5 .75 /
\endpicture
~=~
(q-q^{-1})\left(
\beginpicture
\setcoordinatesystem units <.5cm,.5cm>         
\setplotarea x from -1 to 1, y from -1 to 1    
\plot -0.5 .75   -0.5 -.75 /
\plot  0.5 .75    0.5 -.75 /
\endpicture
-
\beginpicture
\setcoordinatesystem units <.5cm,.5cm>         
\setplotarea x from -1 to 1, y from -1 to 1    
\plot -0.2 0.25  0.2 0.25 /
\plot -0.2 -0.25    0.2 -0.25   /
\setquadratic
\plot  0.2 0.25  0.45 0.45  0.5 0.75 /
\plot  -.5 0.75  -.45 0.45  -0.2 0.25 /
\plot  0.2 -0.25  0.45 -0.45  0.5 -0.75 /
\plot  -.5 -0.75  -.45 -0.45  -0.2 -0.25 /
\endpicture
\right)
\formula$$
$$
\beginpicture
\setcoordinatesystem units <.5cm,.5cm>         
\setplotarea x from -1 to 2, y from -1 to 1    
\plot 1.5 .75  1.5 -0.75 /
\plot 0.8 -1.25  1.2 -1.25 /
\plot 0.8 1.25    1.2 1.25   /
\setquadratic
\plot  1.2 -1.25  1.45 -1.05  1.5 -0.75 /
\plot  0.5 -0.75  .55 -1.05  0.8 -1.25 /
\plot  1.2 1.25  1.45 1.05  1.5 0.75 /
\plot  .5 0.75  .55 1.05  0.8 1.25 /
\setquadratic
\plot  -.5 -.75  -.45 -0.45  -0.1 -0.1 /
\plot  .1 .15  .45 .45  .5 .75 /
\plot -.5 .75  -.45 .45  0 0  .45 -.45  .5 -.75 /
\endpicture
=z^{-1}\;
\beginpicture
\setcoordinatesystem units <.5cm,.5cm>         
\setplotarea x from -0.5 to 0.5, y from -1 to 1    
\plot  0 .75    0 -.75 /
\endpicture
\qquad\hbox{and}\qquad
\beginpicture
\setcoordinatesystem units <.5cm,.5cm>         
\setplotarea x from -1 to 2, y from -1 to 1    
\plot 1.5 .75  1.5 -0.75 /
\plot 0.8 -1.25  1.2 -1.25 /
\plot 0.8 1.25    1.2 1.25   /
\setquadratic
\plot  1.2 -1.25  1.45 -1.05  1.5 -0.75 /
\plot  0.5 -0.75  .55 -1.05  0.8 -1.25 /
\plot  1.2 1.25  1.45 1.05  1.5 0.75 /
\plot  .5 0.75  .55 1.05  0.8 1.25 /
\setquadratic
\plot  -.5 .75  -.45 0.45  -0.1 0.1 /
\plot  .1 -.15  .45 -.45  .5 -.75 /
\plot -.5 -.75  -.45 -.45  0 0  .45 .45  .5 .75 /
\endpicture
= z
\beginpicture
\setcoordinatesystem units <.5cm,.5cm>         
\setplotarea x from -0.5 to 0.5, y from -1 to 1    
\plot  0 .75    0 -.75 /
\endpicture
,\formula$$
$$
\beginpicture
\setcoordinatesystem units <.5cm,.5cm>         
\setplotarea x from -4 to 1.5, y from -2 to 2    
\put{$\hbox{$r$ loops\ }\Bigg\{$} at -3.8 .4   %
\put{$\bullet$} at  0 2      
\put{$\bullet$} at  0 -1.5          
\plot -1.5 2.2 -1.5 1.12 /
\plot -1.5 .88 -1.5 0.12 /
\plot -1.5 -0.12 -1.5 -0.88 /
\plot -1.5 -1.12 -1.5 -1.75 /
\plot -1.25 2.2 -1.25 1.12 /
\plot -1.25 .88 -1.25 0.12 /
\plot -1.25 -0.12 -1.25 -0.88 /
\plot -1.25 -1.12 -1.25 -1.75 /
\ellipticalarc axes ratio 1:1 360 degrees from -1.5 2.2 center 
at -1.375 2.2
\put{$*$} at -1.375 2.2  
\ellipticalarc axes ratio 1:1 180 degrees from -1.5 -1.75 center 
at -1.375 -1.75 
\plot  1 2   1 -1.5 /
\setlinear
\plot -0.3 1.5  -1.1 1.5 /
\ellipticalarc axes ratio 2:1 180 degrees from -1.65 1.5  center 
at -1.65 1.25 
\plot -1.65 1  -1.1 1 /
\ellipticalarc axes ratio 2:1 -180 degrees from -1.1 1  center 
at -1.1 .75 
\ellipticalarc axes ratio 2:1 180 degrees from -1.65 .5  center 
at -1.65 .25 
\ellipticalarc axes ratio 2:1 -180 degrees from -1.1 0  center 
at -1.1 -0.25 
\plot -1.65 0  -1.1 0 /
\ellipticalarc axes ratio 2:1 180 degrees from -1.65 -.5  center 
at -1.65 -.75 
\plot -1.65 -1  -0.3 -1 /
\setquadratic
\plot  -0.3 1.5  -0.05 1.7  -0 2 /
\plot  -0.3 -1  -0.05 -1.2  -0 -1.5 /
\ellipticalarc axes ratio 1:1 180 degrees from 1 2 center 
at 0.5 2
\ellipticalarc axes ratio 1:1 180 degrees from 0 -1.5 center 
at 0.5 -1.5
\endpicture
= Q_r
\beginpicture
\setcoordinatesystem units <.5cm,.5cm>         
\setplotarea x from -5 to -4, y from -1 to 1    
\plot -4.5 1.25 -4.5 -1.25 /
\plot -4.25 1.25 -4.25 -1.25 /
\ellipticalarc axes ratio 1:1 360 degrees from -4.5 1.25 center 
at -4.375 1.25
\put{$*$} at -4.375 1.25  
\ellipticalarc axes ratio 1:1 180 degrees from -4.5 -1.25 center 
at -4.375 -1.25 
\endpicture
\qquad\hbox{and}\qquad 
\beginpicture
\setcoordinatesystem units <.5cm,.5cm>         
\setplotarea x from -4 to 1.5, y from -2 to 2    
\plot -1.5 2.2 -1.5 1.12 /
\plot -1.5 0.88 -1.5 -0.88 /
\plot -1.5 -1.12 -1.5 -1.75 /
\plot -1.25 2.2 -1.25 1.12 /
\plot -1.25 0.88 -1.25 -0.88 /
\plot -1.25 -1.12 -1.25 -1.75 /
\ellipticalarc axes ratio 1:1 360 degrees from -1.5 2.2 center 
at -1.375 2.2
\put{$*$} at -1.375 2.2  
\ellipticalarc axes ratio 1:1 180 degrees from -1.5 -1.75 center 
at -1.375 -1.75 
\plot  1 2   1 0.5 /
\plot  1 -0.5   1 -1.5 /
\setlinear
\plot -0.3 1.5  -1.1 1.5 /
\ellipticalarc axes ratio 2:1 180 degrees from -1.65 1.5  center 
at -1.65 1.25 
\plot -1.65 1  -0.3 1 /
\plot -0.3 -0.5  -1.1 -0.5 /
\ellipticalarc axes ratio 2:1 180 degrees from -1.65 -.5  center 
at -1.65 -.75 
\plot -1.65 -1  -0.3 -1 /
\setquadratic
\plot  -0.3 1  -0.05 0.8  -0 0.5 /
\plot   0 0.5   0.15 0.2   0.7 0 /
\plot   1 -0.5   0.9 -0.15   0.7 0 /
\plot  -0.3 1.5  -0.05 1.7  -0 2 /
\plot   -0.3 -0.5     -.1 -.425  -0.05 -0.325 /
\plot   -0.05 -0.325   0.15 -0.1   0.35 -0.05 /
\plot   0.65 0.15   0.9 0.25   1 0.5 /
\plot  -0.3 -1  -0.05 -1.2  -0 -1.5 /
\ellipticalarc axes ratio 1:1 180 degrees from 1 2 center 
at 0.5 2
\endpicture
=
~z\cdot
\beginpicture
\setcoordinatesystem units <.5cm,.5cm>         
\setplotarea x from -5 to -1.5, y from -1 to 1    
\plot -4.5 1.25 -4.5 -1.25 /
\plot -4.25 1.25 -4.25 -1.25 /
\ellipticalarc axes ratio 1:1 360 degrees from -4.5 1.25 center 
at -4.375 1.25
\put{$*$} at -4.375 1.25  
\ellipticalarc axes ratio 1:1 180 degrees from -4.5 -1.25 center 
at -4.375 -1.25 
\plot -2.7 -0.35    -2.3 -0.35   /
\plot -2 -0.75   -2 -1.25 /
\plot -3 -0.75   -3 -1.25 /
\setquadratic
\plot  -2.3 -0.35  -2.05 -0.45  -2 -0.75 /
\plot  -3 -0.75  -2.95 -0.45  -2.7 -0.35 /
\endpicture
\formula$$
$$
\beginpicture
\setcoordinatesystem units <.5cm,.5cm>         
\setplotarea x from 0 to 1, y from -1 to 1    
\plot 0.3 0.5  0.7 0.5 /
\plot 0.3 -0.5    0.7 -0.5   /
\setquadratic
\plot  0.7 -0.5  0.95 -0.3  1 0 /
\plot  0 0  0.05 -0.3  0.3 -.5 /
\plot  0.7 0.5  0.95 0.3  1 0 /
\plot  0 0  0.05 0.3  0.3 0.5 /
\endpicture
~~=~~
{z-z^{-1}\over q-q^{-1}} + 1 = x.
\formula$$
When working with this algebra it is useful to note that
$$\eqalign{
T_i X^{\varepsilon_{i-1}}T_i X^{\varepsilon_{i-1}}
&= X^{\varepsilon_i} X^{\varepsilon_{i-1}}
=X^{\varepsilon_{i-1}}T_i X^{\varepsilon_{i-1}}T_i,
\qquad\hbox{and, by induction,}\cr
E_i X^{\varepsilon_{i-1}}E_i X^{\varepsilon_{i-1}}
&=E_iT_{i-1}T_iT_i^{-1}X^{\varepsilon_{i-2}}T_{i-1}T_i
X^{\varepsilon_{i-1}} 
=E_i E_{i-1} X^{\varepsilon_{i-2}}T_{i-1}T_i
T_{i-1}^{-1}X^{\varepsilon_{i-1}} \cr
&=E_i E_{i-1} X^{\varepsilon_{i-2}}T_{i-1}
X^{\varepsilon_{i-2}}T_i T_{i-1} 
=z^{-1}E_i E_{i-1}T_i T_{i-1} 
=z^{-1}E_iE_{i-1}E_i \cr
&= z^{-1}E_i. \cr
}
\formula$$

\subsection Schur-Weyl duality for affine and cyclotomic Hecke and BMW algebras

In order to explicitly compute the representations of affine and 
cyclotomic Hecke algebras
and BMW algebras which are obtained by applying the functors
$F_\lambda$ we need to fix notations for working with
the representations of finite dimensional complex semisimple
Lie algebras of classical type.

Let $\fg$ be a complex semisimple Lie algebra of type $A_n$,
$B_n$, $C_n$ or $D_n$ and let $U_h\fg$ be the corresponding 
Drinfeld-Jimbo quantum group.  
Use the notations in [Bou, p. 252-258] for
the root systems of types $A_n$, $B_n$, $C_n$ and $D_n$ so that
$\varepsilon_1,\ldots, \varepsilon_n$ are orthonormal 
(in type $A_n$ also include $\varepsilon_{n+1}$),
$$\matrix{
\hfill \fh^*&=
\left\{\lambda_1\varepsilon_1+\cdots\lambda_{n+1}\varepsilon_{n+1}
\ \big|\ \lambda_i\in \RR, \sum_i\lambda_i=0\right\}, \hfill
&\hbox{in type $A_n$, \quad and}  \hfill\cr
\cr
\hfill \fh^*&=
\left\{\lambda_1\varepsilon_1+\cdots\lambda_n\varepsilon_n
\ |\ \lambda_i\in \RR\right\}, \hfill
&\hbox{in types $B_n$, $C_n$ and $D_n$}, \hfill\cr
}$$
the fundamental weights are given by
$$\matrix{
\hfill \omega_i &= \varepsilon_1+\cdots+\varepsilon_i
-\hbox{$i\over n$}(\varepsilon_1+\cdots+\varepsilon_{n+1}), 
\qquad &1\le i\le n, \hfill &\quad &\hbox{in Type $A_n$}, \hfill \cr
\cr
\hfill \omega_i &= \varepsilon_1+\cdots+\varepsilon_i, 
\qquad 1\le i\le n-1, \hfill &\hbox{in Type $B_n$}, \cr
\hfill \omega_n &= \hbox{$1\over2$}(\varepsilon_1+\cdots+\varepsilon_n), \hfill \cr
\cr
\hfill\omega_i &= \varepsilon_1+\cdots+\varepsilon_i, 
\qquad 1\le i\le n,  \hfill&\hbox{in Type $C_n$}, \cr
\cr
\hfill\omega_i &= \varepsilon_1+\cdots+\varepsilon_i, 
\qquad 1\le i\le n-2,  \hfill \cr
\hfill\omega_{n-1} &= \hbox{$1\over2$}(\varepsilon_1+\cdots
+\varepsilon_{n-1}-\varepsilon_n), \hfill 
&\hbox{in Type $D_n$}, \cr
\hfill\omega_n &= \hbox{$1\over2$}(\varepsilon_1+\cdots
+\varepsilon_{n-1}+\varepsilon_n), \hfill \cr
}$$
and the finite dimensional $U_h\fg$ modules $L(\lambda)$
are indexed by dominant integral weights
$$\matrix{
\lambda = \lambda_1\varepsilon_1+\cdots+\lambda_n\varepsilon_n \hfill
&\ &\lambda_1\ge \lambda_2\ge \cdots\ge \lambda_n\ge 0, \hfill
&\ &\hbox{in Type $A_n$}, \hfill \cr
\qquad\quad -\hbox{$|\lambda|\over n+1$}(\varepsilon_1+\cdots+\varepsilon_{n+1}), 
\hfill
&&\lambda_1,\ldots, \lambda_n\in \ZZ, \hfill \cr
\cr
\lambda = \lambda_1\varepsilon_1+\cdots+\lambda_n\varepsilon_n, \hfill
&&\lambda_1\ge \lambda_2\ge \cdots\ge \lambda_n\ge 0, \hfill \cr
&&\lambda_1,\ldots, \lambda_n\in \ZZ, \hbox{or} \hfill 
&&\hbox{in Type $B_n$}, \hfill \cr
&&\lambda_1,\ldots, \lambda_n\in \hbox{$1\over 2$}+\ZZ, \hfill \cr
\cr
\lambda = \lambda_1\varepsilon_1+\cdots+\lambda_n\varepsilon_n, \hfill
&&\lambda_1\ge \lambda_2\ge \cdots\ge \lambda_n\ge 0, \hfill 
&&\hbox{in Type $C_n$}, \hfill \cr
&&\lambda_1,\ldots, \lambda_n\in \ZZ, \hfill \cr
\cr
\lambda = \lambda_1\varepsilon_1+\cdots+\lambda_n\varepsilon_n, \hfill
&
&\lambda_1\ge \lambda_2\ge \cdots\ge \lambda_{n-1}\ge |\lambda_n|\ge 0, \hfill \cr
&&\lambda_1,\ldots, \lambda_n\in \ZZ, \hbox{or} \hfill 
&&\hbox{in Type $D_n$}, \hfill \cr
&&\lambda_1,\ldots, \lambda_n\in \hbox{$1\over 2$}+\ZZ, \hfill \cr
\cr
}$$
where $|\lambda|=\lambda_1+\cdots+\lambda_n$.  
$$2\rho = \sum_{i=1}^n (y-2i+1)\varepsilon_i,
\qquad\hbox{where}\quad
y=\cases{
n+1, &in type $A_n$, \cr
2n, &in type $B_n$,\cr
2n+1, &in type $C_n$, \cr
2n-1, &in type $D_n$, \cr}
\formula$$
and, in type $A_n$ the sum is over $1\le i\le n+1$ instead of 
$1\le i\le n$.

For all dominant integral weights $\lambda$ in
type $B_n$ and $C_n$ we have
$$L(\lambda)\otimes L(\omega_1)=
\cases{
\displaystyle{\bigoplus_{\lambda^+} L(\lambda^+)}, &in type $A_n$, \cr
\cr
\displaystyle{
L(\lambda)}\bigoplus \left(\bigoplus_{\lambda^\pm} L(\lambda^{\pm})\right),
&in Type $B_n$ with $\lambda_n>0$, \cr
\cr
\displaystyle{\left(\bigoplus_{\lambda^\pm} L(\lambda^{\pm})\right)},
&in types $C_n$ and $D_n$, and\cr
&in type $B_n$ with $\lambda_n=0$,\cr}
\formula$$
where the sum over $\lambda^+$ is a sum over all partitions
(of length $\le n$) obtained by adding a box to $\lambda$, and the
sum over $\lambda^\pm$ denotes a sum over all
dominant weights obtained by adding or removing a box from $\lambda$.
In type $D_n$ addition and removal of a box should include the possibility
of addition and removal of a box marked with a $-$ sign,
and removal of a box from row $n$ when $\lambda_n={1\over2}$ changes
$\lambda_n$ to $-{1\over2}$.

Identify $\lambda$ with the configuration of boxes which has 
$\lambda_i$ boxes in row $i$.  If $\lambda_i\le 0$ put $|\lambda_i|$
boxes in row $i$ and mark them with $-$ signs.  For example
$$\eqalign{
\lambda &= 
\beginpicture
\setcoordinatesystem units <0.25cm,0.25cm>         
\setplotarea x from 0 to 4, y from -3 to 3    
\linethickness=0.5pt                          
\putrule from 0 3 to 5 3          %
\putrule from 0 2 to 5 2          
\putrule from 0 1 to 5 1          %
\putrule from 0 0 to 3 0          %
\putrule from 0 -1 to 3 -1          %
\putrule from 0 -2 to 1 -2          %
\putrule from 0 -3 to 1 -3          %
\putrule from 0 -3 to 0 3        %
\putrule from 1 -3 to 1 3        %
\putrule from 2 -1 to 2 3        %
\putrule from 3 -1 to 3 3        
\putrule from 4 1 to 4 3        %
\putrule from 5 1 to 5 3        %
\endpicture
= \cases{
5\varepsilon_1+5\varepsilon_2
+3\varepsilon_3+3\varepsilon_4+\varepsilon_5+\varepsilon_6
-\hbox{$18\over n+1$}(\varepsilon_1+\cdots+\varepsilon_{n+1}),
\quad
\hbox{in type $A_n$}, \cr
\cr
5\varepsilon_1+5\varepsilon_2
+3\varepsilon_3+3\varepsilon_4+\varepsilon_5+\varepsilon_6,
\quad \hbox{in types $B_n$, $C_n$, and $D_n$}, \cr
}
\cr
\cr
\lambda &=
\beginpicture
\setcoordinatesystem units <0.25cm,0.25cm>         
\setplotarea x from 0 to 4, y from -3 to 3    
\linethickness=0.5pt                          
\putrule from 0 3 to 5.5 3          %
\putrule from 0 2 to 5.5 2          
\putrule from 0 1 to 5.5 1          %
\putrule from 0 0 to 3.5 0          %
\putrule from 0 -1 to 3.5 -1          %
\putrule from 0 -2 to 1.5 -2          %
\putrule from 0 -3 to 1.5 -3          %
\putrule from 0 -3 to 0 3        %
\putrule from 0.5 -3 to 0.5 3        %
\putrule from 1.5 -3 to 1.5 3        %
\putrule from 2.5 -1 to 2.5 3        %
\putrule from 3.5 -1 to 3.5 3        
\putrule from 4.5 1 to 4.5 3        %
\putrule from 5.5 1 to 5.5 3        %
\endpicture
= 
\hbox{$11\over 2$}\varepsilon_1+\hbox{$11\over2$}\varepsilon_2
+\hbox{$7\over2$}\varepsilon_3+\hbox{$7\over2$}\varepsilon_4
+\hbox{$3\over2$}\varepsilon_5+\hbox{$3\over2$}\varepsilon_6,
\qquad\hbox{in Types $B_n$ and $D_n$, and}
\cr
\cr
\lambda &=
\beginpicture
\setcoordinatesystem units <0.25cm,0.25cm>         
\setplotarea x from 0 to 4, y from -3 to 3    
\linethickness=0.5pt                          
\put{$\scriptstyle{-}$} at 0.5 -2.5
\put{$\scriptstyle{-}$} at 1.5 -2.5
\putrule from 0 3 to 6 3          %
\putrule from 0 2 to 6 2          
\putrule from 0 1 to 6 1          %
\putrule from 0 0 to 4 0          %
\putrule from 0 -1 to 4 -1          %
\putrule from 0 -2 to 2 -2          %
\putrule from 0 -3 to 2 -3          %
\putrule from 0 -3 to 0 3        %
\putrule from 1 -3 to 1 3        %
\putrule from 2 -3 to 2 3        %
\putrule from 3 -1 to 3 3        %
\putrule from 4 -1 to 4 3        
\putrule from 5 1 to 5 3        %
\putrule from 6 1 to 6 3        %
\endpicture
= 
6\varepsilon_1+6\varepsilon_2
+4\varepsilon_3+4\varepsilon_4
+2\varepsilon_5-2\varepsilon_6,
\qquad\hbox{in Type $D_6$,}
\cr}$$
If $b$ is a box in position $(i,j)$ of $\lambda$ 
the {\it content} of $b$ is
$$c(b)=j-i = \hbox{the diagonal number of $b$}.
\formula$$
If $\lambda=\lambda_1\varepsilon_1+\cdots\lambda_n\varepsilon_n$,
then
$$\langle \lambda,\lambda+2\rho\rangle 
-\langle \lambda-\varepsilon_i,\lambda-\varepsilon_i+2\rho\rangle 
=2\lambda_i+2\rho_i-1=y+2\lambda_i-2i=y+2c(\lambda/\lambda^-),$$
where $\lambda/\lambda^-$ is the box at the end of row $i$ in $\lambda$.
Note that $c(\lambda/\lambda^-)$ may be a $1\over2$-integer.  Also,
in types $B_n$ and $D_n$,
$$\langle \omega_n,\omega_n+2\rho\rangle 
={n\over4}+{1\over2}\sum_{i=1}^n (y-2i+1) 
= {n\over4}+{n\over2}\cdot y - {n^2\over 2}
=\cases{
{n^2\over2}+{n\over4}, &in type $B_n$, \cr
{n^2\over2}-{n\over4}, &in type $D_n$. \cr
}$$
Using these formulas $\langle \lambda,\lambda+2\rho\rangle$
can easily be computed for all dominant integral weights $\lambda$.
For example
$$\langle \lambda,\lambda+2\rho\rangle
=y|\lambda|+2\sum_{b\in \lambda} c(b)
+\cases{
\hfill \displaystyle{-{|\lambda|^2\over n+1}}\;,\hfill
&in type $A_n$, \cr
\cr
\hfill \displaystyle{0}\;,\hfill 
&in type $C_n$ or in type $B_n$ with $\lambda_i\in \ZZ$, \cr
\cr
\hfill \displaystyle{{n\over4}+{n^2\over 2}}\;,\hfill
&in type $B_n$ with $\lambda_i\in {1\over2}+\ZZ$. \cr
}\formula$$

\thm
Let $\fg$ be the simple complex Lie algebra of classical type,
$U=U_h\fg$ the corresponding quantum group and let $V=L(\omega_1)$ be the 
irreducible of $U_h\fg$ of highest weight $\omega_1$.  
For each $M\in {\cal O}$ let
$\Phi_k\colon\tilde {\cal B}_k\to \End_U(M\otimes V^{\otimes k})$ 
be the affine braid group representation defined in Proposition 3.5.
\smallskip
\item{(a)} If $\fg$ is type $A_n$ then $\Phi_k$ is a representation
of the affine Hecke algebra $\tilde H_k$ with 
$q = e^{h/2}$.
(In this Type $A_n$ case use a different normalization
of the map $\Phi_k$ and set $\Phi_k(T_i) = q^{1/(n+1)}\check R_i$.)
\smallskip
\item{(b)} If $\fg$ is type $A_n$ and if $M=L(\mu)$ where $\mu$ 
is a dominant integral weight then $\Phi_k$ is a representation
of the cyclotomic Hecke algebra $H_{r,1,n}(u_1,\ldots, u_r)$
for any (multi)set of parameters $u_1,\ldots, u_r$ containing the 
(multi)set of values $q^{2c(b)}$ as $b$ runs over the addable boxes of $\mu$.
\smallskip
\item{(c)} If $\fg$ is type $B_n$, $C_n$ or $D_n$ and $M$ is a
highest weight module then there are unique values
$Q_1, Q_2,\ldots \in \CC$, depending only on the central character of $M$,
such that $\Phi_k$ is a representation of the affine BMW algebra $\tilde {\cal Z}_k$ 
with parameters $Q_1, Q_2,\ldots$,
$$q = e^{h/2},
\qquad\hbox{and}\qquad
z=\cases{
q^{2n}, &in Type $B_n$, \cr
-q^{2n+1}, &in Type $C_n$, \cr
q^{2n-1}, &in Type $D_n$.\cr}
$$
\item{(d)} If $\fg$ is type $B_n$, $C_n$ or $D_n$,
and $M=L(\mu)$ where $\mu$ is a dominant integral weight then
$\Phi_k$ is a representation
of the cyclotomic BMW algebra $\tilde {\cal Z}_{r,1,k}$ with $q$ and $z$ as in 
(c),
$$Q_r = \sum_{\mu^{\pm}} q^{r\tilde c(\mu^{\pm},\mu)}
{\dim_q(L(\mu^{\pm}))\over \dim_q(L(\mu))\dim_q(L(\omega_1))},
\qquad r\in \ZZ_{>0},$$
and any (multi)set of parameters $u_1,\ldots, u_r$ containing the 
(multi)set of values $q^{\tilde c(\mu^{\pm},\mu)}$ as $\mu^{\pm}$ runs over 
the dominant integral weights appearing in the decomposition (6.14)
of $L(\mu)\otimes L(\omega_1)$. Here
$$\tilde c(\mu^\pm,\mu) = \cases{
-y, &if $\mu\pm=\mu$, \cr
2c(\mu^\pm/\mu), &if $\mu^\pm\supseteq\mu$, \cr
-2(c(\mu/\mu^\pm)+y), &if $\mu^{\pm}\subseteq \mu$, \cr}$$
where $y$ and $c(b)$ are as defined in (6.13) and (6.15), respectively.
\pf 
(a) It is only necessary to show that $\Phi_k(T_i)=q^{1/(n+1)}\check R_i$ satisfies
$(q^{1/(n+1)}\check R_i)^2=(q-q^{-1})(q^{1/(n+1)}\check R_i)+1$ for $2\le i\le n$.  
This is proved in [LR, Prop. 4.4].  

\smallskip\noindent
(c)  The arguments establishing the relations (6.3a-c) in the definition of the
affine BMW algebra are exactly as in [LR, Prop. 5.10].  It remains to establish
(6.3d-e).   
The element $E_1$ in the affine BMW algebra acts on $V^{\otimes 2}$ as 
$x\cdot {\rm pr}_0$ where ${\rm pr}_0$ is the unique $U_h\fg$-invariant 
projection onto the invariants in $V^{\otimes 2}$ and $x$ is as in (6.5).
Using the identity (6.9) the pictorial equalities
$$
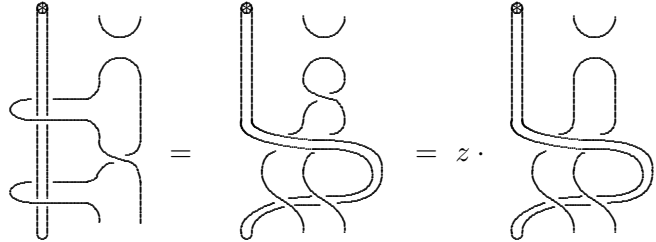
\beginpicture
\setcoordinatesystem units <.5cm,.5cm>         
\setplotarea x from -4 to 1.5, y from -2 to 4    
\plot -1.5 3.7 -1.5 1.12 /
\plot -1.5 0.88 -1.5 -0.88 /
\plot -1.5 -1.12 -1.5 -1.75 /
\plot -1.25 3.7 -1.25 1.12 /
\plot -1.25 0.88 -1.25 -0.88 /
\plot -1.25 -1.12 -1.25 -1.75 /
\ellipticalarc axes ratio 1:1 360 degrees from -1.5 3.7 center 
at -1.375 3.7
\put{$*$} at -1.375 3.7  
\ellipticalarc axes ratio 1:1 180 degrees from -1.5 -1.75 center 
at -1.375 -1.75 
\plot  1 2   1 0.5 /
\plot  1 -0.5   1 -1.5 /
\setlinear
\plot -0.3 1.5  -1.1 1.5 /
\ellipticalarc axes ratio 2:1 180 degrees from -1.65 1.5  center 
at -1.65 1.25 
\plot -1.65 1  -0.3 1 /
\plot -0.3 -0.5  -1.1 -0.5 /
\ellipticalarc axes ratio 2:1 180 degrees from -1.65 -.5  center 
at -1.65 -.75 
\plot -1.65 -1  -0.3 -1 /
\setquadratic
\plot  -0.3 1  -0.05 0.8  -0 0.5 /
\plot   0 0.5   0.15 0.2   0.7 0 /
\plot   1 -0.5   0.9 -0.15   0.7 0 /
\plot  -0.3 1.5  -0.05 1.7  -0 2 /
\plot   -0.3 -0.5     -.1 -.425  -0.05 -0.325 /
\plot   -0.05 -0.325   0.15 -0.1   0.35 -0.05 /
\plot   0.65 0.15   0.9 0.25   1 0.5 /
\plot  -0.3 -1  -0.05 -1.2  -0 -1.5 /
\ellipticalarc axes ratio 1:1 180 degrees from 1 2 center 
at 0.5 2
\ellipticalarc axes ratio 1:1 180 degrees from 0 3.5 center 
at 0.5 3.5
\endpicture
=
\beginpicture
\setcoordinatesystem units <.5cm,.5cm>         
\setplotarea x from -2.5 to 2.5, y from -2 to 4    
\plot -1.5 3.7 -1.5 1 /
\plot -1.25 3.7 -1.25 1 /
\ellipticalarc axes ratio 1:1 360 degrees from -1.5 3.7 center 
at -1.375 3.7
\put{$*$} at -1.375 3.7  
\ellipticalarc axes ratio 1:1 180 degrees from -1.5 -1.75 center 
at -1.375 -1.75 
\setlinear
\setquadratic
\plot   0 2   0.15 1.7   0.7 1.5 /
\plot   1 1   0.9 1.35   0.7 1.5 /
\plot   0 1   0.1 1.3   0.35 1.45 /
\plot   0.65 1.65   0.9 1.75   1 2 /
\plot   -1.5 1   -1.15 0.55   0.215 0.295 /
\plot   1.65 -0.2    1.25 0.2   0.215 0.295 /
\plot   -1.25 1   -0.9 0.7   0.205 0.515 /
\plot    1.9 -0.2    1.5 0.35   0.205 0.515 /
\plot   0.65 0.65   0.9 0.75   1 1 /
\plot   -0.35 0.65  -0.1 0.75   0 1 /
\plot   -1 -0.2   -0.9 0.1   -0.65 0.25 /
\plot   -0 -0.2    0.1 0.1    0.25 0.2 /
\plot   0 -0.2   0.1 -0.6   0.5 -0.95 /
\plot   1 -1.7   0.9 -1.3   0.5 -0.95 /
\plot   -1 -0.2   -0.9 -0.6   -0.5 -0.95 /
\plot   0 -1.7   -0.1 -1.3   -0.5 -0.95 /
\plot  1.9 -0.2    1.65 -0.8    0.8 -1 /
\plot  -1.25 -1.75  -1.05 -1.35  -0.45 -1.1 /   
\plot  -0.35 -0.9    -0.15 -0.875   0.2  -0.85 /
\plot  -0.2 -1.05    0.05 -1.05   0.4  -1.05 /
\plot  1.65 -0.2    1.45 -0.7   0.55  -0.85 /
\plot  -1.5 -1.75  -1.35 -1.35  -0.7 -1 /   
\ellipticalarc axes ratio 1:1 180 degrees from 1 2 center 
at 0.5 2
\ellipticalarc axes ratio 1:1 180 degrees from 0 3.5 center 
at 0.5 3.5
\endpicture
= ~z\cdot
\beginpicture
\setcoordinatesystem units <.5cm,.5cm>         
\setplotarea x from -2 to 1.5, y from -2 to 4    
\plot -1.5 3.7 -1.5 1 /
\plot -1.25 3.7 -1.25 1 /
\ellipticalarc axes ratio 1:1 360 degrees from -1.5 3.7 center 
at -1.375 3.7
\put{$*$} at -1.375 3.7  
\ellipticalarc axes ratio 1:1 180 degrees from -1.5 -1.75 center 
at -1.375 -1.75 
\plot  0 2   0 1 /
\plot  1 2   1 1 /
\setlinear
\setquadratic
\plot   -1.5 1   -1.15 0.55   0.215 0.295 /
\plot   1.65 -0.2    1.25 0.2   0.215 0.295 /
\plot   -1.25 1   -0.9 0.7   0.205 0.515 /
\plot    1.9 -0.2    1.5 0.35   0.205 0.515 /
\plot   0.65 0.65   0.9 0.75   1 1 /
\plot   -0.35 0.65  -0.1 0.75   0 1 /
\plot   -1 -0.2   -0.9 0.1   -0.65 0.25 /
\plot   -0 -0.2    0.1 0.1    0.25 0.2 /
\plot   0 -0.2   0.1 -0.6   0.5 -0.95 /
\plot   1 -1.7   0.9 -1.3   0.5 -0.95 /
\plot   -1 -0.2   -0.9 -0.6   -0.5 -0.95 /
\plot   0 -1.7   -0.1 -1.3   -0.5 -0.95 /
\plot  1.9 -0.2    1.65 -0.8    0.8 -1 /
\plot  -1.25 -1.75  -1.05 -1.35  -0.45 -1.1 /   
\plot  -0.35 -0.9    -0.15 -0.875   0.2  -0.85 /
\plot  -0.2 -1.05    0.05 -1.05   0.4  -1.05 /
\plot  1.65 -0.2    1.45 -0.7   0.55  -0.85 /
\plot  -1.5 -1.75  -1.35 -1.35  -0.7 -1 /   
\ellipticalarc axes ratio 1:1 180 degrees from 1 2 center 
at 0.5 2
\ellipticalarc axes ratio 1:1 180 degrees from 0 3.5 center 
at 0.5 3.5
\endpicture
$$
it follows that $\Phi_2(E_1X^{\varepsilon_1}T_1X^{\varepsilon_1})$
acts as $xz\cdot \check R_{L(0),M}\check R_{M,L(0)}(\id_M\otimes {\rm pr}_0)$.
By (2.11), this is equal to
$$z\cdot(C_M\otimes C_{L(0)})C_{M\otimes L(0)}^{-1}\Phi_2(\id_M\otimes E_1)
=z\cdot C_M C_M^{-1}\Phi_2(\id_M\otimes E_1) = z\cdot\Phi_2(E_1),$$
establishing the relation in (6.3d).  

Since $\Phi_2(E_1)$ acts as $x\cdot(\id_M\otimes {\rm pr}_0)$ on $M\otimes V^{\otimes 2}$
the morphism $\Phi_2(E_1X^{r\varepsilon_1}E_1)$ is a morphism
from $M\otimes L(0)\to M\otimes L(0)$.  Since $M=M\otimes L(0)$ is a 
highest weight module this morphism is $Q_r\cdot\id_M$, for some
$Q_r\in \CC$.  By the results of Drinfeld [Dr] and Reshetikhin [Re] 
(see [Ba, p. 250]),
the action of the morphism $\Phi_2(E_1X^{r\varepsilon_1}E_1)$ corresponds
to the action of a central element of $U_h\fg$ on $M$.  Thus the constant 
$Q_r$ depends only on the central character of $M$.

\smallskip\noindent
(b) Let $b_1,\ldots, b_r$ be the addable boxes of $\mu$ and 
consider the action of
$X^{\varepsilon_1}$ on $M\otimes V=L(\mu)\otimes L(\omega_1)$.
We will show that $\Phi_k(X^{\varepsilon_1})=\check R_0^2$ satisfies the
relation $(\check R_0^2-u_1)\cdots (\check R_0^2-u_r)=0$,
where $u_i=q^{2c(b_i)}$.
By (2.11) and (2.14) it follows that
$$\check R_0^2 = \sum_{\mu^+} 
q^{\langle \mu^+,\;\mu^++2\rho\rangle
-\langle \mu\;,\;\mu+2\rho\rangle
-\langle \omega_1\;,\;\omega_1+2\rho\rangle}
P_{\mu,\omega_1}^{\mu^+}
= \sum_{\mu^+} q^{2c(\mu^+/\mu)} 
P_{\mu,\omega_1}^{\mu^+},$$
where the sum is over all partitions $\mu^+$ obtained by adding a box
to $\mu$, $P_{\mu,\omega_1}^{\mu^+}$ is the projection
onto $L(\mu^+)$ in the tensor product 
$M\otimes V=L(\mu)\otimes L(\omega_1)$, and
$c(\mu^+/\mu)$ is the content of the box $\mu^+/\mu$ which is
added to $\mu$ to get $\mu^+$.  Thus $\check R_0^2$ is a
diagonal operator with eigenvalues $q^{2c(\mu^+/\mu)}$ and so
it satisfies the equation (6.2).

\smallskip\noindent
(d)  Using the appropriate case of the decomposition rule for 
$L(\mu)\otimes L(\omega_1)$, the proof of the relation 
$(X^{\varepsilon_1}-u_1)\cdots (X^{\varepsilon_1}-u_r)=0$
is as in (b).  The values of $\tilde c(\mu^{\pm},\mu)$ are determined
from (6.16).  To compute the value of $Q_r$ note that
$\Phi_2(E_1X^{r\varepsilon_1}E_1)=\Phi_2(\varepsilon_1(X^{r\varepsilon_1})E_1)$,
in the notations of the proof of Theorem 5.4.
Thus $Q_r$ is determined by the formula in Theorem 5.4 
and the decomposition of $L(\mu)\otimes L(\omega_1)$ in (6.14).
\endpf

\medskip\noindent
{\bf Remark.}  The parameters in $Q_1,Q_2,\ldots\in \CC$ needed in 
Theorem 6.17c can be determined 
by using the formula of Baumann [Ba, Theorem 1] 
which characterizes $Q_r$ in terms of the values $Q_1$ given in (5.8).
To do this it is necessary to use formula (5.8) for $Q_1$ several 
times: $\lambda$ is always the highest weight of $M$, but many 
different $\mu$ will be needed.  Note that the proof of the
formula (5.8) for $Q_1$ in [TW] does not require $\lambda$ to be dominant
integral.
\endthm

The following theorem provides an analogue of Schur-Weyl
duality for the affine Hecke algebras,
cyclotomic Hecke algebras, affine BMW algebras and 
cyclotomic BMW algebras.  Alternative
Schur-Weyl dualities have been given by Chari-Pressley [CP2] 
for the case of affine Hecke algebras 
and by Sakamoto and Shoji [SS] for cyclotomic Hecke algebras.
Cherednik [Ch] also used a Schur-Weyl
duality for the affine Hecke algebra which is different from the
Schur-Weyl duality given here.

\thm   Assume that $\fg$ is not of type $D_n$.
Let $\mu$ be a dominant integral weight and let 
$M=L(\mu)$.  In each of the cases given in Theorem 6.17 the 
representation $\Phi_k$ is surjective.
\pf
Part (a) is a consequence of (b) since the representation of $\tilde H_k$
in (a) is the composition of the representation
$\Phi_k\colon H_{r,1,k} \to \End_{U_h\fg}(L(\mu)\otimes V^{\otimes k})$ 
from (b) with the surjective algebra homomorphism
$\tilde H_k \to H_{r,1,k}$ coming from the definition of $H_{r,1,k}$. 
Similarly part (d) is a consequence of part (c).
The proof of the surjectivity of the representation
in Theorem 6.17b and Theorem 6.17d
are exactly the same as the proofs
of [LR, Cor. 4.15]
and [LR, Cor. 5.22], respectively.
The case considered there is the $\mu=0$ case
but all the arguments there generalize verbatim
to the case when $\mu$ is an
arbitrary dominant integral weight.
In [LR, \S 4] the elements $X^{\varepsilon_i}$ in 
the affine braid group are denoted $D_i$.  
The assumption $n>>k$ in [LR] is uncessary for this theorem
if the full decomposition rule given in (6.14) is used.

The main point is that the eigenvalues of $X^{\varepsilon_1},
\ldots, X^{\varepsilon_i}$
separate the components of the decomposition of 
$L(\mu)\otimes V^{\otimes i}$.  By induction it is sufficient
to check that the eigenvalues of $X^{\varepsilon_1}$ distiguish
the components of $L(\lambda)\otimes V$ for all $\lambda$.
By (2.10), (2.11) and (2.14), the eigenvalues
of $X^{\varepsilon_i}$ are of the form $q^{2\tilde c(\lambda^\pm,\lambda)}$
where $\lambda$ is a dominant integral weight
$\tilde c(\lambda^{\pm},\lambda)$ is as in Theorem 6.17d and
$\lambda^{\pm}$ runs over the components in the decomposition
(6.14) of $L(\lambda)\otimes V$.  Different addable boxes for $\lambda$
can never have the same content since they cannot be in the same
diagonal.  Similarly for two different removable boxes.  Let
$b$ be an addable box and $b'$ a removable box for $\lambda$.
Unless $\fg$ is type $D_n$ and $b$ and $b'$ are in row $n$, 
we have $c(b),c(b')\ge -n-1$.  Thus, when $\fg$ is not of type $D_n$,
$c(b) \ne -c(b')-y$ and so
the two eigenvalues coming from these boxes are different.
\endpf

Let $\tilde {\cal Z}_k$ denote the affine Hecke algebra, 
the cyclotomic Hecke algebra, 
the affine BMW algebra or the cyclotomic BMW algebra
corresponding to the case of Theorem 6.17 which is being
considered.  Then, as in the classical Schur-Weyl duality setting,
Theorem 6.18 implies that as $(U_h\fg, \tilde {\cal Z}_k)$ bimodules
$$
L(\mu)\otimes V^{\otimes k} \cong \bigoplus_\lambda
L(\lambda)\otimes {\cal L}^{\lambda/\mu},
\formula$$
where $L(\lambda)$ is the irreducible $U_h\fg$-module
of highest weight $\lambda$ and ${\cal L}^{\lambda/\mu}$ is
the irreducible $\tilde {\cal Z}_k$ module defined by 4.1.

The irreducible $\tilde {\cal Z}_k$ modules ${\cal L}^{\lambda/\mu}$
appearing in (6.19) can be constructed quite explicitly.  All the 
necessary computations for doing this have already been done 
in [LR, \S 4 and 5] which does the case $\mu=0$.
All the arguments in [LR, \S 4 and 5] generalize 
directly to the case when $\mu$ is an
arbitrary dominant integral weight.  The final result
is Theorem 6.20 below.  The result in part (a) of
Theorem 6.20 is due to Cherednik [Ch].

If $\lambda$ and $\mu$ are partitions such that
$\lambda\supseteq\mu$ the {\it skew shape} $\lambda/\mu$ is
the configuration of boxes of in $\lambda$ which
are not in $\mu$. 
Let $\lambda/\mu$ be a skew shape with $k$ boxes.  A 
{\it standard tableau} of shape $\lambda/\mu$ is a filling
$T$ of the boxes of $\lambda/\mu$ with $1,2,\ldots, k$ such that
\smallskip\noindent
\itemitem{(a)} the rows of $T$ are increasing (left to right), and
\smallskip\noindent
\itemitem{(b)} the columns of $T$ are increasing (top to bottom).
\smallskip\noindent
For example, 
$$
\beginpicture
\setcoordinatesystem units <0.5cm,0.5cm>         
\setplotarea x from 0 to 4, y from 0 to 3    
\linethickness=0.5pt                          
\put{11} at 0.5 0.5
\put{6} at 0.5 1.5
\put{8} at 1.5 1.5
\put{2}  at 3.5 2.5
\put{7} at 4.5 3.5
\put{1} at 4.5 4.5
\put{13}  at 5.5 3.5
\put{5}  at 5.5 4.5
\put{3}  at 5.5 5.5
\put{14}  at 6.5 3.5
\put{10}  at 6.5 4.5
\put{4}  at 6.5 5.5
\put{9}  at 7.5 5.5
\put{12}  at 8.5 5.5
\putrule from 5 6 to 9 6          %
\putrule from 4 5 to 9 5          
\putrule from 4 4 to 7 4          %
\putrule from 3 3 to 7 3          %
\putrule from 3 2 to 4 2          %
\putrule from 0 2 to 2 2          %
\putrule from 0 1 to 2 1          %
\putrule from 0 0 to 1 0          %
\putrule from 0 0 to 0 2        %
\putrule from 1 0 to 1 2        %
\putrule from 2 1 to 2 2        %
\putrule from 3 2 to 3 3        
\putrule from 4 2 to 4 5        %
\putrule from 5 3 to 5 6        %
\putrule from 6 3 to 6 6        %
\putrule from 7 3 to 7 6        %
\putrule from 8 5 to 8 6        %
\putrule from 9 5 to 9 6        %
\endpicture
$$
is a standard tableau of shape $\lambda/\mu = (977421)/(5443)$.

For any two partitions $\mu$ and $\lambda$ an {\it up down tableau
of length $k$ from $\mu$ to $\lambda$} is a sequence of 
partitions
$T=\big(\mu=\tau^{(0)},\tau^{(1)},\ldots,\tau^{(k-1)}, \tau^{(k)}=\lambda\big)$
such that
$$
\hbox{(a)\ \  $\tau^{(i)}\supseteq \tau^{(i-1)}$ 
and $\tau^{(i)}/\tau^{(i-1)}=\square\;$,}
 \qquad \hbox{or}\qquad
\hbox{(b)\ \  $\tau^{(i-1)}\supseteq \tau^{(i)}$ 
and $\tau^{(i-1)}/\tau^{(i)}=\square\;$,}
$$
and, in type $B_n$ the situtation
$\tau^{(i-1)}=\tau^{(i)}$ with $\ell(\tau^{(i-1)})=n$ is also allowed.
Note that a standard tableau $\lambda/\mu$ with $k$ boxes is
exactly an up down tableau of length $k$ from $\mu$ to $\lambda$
where all steps in the sequence satisfy condition (a).

\thm 
\item{(a)}  Let $\lambda/\mu$ be a skew shape with $k$ boxes.  
Then the module 
${\cal L}^{\lambda/\mu} = F_\lambda(L(\mu))$
for the affine Hecke algebra $\tilde H_k$
is irreducible and is given by 
$${\cal L}^{\lambda/\mu} = \hbox{span} \{ v_T \ |\ 
\hbox{$T$ standard tableaux of shape $\lambda/\mu$}\}$$
(so that the symbols $v_T$ are a $\CC$-basis of ${\cal L}^{\lambda/\mu}$)
with $\tilde H_k$ action given by
$$\matrix{
\hfill X^{\varepsilon_i}v_T &= q^{2c(T(i))} v_T, \hfill
&&1\le i\le n,\cr
\hfill
T_jv_T &= (T_j)_{TT}v_T + \sqrt{\big(q^{-1}+(T_j)_{TT}\big)
\big(q^{-1}+(T_j)_{s_jT,s_jT}\big) }
\;\;v_{s_jT}, \hfill &\qquad\quad
&1\le j\le n-1, \cr}
$$
where 
\itemitem{} $(T_i)_{TT}$ is the constant
$\displaystyle{
{q-q^{-1}\over 1-q^{2c(T(i))-c(T(i+1))} } }$,
\smallskip
\itemitem{} $c(b)$ denotes the content of the box $b$,
\smallskip
\itemitem{} $T(i)$ is the box containing $i$ in $T$,
\smallskip
\itemitem{} $s_iT$ is the same filling as $T$ except $i$ and $i+1$ are switched,
and
\smallskip
\itemitem{} $v_{s_iT}=0$ if $s_iT$ is not a standard tableau.
\item{(b)}  Let $\lambda/\mu$ be a pair of partitions.
Then the module 
${\cal L}^{\lambda/\mu} = F_\lambda(L(\mu))$
for the affine BMW algebra $\tilde {\cal Z}_k$
is irreducible and is given by
$${\cal L}^{\lambda/\mu} = \hbox{span} \left\{ v_T \ \Big|\ 
\matrix{
\hbox{$T=\big(\mu=\tau^{(0)},\ldots,\tau^{(k)}=\lambda\big)$ an} \cr 
\hbox{up down tableau of length $k$ from $\mu$ to $\lambda$} \cr}
\right\}$$
(so that the symbols $v_T$ are a $\CC$-basis of ${\cal L}^{\lambda/\mu}$)
with $\tilde {\cal Z}_k$ action given by
$$
X^{\varepsilon_i}v_T = q^{\tilde c(\tau^{(i)},\tau^{(i-1)})} v_T, 
\qquad\qquad 1\le i\le n, $$
$$
E_jv_T = \delta_{\tau^{(j+1)},\tau^{(j-1)}}\cdot \sum_S (E_j)_{ST}v_S, 
\qquad\hbox{and}\qquad
T_jv_T = \sum_S (T_j)_{ST}v_S, 
\qquad 1\le j\le n-1, 
$$
where both sums are over up-down tableaux
$S=\big(\mu=\tau^{(0)},\ldots,\tau^{(i-1)},\sigma^{(i)},
\tau^{(i+1)},\ldots, \tau^{(k)}=\lambda\big)$ that are the same as
$T$ except possibly at the $i$th step and
$$\eqalign{
(E_i)_{ST}&=\epsilon\cdot
{\sqrt{\dim_q(L(\tau^{(i)}))\dim_q(L(\sigma^{(i)})) }\over \dim_q(\tau^{(i-1)})},
\cr
\cr
(T_i)_{ST}&=\cases{
\sqrt{\big(q^{-1}+(T_j)_{TT}\big) \big(q^{-1}+(T_j)_{SS}\big) },   
\quad \hbox{if $\tau^{(i-1)}\ne \tau^{(i+1)}$ and $S\ne T$}, \cr
\cr
\displaystyle{
\left({q-q^{-1}\over 
1-\tilde c(\tau^{(i+1)},\sigma^{(i)})\tilde c(\tau^{(i)},\tau^{(i-1)})^{-1} }\right)
(\delta_{ST}-(E_i)_{ST}), }
\quad\hbox{otherwise}, \cr}
\cr
\cr
\tilde c(\tau^{(i)},\tau^{(i-1)}) &= 
\cases{
z^{-1}, &if $\tau^{(i)}=\tau^{(i-1)}$, \cr
q^{2c(\tau^{(i)}/\tau^{(i-1)})}, &if $\tau^{(i)}\supseteq \tau^{(i-1)}$, \cr
z^{-2}q^{-2c(\tau^{(i-1)}/\tau^{(i)})}, &if $\tau^{(i)}\subseteq \tau^{(i-1)}$, \cr
}
\cr
}$$
and $\epsilon = 1$, in type $B_n$ and $D_n$, and
$\epsilon=-1$ in type $C_n$. 
\endthm

\vfill\eject

\subsection Markov traces on affine and cyclotomic Hecke and BMW algebras

If $M=L(\mu)$ where $\mu$ is a dominant integral weight and $V=L(\omega_1)$
then each of the representations $\Phi_k\colon \tilde{\cal Z}_k\to
\End_{U_h\fg}(M\otimes V^{\otimes k})$ (where $\tilde {\cal Z}_k$ is
the affine Hecke algebra, the cyclotomic Hecke algebra, the affine
BMW algebra, or the cyclotomic BMW algebra) gives rise to a Markov
trace via Theorem 5.4.  The parameters and the weights of these
Markov traces are given by Theorems 5.4 and 5.11.  

In type A case $\lambda/\mu$ is a skew shape with
$k$ boxes and the parameters and the weights of (most of) these traces
have been given in terms of partitions in [GIM].  
In [GIM], $\mu$ is a partition of a special form ([GIM, 2.2(*)]) and
so, in their case, the skew shape $\lambda/\mu$ can be viewed as 
an $r$-tuple of partitions.  Their formulas can be recovered from ours by 
rewriting the quantum dimension $\dim_q(L(\lambda))$ from (5.3)
in terms of the partition as in 
[Mac, I \S 3 Ex. 1]: 
$$\dim_q(L(\lambda)) = \prod_{b\in \lambda} {[n+1+c(b)]\over [h(b)]},
\qquad\hbox{in the type $A_n$ case},
\formula$$
where, if $b$ is the box in position $(i,j)$ of $\lambda$,
then $h(b) = \lambda_i-i+\lambda_j'-j+1$ is the 
{\it hook length} at $b$, and
$[d] = (q^d-q^{-d})/(q-q^{-1})$ for a positive integer $d$. 
Thus, the first formula in [GIM, \S 2.3] coincides with 
$\dim_q(L(\lambda))/(\dim_q(V))^{|\lambda|}$ and so the formula
for the weights of the Markov trace on cyclotomic Hecke algebras
which is given in [GIM, Prop. 2.3] coincides exactly with the 
formula in Theorem 5.11.
From Theorem 5.4, (6.14) and (6.16) it follows that the parameters of 
the Markov trace are $z = q/[n+1]$ and 
$$\eqalign{
Q_r &= \sum_{\mu^+} q^{2c(\mu^+/\mu)}{\dim_q(L(\mu^+))
\over \dim_q(L(\mu))\dim_q(V) } \cr
&= \sum_{\mu^+} q^{2c(\mu^+/\mu)}
\left(\prod_{b\in \mu^+} {[n+1+c(b)]\over [h(b)]}\right)
\left(\prod_{b\in \mu} {[h(b)]\over [n+1+c(b)]}\right)
{1\over [n+1]} \cr
&= \sum_{\mu^+} q^{2c(\mu^+/\mu)}
\left({\prod_{b\in \mu} [h(b)]\over \prod_{b\in \mu^+} [h(b)]}\right)
{[n+1+c(\mu^+/\mu)]\over [n+1]} \cr
&= \sum_{\mu^+} q^{2c(\mu^+/\mu)}
\left(\prod_{b'} {[h(b')]\over [h(b')+1]}\right)
\left(\prod_{b''} {[h(b'')]\over [h(b'')+1]}\right)
{[n+1+c(\mu^+/\mu)]\over [n+1]}, \cr
}$$
where, in the last expression, the first product is over boxes
$b'\in \mu$ which are in the same row as the added box $\mu^+/\mu$
and the second product is over $b''\in \mu$ which are in the 
same column as $\mu^+/\mu$.  Then cancellation of the common terms 
in the numerator and denominators of each product yields the 
combinatorial formulas for the parameters of the Markov traces
on cyclotomic Hecke algebras which are given in [GIM, Thm. 2.4].

Lambropoulou [Lb, \S 4]
has proved that there is a {\it unique} Markov trace
on the affine Hecke algebra with a given choice of parameters
$z, Q_1,\ldots, Q_r\in \CC$.  A similar result is true for
the affine BMW algebra.

\thm For each fixed choice of parameters $q$, $z$ and $Q_1, Q_2,\ldots$
there is a unique Markov trace on the affine BMW algebra ${\cal Z}_k$.
\rm\medskip\noindent{\it Sketch of proof. }
Consider the image of an affine braid $b$ in the
affine BMW algebra.  The Markov trace of this braid can be viewed 
pictorially as the closure of the braid $b$.
$$
\mt_k(b) = \mt_1\left(
\beginpicture
\setcoordinatesystem units <.5cm,.5cm>         
\setplotarea x from -2 to 4, y from -4 to 4    
\put{$b$} at 0 0
\plot  -1.5 1   1.5 1 /
\plot  -1.5 -1   1.5 -1 /
\plot  -1.5 -1   -1.5 1 /
\plot  1.5 -1   1.5 1 /
\plot  -.8 1.5    -.8  1 /
\plot  -.4 1.5    -.4  1 /
\plot    0 1.5     0 1 /
\plot  .4 1.5    .4  1 /
\plot  .8 1.5    .8 1 /
\plot  1.2 1.5    1.2 1 /
\plot  -.8 -1.5    -.8  -1 /
\plot  -.4 -1.5    -.4  -1 /
\plot    0 -1.5     0 -1 /
\plot  .4 -1.5    .4  -1 /
\plot  .8 -1.5    .8 -1 /
\plot  1.2 -1.5    1.2 -1 /
\plot  1.9 1.5    1.9 -1.5 /
\plot  2.3 1.5    2.3 -1.5 /
\plot  2.7 1.5    2.7 -1.5 /
\plot  3.1 1.5    3.1 -1.5 /
\plot  3.5 1.5    3.5 -1.5 /
\plot  3.9 1.5    3.9 -1.5 /
\ellipticalarc axes ratio 1:1 180 degrees from 1.9 1.5 center 
at 1.55 1.5
\ellipticalarc axes ratio 1:1 180 degrees from 2.3 1.5 center 
at 1.55 1.5
\ellipticalarc axes ratio 1:1 180 degrees from 2.7 1.5 center 
at 1.55 1.5
\ellipticalarc axes ratio 1:1 180 degrees from 3.1 1.5 center 
at 1.55 1.5
\ellipticalarc axes ratio 1:1 180 degrees from 3.5 1.5 center 
at 1.55 1.5
\ellipticalarc axes ratio 1:1 180 degrees from 3.9 1.5 center 
at 1.55 1.5
\ellipticalarc axes ratio 1:1 -180 degrees from 1.9 -1.5 center 
at 1.55 -1.5
\ellipticalarc axes ratio 1:1 -180 degrees from 2.3 -1.5 center 
at 1.55 -1.5
\ellipticalarc axes ratio 1:1 -180 degrees from 2.7 -1.5 center 
at 1.55 -1.5
\ellipticalarc axes ratio 1:1 -180 degrees from 3.1 -1.5 center 
at 1.55 -1.5
\ellipticalarc axes ratio 1:1 -180 degrees from 3.5 -1.5 center 
at 1.55 -1.5
\ellipticalarc axes ratio 1:1 -180 degrees from 3.9 -1.5 center 
at 1.55 -1.5
\linethickness=1pt
\putrule  from -1.2 1.5 to  -1.2  1 
\putrule  from  -1.2 -1.5  to  -1.2  -1 
\endpicture
\right)
$$
Consider a string in the closure as it winds around the other strings
and the pole.  If the string crosses another string twice
without going around the pole between these two crossings then we can
use the relation
$$\beginpicture
\setcoordinatesystem units <.5cm,.5cm>         
\setplotarea x from -1 to 1, y from -1 to 1    
\setquadratic
\plot  -.5 0  -.45 0.3  -0.1 0.65 /
\plot  .1 .9  .45 1.2  .5 1.5 /
\plot -.5 1.5  -.45 1.2  0 .75  .45 0.3  .5 0 /
\plot  -.5 -1.5  -.45 -1.2  -0.1 -0.85 /
\plot  .1 -0.6  .45 -0.3  .5 0 /
\plot -.5 0  -.45 -0.3  0 -0.75  .45 -1.2  .5 -1.5 /
\endpicture
=
\beginpicture
\setcoordinatesystem units <.5cm,.5cm>         
\setplotarea x from -1 to 1, y from -1.8 to 1.8    
\plot -0.5 .75 -0.5 -.75 /
\plot  0.5 .75  0.5 -.75 /
\endpicture
+(q-q^{-1})\left(
\beginpicture
\setcoordinatesystem units <.5cm,.5cm>         
\setplotarea x from -1 to 1, y from -1 to 1    
\setquadratic
\plot  -.5 -.75  -.45 -0.45  -0.1 -0.1 /
\plot  .1 .15  .45 .45  .5 .75 /
\plot  -.5 .75  -.45 .45  0 0  .45 -.45  .5 -.75 /
\endpicture
-
z\cdot\beginpicture
\setcoordinatesystem units <.5cm,.5cm>         
\setplotarea x from -1 to 1, y from -1.8 to 1.8    
\plot -0.2 -0.5  0.2 -0.5 /
\plot -0.2 -1    0.2 -1  /
\setquadratic
\plot  0.2 -0.5  0.45 -0.3  .5 0 /
\plot  -0.5 0  -0.45 -0.3  -0.2 -0.5 /
\plot  0.2 -1  0.45 -1.2  0.5 -1.5 /
\plot  -0.5 -1.5  -0.45 -1.2  -0.2 -1 /
\setquadratic
\plot  -.5 0  -.45 0.3  -0.1 0.65 /
\plot  .1 .9  .45 1.2  .5 1.5 /
\plot -.5 1.5  -.45 1.2  0 .75  .45 0.3  .5 0 /
\endpicture
\right)
$$
to rewrite the closed braid as a linear combination of 
closed braids with {\it fewer} crossings between strings.
By successive steps of this type we can reduce the computation of the 
Markov trace of a braid to a linear combination of 
$$
\matrix{
\beginpicture
\setcoordinatesystem units <.5cm,.5cm>         
\setplotarea x from -4 to 1.5, y from -2 to 2    
\put{$\hbox{$r_1$ loops\ }\Bigg\{$} at -3.8 .4   %
\plot -1.5 2.2 -1.5 1.12 /
\plot -1.5 .88 -1.5 0.12 /
\plot -1.5 -0.12 -1.5 -0.88 /
\plot -1.5 -1.12 -1.5 -1.75 /
\plot -1.25 2.2 -1.25 1.12 /
\plot -1.25 .88 -1.25 0.12 /
\plot -1.25 -0.12 -1.25 -0.88 /
\plot -1.25 -1.12 -1.25 -1.75 /
\ellipticalarc axes ratio 1:1 360 degrees from -1.5 2.2 center 
at -1.375 2.2
\put{$*$} at -1.375 2.2  
\plot  0 1   0 -0.5 /
\setlinear
\plot -0.3 1.5  -1.1 1.5 /
\ellipticalarc axes ratio 2:1 180 degrees from -1.65 1.5  center 
at -1.65 1.25 
\plot -1.65 1  -1.1 1 /
\ellipticalarc axes ratio 2:1 -180 degrees from -1.1 1  center 
at -1.1 .75 
\ellipticalarc axes ratio 2:1 180 degrees from -1.65 .5  center 
at -1.65 .25 
\ellipticalarc axes ratio 2:1 -180 degrees from -1.1 0  center 
at -1.1 -0.25 
\plot -1.65 0  -1.1 0 /
\ellipticalarc axes ratio 2:1 180 degrees from -1.65 -.5  center 
at -1.65 -.75 
\plot -1.65 -1  -0.3 -1 /
\setquadratic
\plot  -0.3 1.5  -0.05 1.3  -0 1 /
\plot  -0.3 -1  -0.05 -0.8  -0 -0.5 /
\endpicture
\cr
\enspace\quad\vdots\cr
\cr
\beginpicture
\setcoordinatesystem units <.5cm,.5cm>         
\setplotarea x from -4 to 1.5, y from -2 to 2    
\put{$\hbox{$r_k$ loops\ }\Bigg\{$} at -3.8 .4   %
\plot -1.5 2.2 -1.5 1.12 /
\plot -1.5 .88 -1.5 0.12 /
\plot -1.5 -0.12 -1.5 -0.88 /
\plot -1.5 -1.12 -1.5 -1.75 /
\plot -1.25 2.2 -1.25 1.12 /
\plot -1.25 .88 -1.25 0.12 /
\plot -1.25 -0.12 -1.25 -0.88 /
\plot -1.25 -1.12 -1.25 -1.75 /
\ellipticalarc axes ratio 1:1 180 degrees from -1.5 -1.75 center 
at -1.375 -1.75 
\plot  0 1   0 -0.5 /
\setlinear
\plot -0.3 1.5  -1.1 1.5 /
\ellipticalarc axes ratio 2:1 180 degrees from -1.65 1.5  center 
at -1.65 1.25 
\plot -1.65 1  -1.1 1 /
\ellipticalarc axes ratio 2:1 -180 degrees from -1.1 1  center 
at -1.1 .75 
\ellipticalarc axes ratio 2:1 180 degrees from -1.65 .5  center 
at -1.65 .25 
\ellipticalarc axes ratio 2:1 -180 degrees from -1.1 0  center 
at -1.1 -0.25 
\plot -1.65 0  -1.1 0 /
\ellipticalarc axes ratio 2:1 180 degrees from -1.65 -.5  center 
at -1.65 -.75 
\plot -1.65 -1  -0.3 -1 /
\setquadratic
\plot  -0.3 1.5  -0.05 1.3  -0 1 /
\plot  -0.3 -1  -0.05 -0.8  -0 -0.5 /
\endpicture
\cr
}
= {Q_{r_1}\cdots Q_{r_k}\over \dim_q(V)^k}\;\cdot\;\mt
\left(
\beginpicture
\setcoordinatesystem units <.5cm,.5cm>         
\setplotarea x from -5 to -4, y from -1.5 to 1.5    
\plot -4.5 1.25 -4.5 -1.25 /
\plot -4.25 1.25 -4.25 -1.25 /
\ellipticalarc axes ratio 1:1 360 degrees from -4.5 1.25 center 
at -4.375 1.25
\put{$*$} at -4.375 1.25  
\ellipticalarc axes ratio 1:1 180 degrees from -4.5 -1.25 center 
at -4.375 -1.25 
\endpicture
\right)
= {Q_{r_1}\cdots Q_{r_k}\over \dim_q(V)^k}.
\qquad\hbox{\qed}$$

\medskip\noindent
{\it Remark.} For computations it is helpful to note that 
$$\dim_q(L(\omega_1)) = \cases{
[n+1], &in type $A_n$, \cr
[2r]+1, &in type $B_r$, \cr
[2r+1]-1, &in type $C_r$, \cr
[2r-1]+1, &in type $D_r$. \cr}
\formula$$

\vfill\eject

\subsection Standard and simple modules for affine Hecke algebras

The original construction of the irreducible representations
of the affine Hecke algebra of type A is due to Zelevinsky
[Ze2] and is an analogue of the Langlands construction of
admissible representations of real reductive Lie groups.
Zelevinsky used the combinatorics of multisegments
which is easily seen to be equivalent to the combinatorics 
of unipotent-semisimple pairs used later in [KL] (see [Ar]).
Here we show how the construction of affine Hecke algebra
representations via the functors $F_\lambda$  naturally matches 
up with Zelevinsky's indexings by multisegments.
Using the multisegment indexing of representations,
Theorem 6.31 below explicitly matches up the decomposition numbers for
affine Hecke algebras with Kazdhan-Lusztig polynomials. 
Recall that the functor $F_\lambda$ gives representations
of the affine Hecke algebra in the setting of Theorem 6.17a when
$\fg$ is of type $A_n$ and $V=L(\omega_1)$ is the $n$-dimensional
fundamental representation.

Consider an (infinite) sheet of graph paper which has its diagonals
labeled consecutively by $\ldots, -2, -1, 0, 1, 2, \ldots$.  
The {\it content}
$c(b)$ of a box  $b$ on this sheet of graph paper is
$$c(b) = \hbox{the diagonal number of the box $b$}$$
(a natural generalization of the definition of $c(b)$ in (6.15)).
A {\it multisegment} is a collection of rows of boxes (segments) placed
on graph paper.  We can label this multisegment by a pair of weights
$\lambda = \lambda\varepsilon_1+\cdots\lambda_n\varepsilon_n$
and $\mu=\mu_1\varepsilon_1+\cdots+\mu_n\varepsilon_n$
by setting
$$\eqalign{
(\lambda+\rho)_i &= \hbox{content of the last box in row $i$,}
\qquad\hbox{and}\qquad \cr
(\mu+\rho)_i &= \hbox{(content of the first box in row $i$)}-1. \cr
}$$
For example
$$
\beginpicture
\setcoordinatesystem units <0.5cm,0.5cm>         
\setplotarea x from 0 to 7, y from -2 to 3    
\linethickness=0.5pt                          
\put{1} at 1.5 -0.5
\put{2} at 2.5 -0.5
\put{3} at 0.5 2.5
\put{3}  at 1.5 1.5
\put{3} at 3.5 -0.5
\put{3} at 4.5 -1.5
\put{4}  at 1.5 2.5
\put{4}  at 2.5 1.5
\put{4}  at 4.5 -0.5
\put{4}  at 5.5 -1.5
\put{5}  at 2.5 2.5
\put{5}  at 3.5 1.5
\put{5}  at 4.5 0.5
\put{5}  at 5.5 -0.5
\put{5}  at 6.5 -1.5
\put{6}  at 3.5 2.5
\put{6}  at 4.5 1.5
\put{6}  at 5.5 0.5
\put{7}  at 4.5 2.5
\put{7}  at 5.5 1.5
\put{7}  at 6.5 0.5
\putrule from 0 3  to 5 3          
\putrule from 0 2  to 6 2          %
\putrule from 1 1  to 7 1          %
\putrule from 1 0  to 7 0          %
\putrule from 1 -1 to 7 -1          %
\putrule from 4 -2 to 7 -2          %
\putrule from 0 2 to 0 3        %
\putrule from 1 1 to 1 3        %
\putrule from 1 -1 to 1 0        %
\putrule from 2 1 to 2 3        %
\putrule from 2 -1 to 2 0        %
\putrule from 3 1 to 3 3        
\putrule from 3 -1 to 3 0        
\putrule from 4 -2 to 4 3        %
\putrule from 5 -2 to 5 3        %
\putrule from 6 -2 to 6 2        %
\putrule from 7 -2 to 7 -1        %
\putrule from 7 0 to 7 1        %
\endpicture
\qquad\hbox{corresponds to}\qquad
\eqalign{ 
\lambda &= (7,\ 7,\ 7,\ 5,\ 5) \quad\hbox{and}\cr 
\mu &= (2,\ 2,\ 4,\ 0,\ 2) \cr}
\formula$$
(the numbers in the boxes in the picture are the contents of the boxes).
The construction forces the condition
\smallskip
\itemitem{(a)}  $(\lambda+\rho)_i-(\mu+\rho)_i\in \ZZ_{\ge 0}$.
\smallskip\noindent
and since we want to consider {\it unordered} collections of boxes
it is natural to take the following pseudo-lexicographic ordering
on the segments
\smallskip
\itemitem{(b)} $(\lambda+\rho)_i\ge (\lambda+\rho)_{i+1}$,
\smallskip
\itemitem{(c)} $(\mu+\rho)_i\le (\mu+\rho)_{i+1}$ if $(\lambda+\rho)_i=(\lambda+\rho)_{i+1}$, 
\smallskip\noindent
when we denote the multisegment $\lambda/\mu$ by a pair of weights $\lambda,\mu$.
In terms of weights the conditions (a), (b) and (c) can be restated
as (note that in this case both $\lambda$ and $\mu$ are integral)
\smallskip
\itemitem{(a$'$)} $\lambda-\mu$ is a weight of $V^{\otimes k}$, where $k$ is the
number of boxes in $\lambda/\mu$,
\smallskip
\itemitem{(b$'$)} $\lambda$ is integrally dominant,
\smallskip
\itemitem{(c$'$)} $\mu=w\circ \nu$ with $\nu$ integrally
dominant and $w$ maximal length in the coset $W_{\lambda+\rho}wW_{\nu+\rho}$,
\smallskip\noindent
These conditions on the pair of weights $(\lambda,\mu)$
arose previously in Proposition 4.3d and Lemma 4.7.

Let $\lambda/\mu$ be a multisegment with $k$ boxes
and number the boxes of $\lambda/\mu$
from left to right (like a book).  Define
$$\tilde H_{\lambda/\mu} = \hbox{ subalgebra of $\tilde H_k$
generated by\ }
\{X^\lambda, \; T_j \ |\ 
\hbox{$\lambda\in L$, ${\rm box}_j$ is not at the end of its row}\},$$
so that $\tilde H_{\lambda/\mu}$ is the ``parabolic'' subalgebra of $\tilde H_k$
corresponding to the multisegment $\lambda/\mu$.  Define a one-dimensional
$\tilde H_{\lambda/\mu}$ module $\CC_{\lambda/\mu} = \CC v_{\lambda/\mu}$
by setting
$$X^{\varepsilon_i}v_{\lambda/\mu} = q^{2c({\rm box}_i)}v_{\lambda/\mu},
\qquad\hbox{and}\qquad
T_jv_{\lambda/\mu} = qv_{\lambda/\mu},
\formula$$
for $1\le i\le n$ and $j$ such that ${\rm box}_j$ is not at the end of its row.

Let $\fg$ be of type $A_n$ and let $F_\lambda$ be the 
functor $\Hom_{U_h\fg}(M(\lambda), \cdot\otimes V^{\otimes k})$
from the setting of Theorem 6.17a, where $V=L(\omega_1)$.
The {\it standard module} for the affine Hecke algebra $\tilde H_k$
is 
$${\cal M}^{\lambda/\mu}= F_\lambda(M(\mu))\formula$$
as defined in (4.1).  It follows from the above discussion that 
these modules are naturally indexed by multisegments $\lambda/\mu$.
The following proposition shows that this standard module
coincides with the usual standard module for the affine Hecke algebra
as considered by Zelevinsky [Ze2] (see also [Ar], [CG] and [KL]).

\prop  Let $\lambda$ and $\mu$ be integrally dominant weights giving
rise to the multisegment $\lambda/\mu$.  Let $\CC_{\lambda/\mu}$ be the 
one dimensional representation of the parabolic subalgebra
$\tilde H_{\lambda/\mu}$ of the affine Hecke algebra $\tilde H_k$ defined 
in (6.25).  Then
$${\cal M}^{\lambda/\mu} \cong 
\Ind_{\tilde H_{\lambda/\mu}}^{\tilde H_k}(\CC_{\lambda/\mu}).$$
\pf
By Proposition 4.3a, 
${\cal M}^{\lambda/\mu}\cong (V^{\otimes k})_{\lambda-\mu}$ as a vector space.
Let $\{v_1, v_2,\ldots, v_n\}$ be the standard basis of $V=L(\omega_1)$ with
${\rm wt}(v_i)=\varepsilon_i$.
If we let the symmetric group $S_k$ act on $V^{\otimes k}$ by permuting the tensor 
factors then
$$
(V^{\otimes k})_{\lambda-\mu} 
= \hbox{span-}\{ \pi\cdot v^{\otimes(\lambda-\mu)}
\ |\ \pi\in S_k\} 
= \hbox{span-}\{ \pi\cdot v^{\otimes(\lambda-\mu)}
\ |\ \pi\in S_k/S_{\lambda-\mu}\}, 
\qquad\hbox{where}
$$
$$
v^{\otimes(\lambda-\mu)}=\underbrace{v_1\otimes\cdots\otimes v_1}_{\lambda_1-\mu_1}
\otimes\cdots\otimes 
\underbrace{v_n\otimes\cdots\otimes v_n}_{\lambda_n-\mu_n}
\qquad\hbox{and}\qquad
S_{\lambda-\mu}=S_{\lambda_1-\mu_1}\times \cdots\times S_{\lambda_n-\mu_n}$$
is the parabolic subgroup of $S_k$ which stabilizes the vector
$v^{\otimes(\lambda-\mu)}\in V^{\otimes k}$.
This shows that, as vector spaces, 
$${\cal M}^{\lambda/\mu}\cong
\Ind_{\tilde H_{\lambda/\mu}}^{\tilde H_k}(\CC_{\lambda/\mu})
= \hbox{span-} \{ T_\pi\otimes v_{\lambda/\mu}\ |\ \pi\in S_k/S_{\lambda-\mu}\}
\formula$$
are isomorphic.

For notational purposes let
$$b_{\lambda/\mu}
=v_\mu^+\otimes v^{\otimes(\lambda-\mu)}
=v_\mu^+\otimes v_{i_1}\otimes\cdots\otimes v_{i_k}$$
and let $\bar b_{\lambda/\mu}$ be the image of $b_{\lambda/\mu}$
in $(M\otimes V^{\otimes k})^{[\lambda]}$.  Since $\lambda$ is integrally
dominant and $\bar b_{\lambda/\mu}$ has weight $\lambda$ it must be a highest
weight vector.  We will show that $X^{\varepsilon_\ell}$ acts 
on $\bar b_{\lambda/\mu}$ by the constant $q^{c({\rm box}_\ell)}$,
where $c({\rm box}_\ell)$ is the content of the $\ell$th box of the 
multisegment $\lambda/\mu$ (read left to right and top to bottom like a book).

Consider the projections
$${\sl pr}_\ell\colon M(\mu)\otimes V^{\otimes k}
\to (M(\mu)\otimes V^{\otimes \ell})^{[\lambda^{(\ell)}]}
\otimes V^{\otimes (k-\ell)}
\qquad\hbox{where}\qquad
\lambda^{(\ell)}=\mu+\sum_{j\le\ell} {\rm wt}(v_{i_\ell})$$ 
and ${\sl pr}_i$ acts as the identity on the last $k-i$ factors of 
$M(\mu)\otimes V^{\otimes k}$.  Then
$$\bar b_{\lambda/\mu}={\sl pr}_k {\sl pr}_{k-1}\dots {\sl pr}_1 b_{\lambda/\mu},$$
and for each $1\le \ell \le k$,
$pr_{\ell-1}\cdots pr_1(b_{\lambda/\mu})$ is a highest
weight vector of weight $\lambda^{(\ell)}$ in $M\otimes V^{\otimes \ell}$.
It is the ``highest'' highest weight vector of 
$$((M(\mu)\otimes V^{\otimes(\ell-1)})^{[\lambda^{(\ell-1)}]}
\otimes V)^{[\lambda^{(\ell)}]}
\formula$$
with respect to the ordering in Lemma 4.2 and thus it is deepest in the 
filtration constructed there.
Note that the quantum Casimir element acts on the space in (6.29) as the
constant $q^{\langle \lambda^{(\ell)},\lambda^{(\ell)}+2\rho\rangle}$
times a unipotent transformation, and the unipotent transformation
must preserve the filtration coming from Lemma 4.2.  Since
${\sl pr}_\ell(b_\lambda/\mu)$ is the highest weight vector
of the smallest submodule of this filtration
(which is isomorphic to a Verma module by Lemma 4.2b)
it is an eigenvector for the action of the quantum Casimir.
Thus, by (2.11) and (2.13), $X^{\varepsilon_\ell}$ acts on 
${\sl pr}_\ell(b_{\lambda/\mu})$ by the constant
$$
q^{\langle \lambda^{(\ell)},\lambda^{(\ell)}+2\rho\rangle
-\langle \lambda^{(\ell-1)}, \lambda^{(\ell-1)}+2\rho\rangle
-\langle \omega_1, \omega_1+2\rho\rangle}
=q^{c({\rm box}_\ell)}.$$
Since $X^{\varepsilon_\ell}$ commutes with $pr_j$ for $j>\ell$ it
this also specifies the action of $X^{\varepsilon_\ell}$ on 
$\bar b_{\lambda/\mu} = {\sl pr}_\ell(b_{\lambda/\mu})$.

The explicit $R$-matrix $\check R_{VV}\colon V\otimes V\to V\otimes V$
for this case ($\fg$ of type $A$ and $V=L(\omega_1)$) is well known 
(see, for example, the proof of [LR, Prop. 4.4]) and given by
$$\check R_{VV}(v_i\otimes v_j) = 
\cases{ v_j\otimes v_i, &if $i>j$, \cr
(q-q^{-1})v_i\otimes v_j+v_j\otimes v_i, &if $i<j$, \cr
q v_i\otimes v_j, &if $i=j$. \cr
}$$
Since $T_i$ acts by $\check R_{VV}$ on the $i$th and $(i+1)$st
tensor factors of $V^{\otimes k}$ and commutes with the projection
$pr_\lambda$ it follows that
$T_j(\bar b_{\lambda/\mu})=q\,\bar b_{\lambda/\mu}$, 
if ${\rm box}_j$ is not a box at the end of a row of $\lambda/\mu$.
This analysis of the action of $\tilde H_{\lambda/\mu}$ on $\bar b_{\lambda/\mu}$
shows that there is an $\tilde H_k$-homomorphism
$$\matrix{
\Ind_{\tilde H_{\lambda/\mu}}^{\tilde H_k}(\CC v_{\lambda/\mu})
&\longrightarrow &{\cal M}^{\lambda/\mu} \cr
v_{\lambda/\mu} &\longmapsto &\bar b_{\lambda/\mu}. \cr
}$$
This map is surjective since ${\cal M}^{\lambda/\mu}$ is 
generated by $\bar b_{\lambda/\mu}$ (the ${\cal B}_k$ action on $v^{\lambda-\mu}$
generates all of $(V^{\otimes k})_{\lambda-\mu}$).
Finally, (6.28) guarantees that it is an isomorphism.
\endpf

In the same way that each weight $\mu\in \fh^*$ has a 
{\it normal form}
$$\mu = w\circ \tilde\mu,
\qquad\hbox{with}\qquad
\matrix{
\hbox{$\tilde\mu$ integrally dominant, \enspace and\enspace } \hfill\cr
\hbox{$w$ maximal length in the coset $wW_{\tilde \mu+\rho}$,}\hfill \cr
} 
$$
every multisegment $\lambda/\mu$ has a {\it normal form}
$$\lambda/\mu = \nu/(w\circ \tilde \nu),
\qquad\hbox{with}\qquad
\matrix{
\hbox{$\nu$ the sequence of contents of boxes of $\lambda/\mu$},\hfill \cr
\tilde\nu = \nu-(1,1,\ldots, 1),\quad\hbox{and} \hfill \cr
\hbox{$w$ maximal length in $W_{\nu+\rho}wW_{\nu+\rho}$.} \hfill
\cr} 
$$
The element $w$ in the normal form $\nu/(w\circ\tilde\nu)$ 
of $\lambda/\mu$ can be constructed combinatorially by the following
scheme.  We number (order) the boxes of $\lambda/\mu$ in two different ways.
%
\smallskip\noindent
First ordering:
To each box $b$ of $\lambda/\mu$ associate the following triple
$$\big(\hbox{content of the box to the left of $b$}, -\hbox{(content of $b$)},
-\hbox{(row number of $b$)}\big)$$
where, if a box is the leftmost box in a row ``the 
box to its left'' is the rightmost box in the same row.
The lexicographic ordering on these triples induces an
ordering on the boxes of $\lambda/\mu$.
\smallskip\noindent
Second ordering:  To each box $b$ of $\lambda/\mu$ associate the following
pair
$$\big(\hbox{content of $b$}, -\hbox{(the number of box $b$ in the first ordering)}\big)$$
The lexicographic ordering of these pairs induces a second ordering
on the boxes of $\lambda/\mu$.
\smallskip\noindent
Then $w$ is the permutation defined by these two numberings of the boxes.
For example, for the multisegment $\lambda/\mu$ displayed in (6.24)
the numberings of the boxes are given by
$$\matrix{
\beginpicture
\setcoordinatesystem units <0.5cm,0.5cm>         
\setplotarea x from 0 to 7, y from -2 to 3    
\linethickness=0.5pt                          
\put{15} at 1.5 -0.5
\put{1} at 2.5 -0.5
\put{21} at 0.5 2.5
\put{20}  at 1.5 1.5
\put{2} at 3.5 -0.5
\put{14} at 4.5 -1.5
\put{6}  at 1.5 2.5
\put{5}  at 2.5 1.5
\put{4}  at 4.5 -0.5
\put{3}  at 5.5 -1.5
\put{10}  at 2.5 2.5
\put{9}  at 3.5 1.5
\put{19}  at 4.5 0.5
\put{8}  at 5.5 -0.5
\put{7}  at 6.5 -1.5
\put{13}  at 3.5 2.5
\put{12}  at 4.5 1.5
\put{11}  at 5.5 0.5
\put{18}  at 4.5 2.5
\put{17}  at 5.5 1.5
\put{16}  at 6.5 0.5
\putrule from 0 3  to 5 3          
\putrule from 0 2  to 6 2          %
\putrule from 1 1  to 7 1          %
\putrule from 1 0  to 7 0          %
\putrule from 1 -1 to 7 -1          %
\putrule from 4 -2 to 7 -2          %
\putrule from 0 2 to 0 3        %
\putrule from 1 1 to 1 3        %
\putrule from 1 -1 to 1 0        %
\putrule from 2 1 to 2 3        %
\putrule from 2 -1 to 2 0        %
\putrule from 3 1 to 3 3        
\putrule from 3 -1 to 3 0        
\putrule from 4 -2 to 4 3        %
\putrule from 5 -2 to 5 3        %
\putrule from 6 -2 to 6 2        %
\putrule from 7 -2 to 7 -1        %
\putrule from 7 0 to 7 1        %
\endpicture
&\qquad\hbox{and}\qquad 
&\beginpicture
\setcoordinatesystem units <0.5cm,0.5cm>         
\setplotarea x from 0 to 7, y from -2 to 3    
\linethickness=0.5pt                          
\put{1} at 1.5 -0.5
\put{2} at 2.5 -0.5
\put{3} at 0.5 2.5
\put{4}  at 1.5 1.5
\put{6} at 3.5 -0.5
\put{5} at 4.5 -1.5
\put{7}  at 1.5 2.5
\put{8}  at 2.5 1.5
\put{9}  at 4.5 -0.5
\put{10}  at 5.5 -1.5
\put{12}  at 2.5 2.5
\put{13}  at 3.5 1.5
\put{11}  at 4.5 0.5
\put{14}  at 5.5 -0.5
\put{15}  at 6.5 -1.5
\put{16}  at 3.5 2.5
\put{17}  at 4.5 1.5
\put{18}  at 5.5 0.5
\put{19}  at 4.5 2.5
\put{20}  at 5.5 1.5
\put{21}  at 6.5 0.5
\putrule from 0 3  to 5 3          
\putrule from 0 2  to 6 2          %
\putrule from 1 1  to 7 1          %
\putrule from 1 0  to 7 0          %
\putrule from 1 -1 to 7 -1          %
\putrule from 4 -2 to 7 -2          %
\putrule from 0 2 to 0 3        %
\putrule from 1 1 to 1 3        %
\putrule from 1 -1 to 1 0        %
\putrule from 2 1 to 2 3        %
\putrule from 2 -1 to 2 0        %
\putrule from 3 1 to 3 3        
\putrule from 3 -1 to 3 0        
\putrule from 4 -2 to 4 3        %
\putrule from 5 -2 to 5 3        %
\putrule from 6 -2 to 6 2        %
\putrule from 7 -2 to 7 -1        %
\putrule from 7 0 to 7 1        %
\endpicture
\cr
\cr
\cr
\hbox{first ordering of boxes}
&&\hbox{second ordering of boxes} \cr
}$$
and the normal form of $\lambda/\mu$ is 
$$
\eqalign{
\nu &= (7,7,7,6,6,6,5,5,5,5,5,4,4,4,4,3,3,3,3,2,1), \cr
\tilde \nu &= (6,6,6,5,5,5,4,4,4,4,4,3,3,3,3,2,2,2,2,1,0), 
\enspace\hbox{and} \cr
w&=\pmatrix{1 &2 &3 &4 &5 &6 &7 &8 &9 &10 &11 &12 &13 &14 &15 &16 &17 &18 &19 &20 &21 \cr
15 &1 &21 &20 &14 &2 &6 &5 &4 &3 &19 &10 &9 &8 &7 &13 &12 &11 &18 &17 &16 \cr}
\cr}
$$
Let $\fg$ be of type $A_n$ and $V=L(\omega_1)$ and let
$${\cal L}^{\lambda/\mu} = F_\lambda(L(\mu)),\formula$$
as defined in (4.1).  It is known (a consequence of Proposition 6.27 and 
Proposition 4.3c) that ${\cal L}^{\lambda/\mu}$ is always a simple
$\tilde H_k$-module or $0$. Furthermore,
all simple $\tilde H_k$ modules are obtained by this construction.
See [Su] for proofs of these statements.
The following theorem is a reformulation of Proposition 4.12 in terms of the 
combinatorics of our present setting.

\thm  Let $\lambda/\mu$ and $\rho/\tau$ be multisegments with $k$
boxes (with $\mu$ and $\tau$ assumed to be integral) and let
$$\lambda/\mu = \nu/(w\circ\tilde\nu)
\qquad\hbox{and}\qquad
\rho/\tau = \gamma/(v\circ\tilde\gamma)$$
be their normal forms.  Then the
multiplicities of ${\cal L}^{\rho/\tau}$ in a Jantzen filtration of 
${\cal M}^{\lambda/\mu}$ are given by
$$\sum_{j\ge 0} \left[ 
{({\cal M}^{\lambda/\mu})^{(j)}\over
({\cal M}^{\lambda/\mu})^{(j+1)} }:
{\cal L}^{\rho/\tau}\right]
{\tt v}^{{1\over2}(\ell(y)-\ell(w)+j)}
=\cases{
P_{wv}({\tt v}), &if $\nu=\gamma$,\cr
0, &if $\nu\ne \gamma$, \cr
}$$
where $P_{wv}({\tt v})$ is the Kazhdan-Lusztig polynomial
for the symmetric group $S_k$.
\endthm

Theorem 6.31 says that every decomposition number for affine 
Hecke algebra representations is a Kazhdan-Lusztig polynomial.  
The following
is a converse statement which says that every Kazhdan-Lusztig
polynomial for the symmetric group is a decomposition number
for affine Hecke algebra representations.  This statement
is interesting in that Polo [Po] has shown that
{\it every} polynomial in $1+{\tt v}\ZZ_{\ge 0}[{\tt v}]$ is a Kazhdan-Lusztig polynomial
for some choice of $n$ and permutations $v,w\in S_n$.
Thus, the following theorem also shows that every polynomial arises as
a generalized decomposition number for an appropriate
pair of affine Hecke algebra modules.

\prop  
Let $\lambda = (r,r,\ldots,r)=(r^r)$ and 
$\mu = (0,0,\ldots,0)=(0^r)$.
Then, each pair of permutations $v,w\in S_r$,  the Kazhdan-Lusztig polynomial
$P_{vw}({\tt v})$ for the symmetric group $S_r$ is equal to
$$P_{vw}({\tt v}) = 
\sum_{j\ge 0} \left[ 
{({\cal M}^{\lambda/w\circ\mu})^{(j)}\over
({\cal M}^{\lambda/w\circ\mu})^{(j+1)} }:
{\cal L}^{\lambda/v\circ\mu}\right]
{\tt v}^{{1\over2}(\ell(y)-\ell(w)+j)}.$$
\pf
Since $\mu+\rho$ and $\lambda+\rho$ are both
regular, $W_{\lambda+\rho}=W_{\mu+\rho}=1$ and the standard and irreducible
modules ${\cal L}^{\lambda/(w\circ\mu)}$ and ${\cal M}^{\lambda/(v\circ\mu)}$
ranging over all $v,w\in S_k$.  Thus, this statement is a corollary of
Proposition 4.12.  \endpf

\section 7. References

\medskip

\medskip\noindent
\item{[AS]} {\smallcaps T.\ Arakawa and T.\ Suzuki},
{\it Duality between ${\frak{sl}}_n(\CC)$ and the
degenerate affine Hecke algebra of type A}, J.\ Algebra
{\bf 209} (1998), 288-304.

\medskip\noindent
\item{[Ar]} {\smallcaps S.\ Ariki}, 
{\it On the decomposition number of the Hecke algebra of 
$G(m,1,n)$}, J.\ Math.\ Kyoto Univ.\ {\bf 36} (1996), 789-808.

\medskip\noindent
\item{[AK]} {\smallcaps S.\ Ariki and K.\ Koike},
{\it A Hecke algebra of $(\ZZ/r\ZZ)\wr S_n$ an construction of 
its irreducible representations}, Adv.\ in Math.\ {\bf 106} (1994),
216-243.

\medskip\noindent
\item{[Ba]} {\smallcaps P.\ Baumann}, {\it On the center
of quantized enveloping algebras}, J.\ Algebra {\bf 203}
(1998), 244-260.

\medskip\noindent
\item{[Bb]} {\smallcaps D.\ Barbasch},
{\it Filtrations on Verma modules},
Ann.\ Sci.\ \'Ecole Norm.\ Sup.\ (4) {\bf 16} (1983), no.\ 3, (1984) 489--494. 

\medskip\noindent
\item{[BB]} {\smallcaps A.\ Beilinson and J.\ Bernstein},
{\it A proof of Jantzen conjectures},
Adv.\ in Soviet Math.\ {\bf 16} Part 1 (1993), 1-50.

\medskip\noindent
\item{[BGG]} {\smallcaps J.\ Berstein, S.\ Gelfand, I.\ Gelfand},
{\it Differential operators on the base affine space and a study of 
$\fg$-modules}, in {\sl Lie groups and their representations}
(Proc.\ Summer School, Bolya J\'anos Math. Soc., Budapest, 1971, I.M.\ Gelfand
Ed.) London, Hilger, (1975).

\medskip\noindent
\item{[Bou]} {\smallcaps N.\ Bourbaki},
{\sl Groupes et Alg\`ebres de Lie, Chapitres 4, 5 et 6}, Elements de
Math\'ematiques, Hermann, Paris 1968.

\medskip\noindent
\item{[B-tD]} {\smallcaps T.\ Br\"ocker and T.\ tom Dieck},
{\it Representations of compact Lie groups},
Graduate Texts in Mathematics {\bf 98} Springer-Verlag, 
New York-Berlin, 1985. 

\medskip\noindent
\item{[Ch]} {\smallcaps I.\ Cherednik},
{\it A new interpretation of Gelfand-Tzetlin bases}, Duke Math.\ J.\ 
{\bf 54} (1987), 563-577.

\medskip\noindent
\item{[Cr]} {\smallcaps J.\ Crisp}, Ph.D.\ Thesis, Univ.\ of Sydney, 1997.

\medskip\noindent
\item{[CG]} {\smallcaps N.\ Chriss and V.\ Ginzburg},
{\sl Representation theory and complex geometry}, Birkh\"auser,
Boston, 1997.

\medskip\noindent
\item{[CP]} {\smallcaps V.\ Chari and A.\ Pressley},
{\sl A guide to quantum groups},
Cambridge University Press, Cambridge 1994.

\medskip\noindent
\item{[CP2]} {\smallcaps V.\ Chari and A.\ Pressley},
{\it Quantum affine algebras and affine Hecke algebras},
Pacific J.\ Math.\ {\bf 174} (1996), no. 2, 295-326. 

\medskip\noindent
\item{[Dx]} {\smallcaps J.\ Dixmier},
{\sl Enveloping algebras},
 Graduate Studies in Mathematics {\bf 11}, 
American Mathematical Society, Providence, RI, 1996.

\medskip\noindent
\item{[Dr]} {\smallcaps V.\ Drinfeld},
{\it Almost cocommutative Hopf algebras},
Leningrad Math.\ J.\ {\bf 1} (1990), 321-342.

\medskip\noindent
\item{[GJ]} {\smallcaps O.\ Gabber and A.\ Joseph},
{\it On the Bernstien-Gelfand-Gelfand resolution and the
Duflo sum formula}, Compositio Math.\ {\bf 43} (1981), 107-131.

\medskip\noindent
\item{[Gk]} {\smallcaps M.\ Geck}, 
{\it Representations of Hecke algebras at roots of unity},
S\'eminaire Bourbaki, Vol. 1997/98, Ast\'erisque {\bf 252} (1998), 
Exp.\ No.\ 836, 3, 33-55.

\medskip\noindent
\item{[GIM]} {\smallcaps M.\ Geck, L.\ Iancu and G.\ Malle},
{\it Weights of Markov traces and generic degrees},
Indag.\ Math.\ N.S.\ {\bf 11} (2000), 379-397.


\medskip\noindent
\item{[GL]} {\smallcaps M.\ Geck and S.\ Lambropoulou},
{\it Markov traces and knot invariants related to Iwahori-Hecke
algebras of type $B$},
J.\ Reine Angew.\ Math.\ {\bf 482} (1997), 191-213.

\medskip\noindent
\item{[HR]} {\smallcaps T.\ Halverson and A.\ Ram},
{\it Characters of algebras containing a Jones basic construction: 
the Temperley-Lieb, Okada, Brauer, and Birman-Wenzl algebras},
Adv.\ Math.\ {\bf 116} (1995), no. 2, 263-321. 

\medskip\noindent
\item{[H\"a1]} {\smallcaps R.\ H\"aring-Oldenburg},
{\it The reduced Birman-Wenzl algebra of Coxeter type B},
J.\ Algebra {\bf 213} (1999), 437--466. 

\medskip\noindent
\item{[H\"a2]} {\smallcaps R.\ H\"aring-Oldenburg},
{\it An Ariki-Koike like extension of the Birman-Murakami-Wenzl algebra},
preprint 1998.

\medskip\noindent
\item{[Ic]} {\smallcaps L.\ Iancu},
{\it Markov traces and generic degrees in type $B_n$},
J.\ Algebra {\bf 236}, no. 2 (2001), 731-744.

\medskip\noindent
\item{[Jz]} {\smallcaps J.C.\ Jantzen},
{\it Moduln mit einem hochsten Gewicht},
Lecture Notes in Math.\ {\bf 750} (1980) Springer, Berlin.

\medskip\noindent
\item{[Jo1]} {\smallcaps V.\ Jones},
{\it Hecke algebra representations of braid groups and 
link polynomials}, Ann.\ of Math.\ (2) {\bf 126} (1987), 335-388.

\medskip\noindent
\item{[Jo2]} {\smallcaps V.\ Jones},
{\it Index for subfactors}, Invent.\ Math.\ {\bf 72} (1983), 1-25.

\medskip\noindent
\item{[Jo3]} {\smallcaps V.\ Jones},
{\it A quotient of the affine Hecke algebra in the Brauer algebra},
L'Enseignement Math\'ematique {\bf 40} (1994), 313--344.

\medskip\noindent
\item{[KL]} {\smallcaps D.\ Kazhdan and G.\ Lusztig},
{\it Proof of the Deligne-Langlands conjectures for Hecke algebras},
Invent.\ Math.\ {\bf 87} (1987), 153-215.

\medskip\noindent
\item{[Lb]} {\smallcaps S.\ Lambropoulou},
{\it Knot theory related to generalized and cyclotomic 
Hecke algebras of type $B$},
J.\ Knot Theory Ramifications {\bf 8} (1999), 621-658.

\medskip\noindent
\item{[LR]} {\smallcaps R.\ Leduc and A.\ Ram},
{\it A ribbon Hopf algebra approach to the
irreducible representations of centralizer algebras:
The Brauer, Birman-Wenzl and Iwahori-Hecke algebras},
Adv.\ Math.\ {\bf 125} (1997), 1-94.

\medskip\noindent
\item{[Lu]} {\smallcaps G.\ Lusztig},
{\it Affine Hecke algebras and their graded version},
J.\ Amer.\ Math.\ Soc.\ {\bf 2} (1989), 599-635.

\medskip\noindent
\item{[Mac]} {\smallcaps I.G.\ Macdonald},
{\sl Symmetric functions and Hall polynomials},
Second edition, Oxford Mathematical Monographs, Oxford University
Press, New York, 1995.

\medskip\noindent
\item{[Mu]} {\smallcaps J.\ Murakami},
{\it The representations of the $q$-analogue of Brauer's centralizer
algebras and the Kauffman polynomial of links}, Publ.\ RIMS, 
Kyoto Univ.\ {\bf 26} (1990), 935-945.

\medskip\noindent
\item{[Nz]} {\smallcaps M.\ Nazarov},
{\it Young's orthogonal form for Brauer's centralizer algebra},
J.\ Algebra {\bf 182} (1996), no. 3, 664-693. 

\medskip\noindent
\item{[Or]} {\smallcaps R.\ Orellana},
{\it Weights of Markov traces on Hecke algebras},
J reine angew.\ Math.\ {\bf 508} (1999), 157-178.


\medskip\noindent
\item{[Po]} {\smallcaps P.\ Polo}, {\it Construction of 
arbitrary Kazhdan-Lusztig polynomials in symmetric groups},
Representation Theory {\bf 3} (1999), 90-104.

\medskip\noindent
\item{[RR]} {\smallcaps A.\ Ram and J.\ Ramagge},
{\it Affine Hecke algebras, cyclotomic Hecke algebras and Clifford
theory}, preprint 1999.

\medskip\noindent
\item{[Re]} {\smallcaps N.\ Reshetikhin},
{\it Quasitriangle Hopf algebras and invariants of tangles},
Leningrad Math.\ J.\ {\bf 1} (1990), 491-513.

\medskip\noindent
\item{[RT]} {\smallcaps N.\ Reshetikhin and V.\ Turaev},
{\it Invariants of $3$-manifolds via link polynomials and 
quantum groups}, Invent.\ Math.\ {\bf 103} (1991), 547-597.

\medskip\noindent
\item{[SS]} {\smallcaps M.\ Sakamoto and T.\ Shoji},
{\it Schur-Weyl reciprocity for Ariki-Koike algebras},
J.\ Algebra {\bf 221} (1999), no.\ 1, 293-314.

\medskip\noindent
\item{[Su]} {\smallcaps T.\ Suzuki},
{\it Rogawski's conjecture on the Jantzen filtration
for the degenerate affine Hecke algebra of type A},
Representation Theory {\bf 2} (1998), 393-409.

\medskip\noindent
\item{[Tu1]} {\smallcaps V.\ Turaev},
{\it The Yang-Baxter equation and invariants of links},
Invent.\ Math.\ {\bf 92} (1988), 527-553. 

\medskip\noindent
\item{[Tu2]} {\smallcaps V.\ Turaev},
{\it The Conway and Kauffman modules of the solid torus},
J.\ Soviet Math.\ {\bf 52}, no.\ 1 (1990), 2806--2807;
also in {\sl Progress in knot theory and related topics}, 90--102, 
Travaux en Cours {\bf 56}, Hermann, Paris, 1997. 

\medskip\noindent
\item{[TW]} {\smallcaps V.\ Turaev and H.\ Wenzl},
{\it Quantum invariants of $3$-manifolds associated with 
classical simple Lie algebras},
Internat.\ J.\ Math.\ {\bf 4} (1993), 323-358.

\medskip\noindent
\item{[Wz]} {\smallcaps H.\ Wenzl},
{\it Hecke algebras of type $A\sb n$ and subfactors}, 
Invent.\ Math.\ {\bf 92} (1988), 349-383.

\medskip\noindent
\item{[Wz2]} {\smallcaps H.\ Wenzl},
{\it Quantum groups and subfactors of type B, C, and D},
Comm.\ Math.\ Phys.\ {\bf 133} (1990), 383-432.

\medskip\noindent
\item{[Ze]} {\smallcaps A.\ Zelevinsky},
{\it Resolvents, dual pairs and character formulas},
Funct.\ Anal.\ Appl.\ {\bf 21} (1987), 152-154.

\medskip\noindent
\item{[Ze2]} {\smallcaps A.\ Zelevinsky},
{\it Induced representations of reductive p-adic groups II},
Ann.\ Sci.\ \'Ecole Norm.\ Sup.\ (4) Serie {\bf 13} (1980), 
165-210.

\vfill\eject
\end

\subsection The affine Temperley-Lieb algebras ${\widetilde{TL}}_k$

Let $T_0$ be the element of the affine Hecke algebra given by
$$T_0=PICTURE.$$
The {\it affine Temperley-Lieb algebra} ${\widetilde{TL}}_k$ is the
quotient of the affine Hecke algebra by the relations
$$1+qT_i+qT_{i+1}+q^2T_iT_{i+1}+q^2T_{i+1}T_i+q^3T_iT_{i+1}T_i = 0,
\qquad 0\le i\le n-1,$$
where the indices are taken modulo $n$.  Let $E_i = -(T_i+q^{-1})$.
By representing the element $E_i$ as the affine diagram
$$E_i=PICTURE$$
the algebra $TL_k$ can be viewed as the linear span of the affine
planar diagrams (diagrams in the annulus with no string crossings)
with multiplication
$$a_1a_2 = \delta^{\rm (\# of closed loops)} a_1\circ a_2,
\qquad\hbox{where $\delta = -(q+q^{-1})$,}$$
and $a_1\circ a_2$ is the composition of the two affine diagrams $a_1$ and $a_2$
with all loops removed.  

The {\it annular algebra} [Jo] is the quotient of the
affine Temperley-Lieb algebra by the relations
$$\omega^n = 1 \qquad\hbox{and}\qquad \alpha_0=\delta,$$
where $\alpha_0$ is the loop around the center of the annulus.
The {\it Brauer algebra} is the span of rectangular diagrams with multiplication
given by 
$$a_1a_2 = \delta^{\rm (\# of closed loops)} a_1\circ a_2,$$
where $a_1\circ a_2$ is the composition of the two rectangular diagrams 
$a_1$ and $a_2$ with all loops removed.  
The annular algebra can be realized as a subalgebra of the Brauer
algebra by identifying $E_0$ with the diagram
$$E_0=PICTURE\in B_k(\delta).$$
The representation theory of the annular algebra has been studied
in [Jo], [GL1-2], [Grn].